\def\ann{\mathop{\rm ann}\nolimits}
\def\Ann{\mathop{\rm Ann}\nolimits}
\def\ara{\mathop{\rm ara}\nolimits}
\def\Ass{\mathop{\rm Ass}\nolimits}
\def\Att{\mathop{\rm Att}\nolimits}
\def\Assh{\mathop{\rm Assh}\nolimits}
\def\codim{\mathop{\rm codim}\nolimits}
\def\cd{\mathop{\rm cd}\nolimits}
\def\depth{\mathop{\rm depth}\nolimits}
\def\End{\mathop{\rm End}\nolimits}
\def\Ndim{\mathop{\rm N.dim}\nolimits}
\def\width{\mathop{\rm width}\nolimits}
\def\height{\mathop{\rm height}\nolimits}
\def\Hom{\mathop{\rm Hom}\nolimits}
\def\id{\mathop{\rm id}\nolimits}
\def\im{\mathop{\rm im}\nolimits}
\def\dirlim{\mathop{\rm dirlim}\nolimits}
\def\invlim{\mathop{\rm invlim}\nolimits}
\def\ker{\mathop{\rm ker}\nolimits}
\def\length{\mathop{\rm length}\nolimits}
\def\Ext{\mathop{\rm Ext}\nolimits}
\def\InjH{\mathop{\rm E}\nolimits}
\def\LCMo{\mathop{\rm H}\nolimits}
\def\Min{\mathop{\rm Min}\nolimits}
\def\Spec{\mathop{\rm Spec}\nolimits}
\def\Supp{\mathop{\rm Supp}\nolimits}
\def\char{\mathop{\rm char}\nolimits}
\def\Naturalsign{{\rm l\kern-.23em N}}
\baselineskip=15pt
\parindent=0pt
\font \normal=cmr10 scaled \magstep0  \font \gross=cmr10 scaled \magstep5
\input amssym.def
\input amssym.tex
\pageno=1
\def\leaderfill{\leaders\hbox to 1em{\hss.\hss}\hfill}
{\bf Local Cohomology and Matlis duality -- Table of contents}
\bigskip
\bigskip
\line{0 Introduction\leaderfill 2} \line{1 Motivation and General
Results\leaderfill 6} \line{\ \ \ 1.1 Motivation\leaderfill 6}
\line{\ \ \ 1.2 Conjecture (*) on the structure of $\Ass _R(D(\LCMo
^h_{(x_1,\dots ,x_h)R}(R)))$\leaderfill 10} \line{\ \ \ 1.3 Regular
sequences on $D(\LCMo ^h_I(R))$ are well-behaved in some
sense\leaderfill 13} \line{\ \ \ 1.4 Comparison of two Matlis
Duals\leaderfill 15} \line{2 Associated primes -- a constructive
approach\leaderfill 19} \line{3 Associated primes -- the
characteristic-free approach\leaderfill 26} \line{\ \ \ 3.1
Characteristic-free versions of some results\leaderfill 26}\line{\ \
\ 3.2 On the set $\Ass _R(D(\LCMo ^{\dim (R)-1}_I(R)))$\leaderfill
28}\line{4 The regular case and how to reduce to it\leaderfill 35}
\line{\ \ \ 4.1 Reductions to the regular case\leaderfill 35}
\line{\ \ \ 4.2 Results in the general case, i. e. $h$ is
arbitrary\leaderfill 36} \line{\ \ \ 4.3 The case $h=\dim (R)-2$, i.
e. the set $\Ass _R(D(\LCMo ^{n-2}_{(X_1,\dots
,X_{n-2})R}(k[[X_1,\dots ,X_n]])))$\leaderfill 38} \line{5 On the
meaning of a small arithmetic rank of a given ideal\leaderfill 43}
\line{\ \ \ 5.1 An Example\leaderfill 43} \line{\ \ \ 5.2 Criteria
for $\ara (I)\leq 1$ and $\ara (I)\leq 2$\leaderfill 44} \line{\ \ \
5.3 Differences between the local and the graded case\leaderfill
48}\line{6 Applications\leaderfill 50} \line{\ \ \ 6.1
Hartshorne-Lichtenbaum vanishing\leaderfill 50} \line{\ \ \ 6.2
Generalization of an example of Hartshorne\leaderfill 52} \line{\ \
\ 6.3 A necessary condition for set-theoretic complete
intersections\leaderfill 54} \line{\ \ \ 6.4 A generalization of
local duality\leaderfill 55} \line{7 Further Topics\leaderfill 57}
\line{\ \ \ 7.1 Local Cohomology of formal schemes\leaderfill 57}
\line{\ \ \ 7.2 $D(\LCMo ^i_I(R))$ has a natural $D$-module
structure\leaderfill 57} \line{\ \ \ 7.3 The zeroth Bass number of
$D(\LCMo ^i_I(R))$ (w. r. t. the zero ideal) is not finite in
general\leaderfill 58} \line{\ \ \ 7.4 On the module $\LCMo
^h_I(D(\LCMo ^h_I(R)))$\leaderfill 61} \line{8 Attached prime ideals
and local homology\leaderfill 65} \line{\ \ \ 8.1 Attached prime
ideals -- basics\leaderfill 65} \line{\ \ \ 8.2 Attached prime
ideals -- results\leaderfill 68} \line{\ \ \ 8.3 Local homology and
a necessary condition for Cohen-Macaulayness\leaderfill 71} \line{\
\ \ 8.4 Local homology and Cohen-Macaulayfications\leaderfill
74}\line{References\leaderfill 76}\line{Summary in German (deutsche
Zusammenfassung)\leaderfill 78}
 \vfil \eject \gross
Local Cohomology and Matlis duality \normal
\bigskip
\bigskip
{\bf 0 Introduction}
\bigskip
\bigskip
{\parindent=15pt In algebraic geometry, a (set-theoretic) complete
intersection is a variety $Y$ (say, in affine or projective space
over a field) that can be cut out be $\codim (Y)$ many equations.
For example, every curve in affine $n$-space over a field of
positive characteristic is a set-theoretic complete intersection
(see [CN]). On the other hand, many questions on complete
intersections are still open: Is every closed point in ${\bf
P}^2_{\bf Q}$ (projective 2-space over $\bf Q$, the rationals) a
set-theoretic complete intersection? Is every irreducible curve in
${\bf A}^3_{\bf C}$ (affine 3-space over $\bf C$, the complex
numbers) a set-theoretic complete intersection? See [Ly2] for a
survey on these and other questions.
\par
Here is another example: Over an algebraically closed field $k$, let
$C_d\subseteq {\bf P}_k^3$ be the curve parameterized by
$$(u^d:u^{d-1}v:uv^{d-1}:v^d)$$
(for $(u:d)\in {\bf P}_k^1$). Hartshorne has shown (see [Ha2,
Theorem.*]) that, in positive characteristic, every curve $C_d$ is a
set-theoretic complete intersection. In characteristic zero, the
question is open. It is even unknown if $C_4$ is a set-theoretic
complete intersection or not. An obvious obstruction for $C_4$ to be
a set-theoretic complete intersection would be $\LCMo ^3_I(R)\neq 0$
($I\subseteq R=k[X_0,X_1,X_2,X_3]$ the vanishing ideal of
$C_4\subseteq {\bf P}_k^3$), but, as is well-known, one has
$$\LCMo ^3_I(R)=0\ \ .$$
Thus, if we define the so-called arithmetic rank of $I$,
$$\ara(I):=\min \{ l\in \Naturalsign \vert \exists r_1,\dots ,r_l\in
R:\sqrt {I}=\sqrt {(r_1,\dots ,r_l)R}\} \ \ ,$$ it seems that
(non-)vanishing of the modules $\LCMo ^i_I(R)$ does not carry enough
information to determine $\ara (I)$ (because our example 5.1 shows
that this can really happen in the sense that $\cd (I)<\ara (I)$,
where $\cd $ is the (local) cohomological dimension of $I$).
\par
It is interesting that, although the vanishing of $\LCMo ^3_I(R)$
does not seem to help in the case of $C_4$, the Matlis dual $D(\LCMo
^2_I(R))$ (note that $D$ will stand for the Matlis dual functor,
also see the end of this introduction for more notation) of the
module $\LCMo ^2_I(R)$ "knows" whether we have a set-theoretic
complete intersection or not, in the following sense (take $h=2$):}
\bigskip {\bf (1.1.4 Corollary)}
\par
Let $(R,\goth m)$ be a noetherian local ring, $I$ a proper ideal of
$R$, $h\in \Naturalsign $ and $\underline f=f_1,\dots ,f_h\in I$ an
$R$-regular sequence. The following statements are equivalent:
\par
(i) $\sqrt {\underline fR}=\sqrt I$, i. e. $I$ is -- up to radical
-- the set-theoretic complete intersection ideal $\underline fR$; in
particular, it is a set-theoretic complete intersection ideal
itself.
\par
(ii) $\LCMo ^l_I(R)=0$ for every $l>h$ and the sequence $\underline
f$ is quasi-regular on $D(\LCMo ^h_I(R))$.
\par
(iii) $\LCMo ^l_I(R)=0$ for every $l>h$ and the sequence $\underline
f$ is regular on $D(\LCMo ^h_I(R))$. \bigskip {\parindent=15pt This
result gives motivation to study modules of the form $D(\LCMo
^l_I(R))$, in particular its associated primes (as they determine
which elements operate injectively on $D(\LCMo ^l_I(R))$). Modules
of the form $D(\LCMo ^l_I(R))$ and their associated prime ideals
have been studied in [H2], [H3], [H4], [H5], [HS1] and [HS2].}
\bigskip
\bigskip
{\bf Main results}
\bigskip
In the sequel we will list the main results of this work. In this
context we would like to remark that conjecture (*) (1.2.2) is a
central theme of this work. We also remark that many of the results
listed below (e. g. (1.2.1), (2.7) (i), (3.1.3) (ii) and (iii),
(3.2.7), (4.1.2), but also (8.2.6) (iii) ($\zeta $)) give evidence
for it, though we are not able to prove conjecture (*).
\bigskip
We also note that the results in this work lead to various
applications. These applications are collected in section 6, they
are not listed here.
\bigskip
{\bf (1.2.1 Remark)}
\par
Let $(R,\goth m)$ be a noetherian local ring and $\underline
x=x_1,\dots ,x_h$ a sequence in $R$. Then one has
$$\Ass _R(D(\LCMo ^h_{(x_1,\dots ,x_h)R}(R)))\subseteq \{ \goth p\in
\Spec (R)\vert \LCMo ^h_{(x_1,\dots ,x_h)R}(R/\goth p)\neq 0\} \ \
.$$
\bigskip
Though easy, remark 1.2.1 is crucial for many proofs in this work;
it seems reasonable to conjecture:
\bigskip {\bf (1.2.2
Conjecture)}
\par
If $(R,\goth m)$ is a noetherian local ring, $h>0$ and $x_1,\dots
,x_h$ are elements of $R$,
$$\Ass _R(D(\LCMo ^h_{(x_1,\dots ,x_h)R}(R)))=\{ \goth p\in
\Spec (R)\vert \LCMo ^h_{(x_1,\dots ,x_h)R}(R/\goth p)\neq 0\} \ \
.\leqno{(*)} $$ Besides remark 1.2.1, there is more evidence for
conjecture (*), e. g.:
\bigskip
{\bf (3.1.3 Theorem, statements (ii) and (iii))}
\par
Let $(R,\goth m)$ be a noetherian local ring, $\underline
x=x_1,\dots ,x_m$ a sequence in $R$ and $M$ a finitely generated
$R$-module. Then $$\{ \goth p\in \Supp _R(M)\vert x_1,\dots
,x_m\hbox { is part of a system of parameters of }R/\goth p\}
\subseteq \Ass _R(D(\LCMo ^m_{(x_1,\dots ,x_m)R}(M)))$$ holds. Now,
if we assume furthermore that $R$ is a domain and $\underline x$ is
part of a system of parameters of $R$, we have $\{ 0\} \in \Ass
_R(D:=D(\LCMo ^m_{\underline xR}(R)))$. Therefore, it is natural to
ask for the zeroth Bass number of $D$ with respect to the zero
ideal. We will see that, in general, this number is not finite
(theorem 7.3.2). In the special case $m=1$ we can completely compute
the associated prime ideals: Namely, for every $x\in R$, one has
$$\Ass _R(D(\LCMo ^1_{xR}(R)))=\Spec (R)\setminus {\cal V}(x)$$
(${\cal V}(x)$ is the set of all prime ideals of $R$ containing
$x$). In particular, the set
$$\Ass _R(D(\LCMo ^m_{(x_1,\dots ,x_m)R}(R)))$$
is, in general, infinite. Here is further evidence for conjecture
(*):
\bigskip
{\bf (3.2.7 Theorem)}
\par
Let $(R,\goth m)$ be a $d$-dimensional local complete ring and
$J\subseteq R$ an ideal such that $\dim (R/J)=1$ and $\LCMo
^d_J(R)=0$. Then
$$\Ass _R(D(\LCMo ^{d-1}_J(R))=\{ P\in \Spec(R)\vert \dim (R/P)=d-1,\dim
(R/(P+J))=0\} \cup \Assh (R)$$ holds. Here $\Assh (R)$ denotes the
set of all associated prime ideals of $R$ of highest dimension.
Further evidence for (*) can be found in section 8.2 in connection
with attached primes (see 8.2.6 (iii) ($\zeta $) for details).
{\parindent=15pt \par We continue our list of main results on Matlis
duals of local cohomology modules:} \bigskip {\bf (4.1.2 Theorem)}
\par
Let $(R,\goth m)$ be a noetherian local complete ring with
coefficient field $k\subseteq R$, $l\in \Naturalsign ^+$ and
$x_1,\dots ,x_l\in R$ a part of a system of parameters of $R$. Set
$I:=(x_1,\dots ,x_l)R$. Let $x_{l+1},\dots ,x_d\in R$ be such that
$x_1,\dots ,x_d$ is a system of parameters of $R$. Denote by $R_0$
the (regular) subring $k[[x_1,\dots ,x_d]]$ of $R$. Then if $\Ass
_{R_0}(D(\LCMo ^l_{(x_1,\dots ,x_l)R_0}(R_0)))$ is stable under
generalization, $\Ass _R(D(\LCMo ^l_I(R)))$ is also stable under
generalization. \bigskip (A set $X$ of prime ideals of a ring is
stable under generalization, if $\goth p\in X$ implies $\goth p_0\in
X$ for every $\goth p_0\subseteq \goth p$.) Clearly, 4.1.2 can be
helpful when we want to reduce from a general (complete) to a
regular (complete) case. {\parindent=15 pt
\par
The next result shows that the question when $\LCMo ^{\dim
(R)-1}_I(R)$ is zero (for an ideal $I$ in a local regular ring $R$)
is related to the question which prime ideals are associated to the
Matlis dual of a certain local cohomology module:}
\bigskip
{\bf (4.3.1 Corollary)}
\par
Let $R_0$ be a noetherian local complete equicharacteristic ring,
let $\dim (R_0)=n-1$, $k\subseteq R_0$ a coefficient field of $R_0$.
Let $x_1,\dots ,x_n$ be elements of $R_0$ such that $\sqrt
{(x_1,\dots ,x_n)R_0}=\goth m_0$. Set $I_0:=(x_1,\dots
,x_{n-2})R_0$. Let $R:=k[[X_1,\dots ,X_n]]$ be a power series
algebra over $k$ in the variables $X_1,\dots ,X_n$, $I:=(X_1,\dots
,X_{n-2})R$. Then the $k$-algebra homomorphism $R\to R_0$ determined
by $X_i\mapsto x_i$ $(i=1,\dots n$) induces a module-finite ring map
$\iota :R/fR\to R_0$ for some prime element $f\in R$. Furthermore,
suppose that $R_0$ is regular and $\height (I_0)<h$; then we have
$$fR\in \Ass _R(D(\LCMo ^h_I(R)))\iff \LCMo ^{n-2}_{I_0}(R_0)\neq 0\ \
.$$ In this case, $fR$ is a maximal element of $\Ass _R(D(\LCMo
^h_I(R)))$. By [HL, Theorem 2.9] the latter holds if and only if
$\dim (R_0/I_0)\geq 2$ and $\Spec (\overline {R_0}/I_0\overline
{R_0})\setminus \{ \goth m_0(\overline {R_0}/I_0\overline {R_0})\} $
is connected, where $\overline {R_0}$ is defined as the completion
of the strict henselization of $R_0$; this means that $\overline
{R_0}$ is obtained from $R_0$ by replacing the coefficient field $k$
by its separable closure in any fixed algebraic closure of $k$.
{\parindent=15 pt
\bigskip
It was shown in [Ly1, Example 2.1. (iv)] that every local cohomology
module $\LCMo ^i_I(R)$ has a natural $D$-module structure, where
$$D:=D(R,k)\subseteq \End _k(R)$$ is the subring generated by all
$k$-linear derivations (from $R$ to $R$) and the multiplications by
elements of $R$ (here $k\subseteq R$ is any subring). We show in
section 7.2 that, at least if $R=k[[X_1,\dots ,X_n]]$ is a formal
power series ring over $k$, the Matlis dual
$$D(\LCMo ^i_I(R))$$
has a canonical $D$-module structure, too (for every ideal
$I\subseteq R$); furthermore, we will see that, with respect to this
$D$-module structure, $D(\LCMo ^i_I(R))$ is not finitely generated,
in general; in particular, it is not holonomic (see [Bj, in
particular sections 1 and 3] for the notion of holonomic
$D$-modules).
\par
We will use the $D$-module structure on $D(\LCMo ^i_I(R))$ to show
}\bigskip {\bf (7.4.1 and 7.4.2 Theorems (special cases))}
\par
Let $(R,\goth m)$ be a noetherian local complete regular ring of
equicharacteristic zero, $I\subseteq R$ an ideal of height $h\geq 1$
such that $\LCMo ^l_I(R)=0$ for every $l>h$, and $\underline
x=x_1,\dots ,x_h$ an $R$-regular sequence in $I$; then, in general,
$$\LCMo ^h_I(D(\LCMo ^h_I(R)))$$
is either $\InjH _R(R/\goth m)$ or zero; if we assume
$$I=(x_1,\dots ,x_h)R$$
in addition, we have
$$\LCMo ^h_I(D(\LCMo ^h_I(R)))=\InjH _R(R/\goth m)\ \ .$$
{\parindent=15pt
\par
Further main results of this work are contained in section 6, in
which we collect various applications of our theory, namely new
proofs for Hartshorne-Lichtenbaum vanishing (6.1), a generalization
of an example of a non-artinian but zero-dimensional local
cohomology module (the original example, which is more special, is
from Hartshorne) (6.2), a new necessary condition for an ideal to be
a set-theoretic complete intersection ideal (6.3) and a
generalization of local duality (6.4).
\bigskip
}{\bf Notation} \bigskip {\parindent=15pt If $I$ is an ideal of a
ring $R$ and $M$ is an $R$-module, we denote by $\LCMo ^l_I(M)$ the
$l$-th local cohomology of $M$ supported in $I$; material on local
cohomology can e. g. be found in [BH], [BS], [Gr] and [Hu]. If
$(R,\goth m)$ is a noetherian local ring, $\InjH _R(R/\goth m)$
stands for any (fixed) $R$-injective hull of the $R$-module $R/\goth
m$; see, for example, [BH] and [Ms] for more on injective modules.
Finally, $D$ is the Matlis dual functor from the category of
$R$-modules to itself, i. e.
$$D(M):=\Hom _R(M,\InjH _R(R/\goth m))$$
for every $R$-module $M$. The term "Matlis dual of $M$" will always
stand for $D(M)$ (and therefore, will only be used over a local ring
$(R,\goth m)$). Sometimes we will write $D_R$ instead of $D$ to
avoid misunderstandings. References for general facts from
commutative algebra are [Ei], [Ma].}
\bigskip
\bigskip
{\bf Acknowledgement}
\bigskip
I thank J\"urgen St\"uckrad and Gennady Lyubeznik for many helpful
discussions. \vfil \eject {\bf 1 Motivation and General Results}
\bigskip
\bigskip
{\bf 1.1 Motivation} \bigskip Let $I$ be an ideal of a noetherian
ring $R$.
$$\ara (I):=\min \{ l\in \Naturalsign \vert \exists
x_1,\dots ,x_l\in I:\sqrt I=\sqrt {(x_1,\dots ,x_l)R}\} $$ denotes
the arithmetic rank of $I$. Geometrically, it is the (minimal)
number of equations needed to cut out a given algebraic set (say, in
an affine space). It is well-known (and follows by using
\v{C}ech-cohomology) that one has
$$\LCMo ^l_I(R)=0\ \ (l>\ara (I))\ \ .$$
But, conversely, it is in general not true that $\ara (I)$ is
determined by these vanishing conditions, see Example 5.1 for a
counterexample. Assume that $I$ is generated up to radical by a
regular sequence $\underline f=f_1,\dots ,f_h$ in $R$. Then
$\underline f$ is also a regular sequence on $D(\LCMo ^h_I(R))$
(this follows from theorem 1.1.2 resp. corollary 1.1.4 below, see
definition 1.1.1 below for a definition of regular sequences in this
context). It is an interesting fact that the reversed statement also
holds: If $\underline f$ is a $D(\LCMo ^h_I(R))$-regular sequence
then
$$\sqrt I=\sqrt {\underline f}$$
holds. This fact is one of the main motivations for the study of
Matlis duals of local cohomology modules (see theorem 1.1.3 resp.
corollary 1.1.4 for details and the precise statement).
\bigskip {\bf 1.1.1 Definition}
\par Let $R$ be a ring, $M$ an $R$-module, $h\in \Naturalsign $ and
$\underline f=f_1,\dots ,f_h$ a sequence of elements of $R$.
$\underline f$ is called a quasi-regular sequence on $M$ if
multiplication by $f_i$ is injective on $M/(f_1,\dots ,f_{i-1})M$
for every $i=1,\dots ,h$. $\underline f$ is called a regular
sequence on $M$ if $\underline f$ is quasi-regular on $M$ and
$M/\underline fM\neq 0$ holds, in addition.
\bigskip
Before we show the statements on regular sequences mentioned in the
introduction of this section (corollary 1.1.4), we prove something
slightly more general (namely theorems 1.1.2 and 1.1.3); corollary
1.1.4 then simply combines the most interesting special cases from
these two theorems.
\bigskip
{\bf 1.1.2 Theorem}
\par
Let $(R,\goth m)$ be a noetherian local ring, $I$ an ideal of $R$,
$h\geq 1$ and $\underline f=f_1,\dots ,f_h\in I$ a sequence of
elements such that $\sqrt {\underline fR}=\sqrt I$ and
$$\LCMo ^{h-1-l}_I(R/(f_1,\dots ,f_l)R)=0\ \ (l=0,\dots ,h-3)$$
hold (of course, for $h\leq 2$, this condition is void). Then
$\underline f$ is a quasi-regular sequence on $D(\LCMo ^h_I(R))$.
\par
Proof:
\par
By induction on $h$: $h=1$: the functor $\LCMo ^1_I$ is right-exact
because $\LCMo ^2_I=\LCMo ^2_{f_1R}=0$. Hence the exact sequence
$$R\buildrel f_1\over \to R\to R/f_1R\to 0$$
induces an exact sequence
$$\LCMo ^1_I(R)\buildrel f_1\over \to \LCMo ^1_I(R)\to \LCMo
^1_I(R/f_1R)=\LCMo ^1_{f_1R}(R/f_1R)=0$$ (here $f_1$ stands for
multiplication by $f_1$ on $R$ resp. $\LCMo ^1_I(R)$). Application
of $D$ to the last sequence yields injectivity of $f_1$ on $D(\LCMo
^1_I(R))$.
\par
$h=2$: We have $\LCMo ^2_I(R/f_1R)=\LCMo ^2_{(f_1,f_2)R}(R/f_1R)=0$
This implies both $\LCMo ^l_I(M)=0$ for every $l\geq 2$ and every
$R/f_1R$-module $M$ and the fact that $f_1$ operates injectively on
$D(\LCMo ^2_I(R))$. Now the exact sequence
$$0\to (0:_Rf_1)\to R\buildrel \alpha \over \to f_1R\to 0$$
(where $\alpha $ is induced by multiplication by $f_1$) induces an
exact sequence
$$\LCMo ^2_I((0:_Rf_1))\to \LCMo ^2_I(R)\buildrel \LCMo ^2_I(\alpha
)\over \to \LCMo ^2_I(f_1R)\to 0\ .$$ But $(0:_Rf_1R)$ is an
$R/f_1R$-module and so $\LCMo ^2_I((0:_Rf_1R))=0$, showing that
$\LCMo ^2_I(\alpha )$ is an isomorphism. On the other hand the exact
sequence
$$0\to f_1R\buildrel \beta \over \to R\to R/f_1R\to
0$$ (where $\beta $ is an inclusion map) induces an exact sequence
$$\LCMo ^1_I(R/f_1R)\to \LCMo ^2_I(f_1R)\buildrel \LCMo ^2_I(\beta
)\over \to \LCMo ^2_I(R)\to 0\ ,$$ which shows the existence of a
natural epimorphism
$$\eqalign {\LCMo ^1_I(R/f_1R)&\to \ker (\LCMo ^2_I(\beta ))\cr &\cong \ker (\LCMo
^2_I(\beta )\circ \LCMo ^2_I(\alpha ))\cr &=\ker (\LCMo ^2_I(\beta
\circ \alpha ))\cr &=\ker (f_1)\ \ ,\cr }$$ where $f_1$ denotes
multiplication by $f_1$ on $\LCMo ^2_I(R)$. This means that we have
a surjection
$$\LCMo ^1_I(R/f_1R)\to \Hom _R(R/f_1R,\LCMo ^2_I(R))$$
and hence an injection
$$D(\LCMo ^2_I(R))/f_1D(\LCMo ^2_I(R))=D(\Hom _R(R/f_1R,\LCMo
^2_I(R)))\to D(\LCMo ^1_I(R/f_1R))\ \ .$$ Note that the first
equality follows formally from the exactness of $D$; note also that
it does not make any difference if one takes the last Matlis dual
with respect to $R$ or $R/f_1R$. For this reason the case $h=1$
shows that $f_2$ operates injectively on $D(\LCMo ^1_I(R/f_1R))$ and
thus also on $D(\LCMo ^2_I(R))/f_1D(\LCMo ^2_I(R))$. \par Now we
consider the general case $h\geq 3$: Similar to the case $h=2$ we
see that $\LCMo ^l_I(M)=0$ for every $l\geq h$ and every
$R/f_1R$-module $M$ and that $f_1$ operates injectively on $D(\LCMo
^h_I(R))$. The short exact sequence
$$0\to (0:_Rf_1)\to R\buildrel \alpha \over \to f_1R\to 0$$
(where, again, $\alpha $ is induced by multiplication by $f_1$ on
$R$) induces an exact sequence
$$\LCMo ^h_I((0:_Rf_1))\to \LCMo ^h_I(R)\buildrel \LCMo ^h_I(\alpha
)\over \to \LCMo ^h_I(f_1R)\to 0\ .$$ But $(0:_Rf_1)$ is an
$R/f_1R$-module and therefore $\LCMo ^h_I((0:_Rf_1))=0$, showing
that $\LCMo ^h_I(\alpha )$ is an isomorphism. On the other hand the
short exact sequence
$$0\to f_1R\buildrel \beta \over \to R\to R/f_1R\to
0$$ (where $\beta $ is an inclusion map) induces an exact sequence
$$0=\LCMo ^{h-1}_I(R)\to \LCMo ^{h-1}_I(R/f_1R)\to \LCMo
^h_I(f_1R)\buildrel \LCMo ^h_I(\beta )\over \to \LCMo ^h_I(R)\to 0$$
(here we use the fact that $h\geq 3$ and therefore $\LCMo
^{h-1}_I(R)=0$). We conclude
$$\LCMo ^{h-1}_I(R/f_1R)=\ker (\LCMo ^h_I(\beta ))\cong \ker (\LCMo
^h_I(\beta )\circ \LCMo ^h_I(\alpha ))=\Hom _R(R/f_1R,\LCMo
^h_I(R))$$ and, by Matlis duality,
$$D(\LCMo ^{h-1}_I(R/f_1R))=D(\Hom _R(R/f_1R,\LCMo ^h_I(R)))=D(\LCMo
^h_I(R))/f_1D(\LCMo ^h_I(R))\ .$$ Because of our hypothesis, we can
apply the induction hypothesis (to the ring $R/f_1R$) which says
that $f_2,\dots ,f_h$ is a quasi-regular sequence on $D(\LCMo
^{h-1}_I(R/f_1R))$; thus, by the last formula, $\underline f$ is a
quasi-regular sequence on $D(\LCMo ^h_I(R))$.
\bigskip
{\bf 1.1.3 Theorem}
\par
Let $I$ be an ideal of a noetherian local ring $(R,\goth m)$, $h\geq
1$ and $f_1,\dots ,f_h\in I$ be such that
$$\LCMo ^l_I(R)=0\ \ (l>h)$$
and
$$\LCMo ^{h-1-l}_I(R/(f_1,\dots ,f_l)R)=0\ \ (l=0,\dots ,h-2)$$
hold (of course, for $h<2$, this condition is void) and such that
the sequence $\underline f=f_1,\dots ,f_h$ is quasi-regular on
$D(\LCMo ^h_I(R))$. Then $\sqrt I=\sqrt {(f_1,\dots ,f_h)R}$ holds.
\par
Proof:
\par
By induction on $h$: $h=1$: By our hypothesis, the functor $\LCMo
^1_I$ is right-exact. Therefore the exact sequence
$$R\buildrel f_1\over \to R\to R/f_1R\to 0$$
induces an exact sequence
$$\LCMo ^1_I(R)\buildrel f_1\over \to \LCMo ^1_I(R)\to \LCMo
^1_I(R/f_1R)\to 0\ \ ,$$ where $f_1$ stands for multiplication by
$f_1$ on $R$ resp. on $\LCMo ^1_I(R)$. But multiplication by $f_1$
is injective on $D(\LCMo ^1_I(R))$ and so we get $\LCMo
^1_I(R/f_1R)=0$; by our hypothesis, we have $\LCMo ^l_I(R/f_1R)=0$
for every $l\geq 1$. It is well-known that the latter condition is
equivalent to $\LCMo ^l_I(R/\goth p)=0$ for every $l\geq 1$ and
every prime ideal $\goth p$ of
$${\cal V}(f_1R):=\{ \goth p\in \Spec
(R)\vert f_1R\subseteq \goth p\} \ \ .$$
Thus, we must have
$I\subseteq \sqrt {f_1R}$ and, therefore, $\sqrt I=\sqrt {f_1R}$.
\par
$h\geq 2$: Similar to the case $h=1$ we see that $\LCMo
^h_I(R/f_1R)=0$ holds. By our hypothesis, we get $\LCMo ^l_I(M)=0$
for every $l\geq h$ and every $R/f_1R$-module $M$. The short exact
sequence
$$0\to (0:_Rf_1)\to R\buildrel \alpha \over \to f_1R\to 0$$
induces an exact sequence
$$\LCMo ^h_I((0:_Rf_1))\to \LCMo ^h_I(R)\buildrel \LCMo
^h_I(\alpha )\over \to \LCMo ^h_I(f_1R)\to 0\ \ .$$ But $\LCMo
^h_I((0:_Rf_1))=0$, because $(0:_Rf_1)$ is an $R/f_1R$-module. Thus
$\LCMo ^h_I(\alpha )$ is an isomorphism. Now the short exact
sequence
$$0\to f_1R\buildrel \beta \over \to R\to R/f_1R\to
0$$ (where $\beta $ is an inclusion map) induces an exact sequence
$$0=\LCMo ^{h-1}_I(R)\to \LCMo ^{h-1}_I(R/f_1R)\to \LCMo
^h_I(f_1R)\buildrel \LCMo ^h_I(\beta )\over \to \LCMo ^h_I(R)\to 0\
\ ,$$ from which we conclude
$$\LCMo ^{h-1}_I(R/f_1R)=\ker (\LCMo ^h_I(\beta ))\cong \ker
(\LCMo ^h_I(\beta )\circ \LCMo ^h_I(\alpha ))=\Hom _R(R/f_1R,\LCMo
^h_I(R))\ \ .$$ Here we used the facts the $\LCMo ^h_I(\alpha )$ is
an isomorphism and that $\beta \circ \alpha $ is multiplication by
$f_1$ on $R$. Application of the functor $D$ shows
$$D(\LCMo ^{h-1}_I(R/f_1R))=D(\LCMo ^h_I(R))/f_1D(\LCMo ^h_I(R))\ \
.$$ Note that, again, it is irrelevant whether we take the first
functor $D$ here with respect to $R$ or with respect to $R/f_1R$ and
so our induction hypothesis (applied to $R/f_1R$) implies that
$f_2,\dots ,f_h$ is a quasi-regular sequence on $D(\LCMo
^{h-1}_I(R/f_1R))$ and that
$$\sqrt {I(R/f_1R)}=\sqrt {(f_2,\dots ,f_h)\cdot (R/f_1R)}$$
holds. The statement $\sqrt I=\sqrt {(f_1,\dots ,f_h)R}$ follows.
\bigskip
Now it is easy to specialize to the following statement:
\bigskip
{\bf 1.1.4 Corollary}
\par
Let $(R,\goth m)$ be a noetherian local ring, $I$ a proper ideal of
$R$, $h\in \Naturalsign $ and $\underline f=f_1,\dots ,f_h\in I$ an
$R$-regular sequence. The following statements are equivalent:
\par
(i) $\sqrt {\underline fR}=\sqrt I$.
\par
(ii) $\LCMo ^l_I(R)=0$ for every $l>h$ and the sequence $\underline
f$ is quasi-regular on $D(\LCMo ^h_I(R))$.
\par
(iii) $\LCMo ^l_I(R)=0$ for every $l>h$ and the sequence $\underline
f$ is regular on $D(\LCMo ^h_I(R))$.
\par
(The case $h=0$ means
$$\eqalign {\sqrt I=\sqrt 0&\iff \LCMo ^l_I(R)=0 \hbox{ for every }l>0\cr &\iff \LCMo ^l_I(R)=0 \hbox{ for every }l>0 \hbox { and } \Gamma _I(R)\neq 0\ \ ).\cr }$$
Proof:
\par
$h=0$: Clearly the condition $\sqrt I=\sqrt 0$ implies $\LCMo
^l_I(R)=0$ for every $l>0$. On the other hand, if we have $\LCMo
^l_I(R)=0$ for every $l>0$, then, by a well-known theorem, one also
has $\LCMo ^l_I(R/\goth p)=0$ for every prime ideal $\goth p$ of $R$
and for every $l>0$; thus $I\subseteq \goth p$ for every prime ideal
$\goth p$ of $R$. But then it is also true that $\Gamma _I(R)=R\neq
0$ holds.
\par
$h\geq 1:$ The fact that (i) and (ii) are equivalent follows from
theorems 1.1.2 and 1.1.3. Thus we only have to show that (i) implies
$$D(\LCMo ^h_I(R))/(f_1,\dots ,f_h)D(\LCMo ^h_I(R)))\neq 0\ \ :$$But, by
general Matlis duality theory, the last module is
$$D(\Hom _R(R/(f_1,\dots ,f_h)R,\LCMo ^h_I(R)))\ \ ;$$
furthermore, every element
of $\LCMo ^h_I(R)$ is annihilated by a power of $I$ and so it
suffices to show $\LCMo ^h_I(R)\neq 0$, which is clear, because $I$
is generated up to radical by the regular sequence $f_1,\dots ,f_h$.
\bigskip
\bigskip
{\bf 1.2 Conjecture (*) on the structure of $\Ass _R(D(\LCMo
^h_{(x_1,\dots ,x_h)R}(R)))$}
\bigskip
Now, we present an easy property of associated primes of Matlis
duals of certain local cohomology modules; this property will
naturally lead us to a conjecture on the structure of the set of
associated prime ideals of such modules.
\bigskip
{\bf 1.2.1 Remark}
\par
Let $(R,\goth m)$ be a noetherian local ring, $M$ an $R$-module,
$h\in \Naturalsign $ and $I\subseteq R$ an ideal such that $\LCMo
^l_I(M)=0$ holds for every $l>h$ and suppose that we have
$$\goth p\in \Ass _R(D(\LCMo
^h_I(M)))\ \ :$$
This condition clearly implies $\Ann _R(M)\subseteq
\goth p$ (because $\Ann _R(M)$ is contained in the annihilator of
every element of $D(\LCMo ^h_I(M))$) and
$$\eqalign {0&\neq \Hom _R(R/\goth p,D(\LCMo ^h_I(M)))\cr &=D(\LCMo
^h_I(M)\otimes _R(R/\goth p))\cr &=D(\LCMo ^h_I(M/\goth pM))\ \ ,\cr
}$$ i. e. $\LCMo ^h_I(M/\goth pM)\neq 0$. In particular $\dim (\Supp
_R(M/\goth pM))\geq h$.\smallskip As a special case we get the
implication
$$\goth p\in \Ass _R(D(\LCMo ^h_{(x_1,\dots ,x_h)R}(R)))\Rightarrow
\dim (R/\goth p)\geq h$$ for every sequence $x_1,\dots ,x_h\in R$.
\smallskip
Furthermore, as we have seen,
$$\Ass _R(D(\LCMo ^h_{(x_1,\dots ,x_h)R}(R)))\subseteq \{ \goth p\in
\Spec (R)\vert \LCMo ^h_{(x_1,\dots ,x_h)R}(R/\goth p)\neq 0\} $$
holds for every sequence $x_1,\dots ,x_h\in R$. \bigskip {\bf 1.2.2
Conjecture}
\par
If $(R,\goth m)$ is a noetherian local ring, $h>0$ and $x_1,\dots
,x_h$ are elements of $R$,
$$\Ass _R(D(\LCMo ^h_{(x_1,\dots ,x_h)R}(R)))=\{ \goth p\in
\Spec (R)\vert \LCMo ^h_{(x_1,\dots ,x_h)R}(R/\goth p)\neq 0\}
\leqno{(*)} $$ holds. We denote this conjecture by (*). It is one of
the central themes of this work. The next theorem 1.2.3 presents
some equivalent characterizations of conjecture (*); one of them is
stableness under generalization of the set of associated primes of
the Matlis dual of the local cohomology module in question
(condition (ii) from theorem 1.2.3). The theorem also shows
(condition (iv)) that (*) is actually equivalent to a similar
statement, where $R$ is replaced by a finite $R$-module $M$, i. e.
if (*) holds, a version of (*) also holds for finite $R$-modules.
\bigskip
{\bf 1.2.3 Theorem}
\par
The following statements are equivalent:
\par
(i) Conjecture (*) holds, i. e. for every noetherian local ring
$(R,\goth m)$, every $h>0$ and every sequence $x_1,\dots ,x_h$ of
elements of $R$ the equality
$$\Ass _R(D(\LCMo ^h_{(x_1,\dots ,x_h)R}(R)))=\{ \goth p\in
\Spec (R)\vert \LCMo ^h_{(x_1,\dots ,x_h)R}(R/\goth p)\neq 0\} $$
holds.
\par
(ii) For every noetherian local ring  $(R,\goth m)$, every $h>0$ and
every sequence $x_1,\dots ,x_h$ of elements of $R$ the set
$$Y:=\Ass _R(D(\LCMo ^h_{(x_1,\dots ,x_h)}(R)))$$
is stable under generalization, i. e. the implication
$$\goth p_0,\goth p_1\in \Spec (R),\goth p_0\subseteq \goth p_1,
\goth p_1\in Y\Longrightarrow \goth p_0\in Y$$ holds.
\par
(iii) For every noetherian local domain $(R,\goth m)$, every $h>0$
and every sequence $x_1,\dots ,x_h$ of elements of $R$ the
implication
$$\LCMo ^h_{(x_1,\dots ,x_h)}(R)\neq 0\Longrightarrow \{ 0\} \in
\Ass _R(D(\LCMo ^h_{(x_1,\dots ,x_h)R}(R)))$$ holds.
\par
(iv) For every noetherian local ring $(R,\goth m)$, every finitely
generated $R$-module $M$, every $h>0$ and every sequence  $x_1,\dots
,x_h$ of elements of $R$ the equality
$$\Ass _R(D(\LCMo ^h_{(x_1,\dots ,x_h)R}(M)))=\{ \goth p\in
\Supp_R(M)\vert \LCMo ^h_{(x_1,\dots ,x_h)R}(M/\goth pM)\neq 0\}
\leqno{(1)}$$ holds.
\par
Proof:
\par
First we show that (i) -- (iii) are equivalent. (i) $\Longrightarrow
$ (ii): In the given situation we have
$$\Hom _R(R/\goth p_1,D(\LCMo
^h_{(x_1,\dots ,x_h)R}(R)))\neq 0\ \ ;$$
this implies
$$\eqalign {0&\neq \Hom _R(R/\goth p_0,D(\LCMo ^h_{(x_1,\dots
,x_h)R}(R)))\cr &= \Hom _R(\LCMo ^h_{(x_1,\dots ,x_h)R}(R)\otimes
_R(R/\goth p_0) ,\InjH _R(R/\goth m))\cr &=D(\LCMo ^h _{(x_1,\dots
,x_h)R}(R/\goth p_0))\ \ .\cr }$$ Thus conjecture (*) implies that
$\goth p_0$ is associated to $D(\LCMo ^h_{(x_1,\dots ,x_h)R}(R))$.
\par
(ii) $\Longrightarrow $ (iii): We assume that $\LCMo ^h_{(x_1,\dots
,x_h)R}(R)\neq 0$. This implies $D(\LCMo ^h_{(x_1,\dots
,x_h)R}(R))\neq 0$ and hence $\Ass _R(D(\LCMo ^h_{(x_1,\dots
,x_h)R}(R)))\neq \emptyset $; now (ii) shows $\{ 0\} \in \Ass
_R(D(\LCMo ^h_{(x_1,\dots ,x_h)R}(R)))$.
\par
(iii) $\Longrightarrow $ (i): We have seen above that the inclusion
$\subseteq $ holds; we take a prime ideal $\goth p$ of $R$ such that
$\LCMo ^h_{(x_1,\dots ,x_h)R}(R/\goth p)\neq 0$ and we have to show
$\goth p\in \Ass _R(D(\LCMo ^h_{(x_1,\dots ,x_h)R}(R)))$: We apply
(iii) to the domain $R/\goth p$ and get an $R$-linear injection
$$\eqalign {R/\goth p&\to D(\LCMo ^h_{(x_1,\dots ,x_h)(R/\goth
p)}(R/\goth p))\cr &=\Hom _R(\LCMo ^h_{(x_1,\dots ,x_h)R}(R/\goth
p),\InjH _R(R/\goth m))\cr &=\Hom _R(\LCMo ^h_{(x_1,\dots
,x_h)R}(R)\otimes _RR/\goth p,\InjH _R(R/\goth p))\cr &=\Hom
_R(R/\goth p,D(\LCMo ^h_{(x_1,\dots ,x_h)R}(R)))\cr &\subseteq
D(\LCMo ^h_{(x_1,\dots ,x_h)R}(R))\ \ .\cr }$$ Note that we used
$\LCMo ^h_{(x_1,\dots ,x_h)(R/\goth p)}(R/\goth p)=\LCMo
^h_{(x_1,\dots ,x_h)R}(R/\goth p)$ and the fact that $\Hom
_R(R/\goth p,\InjH _R(R/\goth m))$ is an $R/\goth p$-injective hull
of $R/\goth m$.\par Now it is clearly sufficient to show that (i)
implies (iv): $\subseteq $: Every element $\goth p$ of the left-hand
side of identity (1) must contain $\Ann _R(M)$ and hence is an
element of $\Supp _R(M)$; furthermore it satisfies
$$\eqalign {0&\neq \Hom _R(R/\goth p,D(\LCMo ^h_{(x_1,\dots
,x_h)R}(M)))\cr &=\Hom _R(R/\goth p\otimes _R\LCMo ^h_{(x_1,\dots
,x_h)R}(M),\InjH _R(R/\goth m))\cr &=D(\LCMo ^h_{(x_1,\dots
,x_h)R}(M/\goth pM))\ \ .\cr }$$ $\supseteq $: Let $\goth p$ be an
element of the support of $M$ such that $\LCMo ^h_{(x_1,\dots,
x_h)R}(M/\goth pM)$ is not zero. We set $\overline R:=R/\Ann _R(M)$,
$M$ is an $\overline R$-module. $\goth p\supseteq \Ann _R(M)$, we
set $\overline \goth p:=\goth p/\Ann _R(M)$. Clearly our hypothesis
implies that $\LCMo ^h_{(x_1,\dots ,x_h)\overline R}(\overline
R)\neq 0$. We apply (i) to $\overline R$ and deduce
$$\overline \goth p\in \Ass _{\overline R}(D(\LCMo ^h_{(x_1,\dots
,x_h)\overline R}(\overline R)))\ \ .$$
Hence there is an $R$-linear injection
$$0\to R/\goth p=\overline R/\overline \goth p\to D(\LCMo ^h_{(x_1,\dots
,x_h)\overline R}(\overline R))\ \ ,$$
which induces an $R$-linear injection
$$\eqalign {0\to \Hom _R(M,R/\goth p)&\to \Hom _R(M,D(\LCMo ^i_{(x_1,\dots
,x_h)\overline R}(\overline R)))\cr &=\Hom _{\overline R}(M,D(\LCMo
^h_{(x_1,\dots ,x_h)\overline R}(\overline R)))\cr &=D(\LCMo
^h_{(x_1,\dots ,x_h)\overline R}(M))\cr &=D(\LCMo ^h_{(x_1,\dots
,x_h)R}(M))\ \ .\cr }$$ Note that for the second equality we have
used Hom-Tensor adjointness and for the last equality the facts that
$M$ is an $\overline R$-module and that $\Hom _R(\overline R,\InjH
_R(R/\goth m))$ is an $\overline R$-injective hull of $R/\goth m$;
It is sufficient to show $\goth p\in \Ass _R(\Hom _R(M,R/\goth p))$;
but $M$ is finite and so we have
$$(\Hom _R(M,R/\goth p))_\goth p=\Hom
_{R_\goth p}(M_\goth p,R_\goth p/\goth pR_\goth p)\neq 0\ \ ,$$
which
shows that $\goth pR_\goth p$ is associated to the $R_\goth
p$-module $(\Hom _R(M,R/\goth p))_\goth p$. Thus $\goth p\in \Ass
_R(\Hom _R(M,R/\goth p))$.
\bigskip
{\bf 1.2.4 Remark}
\par
In [HS1, section 0, conjecture (+)] more was conjectured, namely:
\smallskip
If $(R,\goth m)$ is a noetherian local ring, $h\geq 1$ and
$x_1,\dots ,x_h$ is a sequence of elements of $R$, then all prime
ideals $\goth p$ maximal in $\Ass _R(D(\LCMo ^h_{(x_1,\dots
,x_h)R}(R)))$ have the same dimension, namely $\dim (R/\goth p)=h$.
\smallskip
This conjecture is false, here is a counterexample:
\par
Let $\bf Q$ denote the rationals and $R={\bf Q}[[X_1,X_2,X_3,X_4,X_5]]$ a power
series algebra over $\bf Q$ in the variables $X_1,\dots ,X_5$. Set $h=3$ and $x_1=X_1,x_2=X_2,x_3=X_3$.
Then
$$p:=-X_2X_4^2+X_3X_4X_5-X_1X_5^2+4X_1X_2-X_3^2\in R$$
is a prime element of $R$; in fact, $pR$ is a maximal element of
$\Ass _R(D(\LCMo ^h_{(x_1,x_2,x_3)R}(R)))$, but $\dim (R/fR)=4\neq
3$. These statements will be proved in remark 4.3.2 (ii), here we
explain where $f$ comes from:
\par
We define a ring
$$S:={\bf Q}[[y_1,y_2,y_3,y_4]]$$
and a module-finite $\bf Q$-algebra homomorphism
$$f:R\to S$$
such that
$$f(X_1)=y_1y_3, f(X_2)=y_2y_4, f(X_3)=y_1y_4+y_2y_4, f(X_4)=y_1+y_3,
f(X_5)=y_2+y_4)\ \ .$$ As we will see in remark 4.3.2 (ii), $pR$ is
the kernel of $f$; the crucial point here is that the radical of the
extension ideal of $(x_1,x_2,x_3)R$ in $R_0$ is
$$I_0=(y_1,y_2)R_0\cap (y_3,y_4)R_0$$
and $\LCMo ^3_{I_0}(R_0)\neq 0$,
although $I_0$ has height two (again, see section 4.3 and, in
particular, remark 4.3.2 (ii) for details).
\bigskip
\bigskip
{\bf 1.3 Regular sequences on $D(\LCMo ^h_I(R))$ are well-behaved in
some sense}
\bigskip
Obviously we are dealing here with the notion of regular sequences
on modules which are, in general, not finitely generated. Such
regular sequences do not have all the good properties of regular
sequences on finite modules. However, in our situation, some kind of
well-behavior holds, here is the idea (see theorem 1.3.1 below for
the precise statement): For finite modules, the following is
well-known: If $(R,\goth m)$ is a noetherian local ring, $M$ a
finite $R$-module and $r_1,\dots ,r_h\in R$ an $M$-regular sequence
then $r_1^\prime ,\dots ,r_h^\prime \in R$ is also an $M$-regular
sequence provided $\sqrt {(r_1^\prime ,\dots ,r_h^\prime
)R}=\sqrt{(r_1,\dots ,r_h)R}$ holds (because $R$ is local). In our
case, if $(R,\goth m)$ is a noetherian local ring and $I\subseteq R$
an ideal of $R$ such that $\LCMo ^l_I(R)\neq 0\iff l=h$ holds, it is
clear that if an $R$-regular sequence $r_1,\dots ,r_h\in I$ is a
$D(\LCMo ^h_I(R))$-regular sequence then an $R$-regular-sequence
$r_1^\prime ,\dots ,r_h^\prime \in I$ is also $D(\LCMo
^h_I(R))$-regular if $\sqrt{(r_1^\prime ,\dots ,r_h^\prime
)R}=\sqrt{(r_1,\dots ,r_h)R}$ holds (simply because of $\sqrt
{(r_1^\prime ,\dots ,r_h^\prime )R}=\sqrt I$ and corollary 1.1.4).
But a more sophisticated statement is also true:
\bigskip
{\bf 1.3.1 Theorem}
\smallskip
Let $(R,\goth m)$ be a noetherian local ring, $h\geq 1$ and
$I\subseteq R$ an ideal such that $\LCMo ^l_I(R)\neq 0\iff l=h$
holds. Furthermore, let $1\leq h^\prime \leq h$ and let $r_1,\dots
,r_{h^\prime }\in I$ be an $R$-regular sequence that is also
$D(\LCMo ^h_I(R))$-regular. Furthermore, let $r_1^\prime ,\dots
,r_{h^\prime }^\prime \in I$ be such that $\sqrt{(r_1^\prime , \dots
,r_{h^\prime }^\prime )R}=\sqrt{(r_1,\dots ,r_{h^\prime })}R$ holds.
Then $r_1^\prime ,\dots ,r_{h^\prime }^\prime $ is a $D(\LCMo
^h_I(R))$-regular sequence. In particular, any permutation of
$r_1,\dots ,r_{h^\prime }$ is again a $D(\LCMo ^h_I(R))$-regular
sequence.
\par
Proof:
\par
$R$ is local, and thus it is clear that $r_1^\prime ,\dots
,r_{h^\prime }^\prime $ is an $R$-regular sequence. By induction on
$s\in \{ 1,\dots ,h^\prime \} $ we show two statements:
$$\LCMo ^l_I(R/(r_1,\dots ,r_s)R)\neq 0\iff l=h-s$$
and
$$D(\LCMo ^{h-s}_I(R/(r_1,\dots ,r_s)R))=D(\LCMo ^h_I(R))/(r_1,\dots
,r_s)D(\LCMo ^h_I(R))\ \ :$$ $s=1$: The short exact sequence
$$0\to R\buildrel r_1\over \to R\to R/r_1R\to 0$$
induces a short exact sequence
$$0\to \LCMo ^{h-1}_I(R/r_1R)\to \LCMo ^h_I(R)\buildrel r_1\over \to \LCMo
^h_I(R)\to 0$$ and we conclude, therefore, that
$$\LCMo ^l_I(R/r_1R)\neq 0\iff l=h-1$$
holds. Now, the statement
$$D(\LCMo ^{h-1}_I(R/r_1R))=D(\LCMo
^h_I(R))/r_1D(\LCMo ^h_I(R))$$
follows from the exactness of $D$.
\par
$s>1$: The short exact sequence
$$0\to R/(r_1,\dots ,r_{s-1})R\buildrel r_s\over
\to R/(r_1,\dots ,r_{s-1})R\to R/(r_1,\dots ,r_s)R\to 0$$ induces,
by our induction hypothesis, an exact sequence
$$0\to \LCMo ^{h-s}_I(R/(r_1,\dots ,r_s)R)\to \LCMo ^{h-(s-1)}_I(R/(r_1,\dots
,r_{s-1})R)\buildrel r_s\over \to \LCMo ^{h-(s-1)}_I(R/(r_1,\dots
,r_{s-1})R)\ \ .$$ By induction hypothesis,
$$D(\LCMo
^{h-(s-1)}_I(R/(r_1,\dots ,r_{s-1})R))=D(\LCMo ^h_I(R))/(r_1,\dots
,r_{s-1})D(\LCMo ^h_I(R))$$
and so, by assumption, $r_s$ operates
surjectively on $\LCMo ^{h-(s-1)}_I(R/(r_1,\dots ,r_{s-1})R)$ and we
get
$$\LCMo ^l_I(R/(r_1,\dots ,r_s)R)\neq 0\iff l=r-s$$
and
$$\eqalign {D(\LCMo ^{h-s}_I(R/(r_1,\dots
,r_s)R))&=D(\LCMo ^{h-(s-1)}_I(R/(r_1,\dots ,r_{s-1})R))/r_sD(\LCMo
^{h-(s-1)}_I(R/(r_1,\dots ,r_{s-1})R))\cr &=D(\LCMo
^h_I(R))/(r_1,\dots ,r_s)D(\LCMo ^h_I(R))\ \ .\cr }$$ In particular
for $s=h^\prime $ we have
$$\LCMo ^l_I(R/(r_1,\dots ,r_{h^\prime })R)\neq 0\iff l=h-h^\prime \ \ .$$
Note that, because of
$$\depth (I,R/(r_1,\dots ,r_{h^\prime
})R)=\depth (I,R)-h^\prime =\depth (I,R/(r_1^\prime ,\dots
,r_{h^\prime }^\prime )R)$$ and
$$\Supp _R(R/(r_1,\dots ,r_{h^\prime
})R)=\Supp _R(R/(r_1^\prime ,\dots ,r_{h^\prime }^\prime )R$$ this
implies
$$\LCMo ^l_I(R/(r_1^\prime ,\dots ,r_{h^\prime }^\prime )R)\neq 0\iff
l=h-h^\prime \leqno{(1)} $$ (the $\depth $-argument shows vanishing
for $l<h-h^\prime $ and the $\Supp $-argument shows that $h-h^\prime
$ is the largest number such that $\LCMo ^l_I(R/(r_1^\prime ,\dots
,r_{h^\prime }^\prime )R)\neq 0$). Now, by descending induction on
$s\in \{ 0,\dots ,h^\prime -1\} $, we will prove the following three
statements:
$$r_{s+1}^\prime \hbox { operates surjectively on } \LCMo
^{h-s}_I(R/(r_1^\prime ,\dots ,r_s^\prime )R)\ \ ,$$
$$\LCMo ^{h-l}_I(R/(r_1^\prime ,\dots ,r_s^\prime )R)\neq 0\iff l=s$$
and
$$D(\LCMo ^{h-(s+1)}_I(R/(r_1^\prime ,\dots ,r_{s+1}^\prime )R))=D(\LCMo
^{h-s}_I(R/(r_1^\prime ,\dots ,r_s^\prime )R))/r_{s+1}^\prime
D(\LCMo ^{h-s}_I(R/(r_1^\prime ,\dots ,r_s^\prime )R))\ \ :$$
$s=h^\prime -1$: We consider the long exact $\Gamma _I$-sequence
belonging to the short exact sequence
$$0\to R/(r_1^\prime ,\dots ,r_{h^\prime -1}^\prime )R\buildrel r_{h^\prime
}^\prime \over \to R/(r_1^\prime ,\dots ,r_{h^\prime -1}^\prime
)R\to R/(r_1^\prime ,\dots ,r_{h^\prime }^\prime )R\to 0\ \ :$$
Then, the surjectivity of $r_{h^\prime }^\prime $ on $\LCMo
^{h-(h^\prime -1)}_I(R/(r_1^\prime ,\dots ,r_{h^\prime -1}^\prime
)R)$ follows from (1) and the other statements from the fact that
for $l\neq h-(h^\prime -1)$ we have injectivity of $r_{h^\prime
}^\prime $ on $\LCMo ^l_I(R/(r_1^\prime ,\dots ,r_{h^\prime
-1}^\prime )R)$, hence
$$\LCMo ^l_I(R/(r_1^\prime ,\dots
,r_{h^\prime -1}^\prime )R)=0$$ as $r_{h^\prime }^\prime \in I$.
\par
$s<h^\prime -1$: We consider the long exact $\Gamma _I$-sequence
belonging to the short exact sequence
$$0\to R/(r_1^\prime ,\dots ,r_s^\prime )R\buildrel r_{s+1}^\prime \over \to
R/(r_1^\prime ,\dots ,r_s^\prime )R\to R/(r_1^\prime ,\dots
,r_{s+1}^\prime )R\to 0\ \ :$$ Then, our induction hypothesis shows
that multiplication by $r_{s+1}^\prime $ is surjective on $\LCMo
^{h-s}_I(R/(r_1^\prime ,\dots ,r_s^\prime))$. Like before, the two
other statements follow from the fact that, for $l\neq h-s$,
multiplication by $r_{s+1}^\prime $ is injective on $\LCMo
^{h-l}_I(R/(r_1^\prime ,\dots ,r_s^\prime ))$ and so $\LCMo
^{h-l}_I(R/(r_1^\prime ,\dots ,r_s^\prime))$ is trivial. It is clear
that these three statements prove the theorem (in fact, the first
and the third statement are sufficient here, the second is used for
technical reasons).
\bigskip
\bigskip
{\bf 1.4 Comparison of two Matlis Duals} \bigskip For a noetherian
local ring $(R,\goth m)$, the Matlis dual functor clearly depends on
$R$. In this section we will have a local subring $R_0$ of $R$.
Given any local cohomology module over $R$, we will take its Matlis
dual both with respect to $R$ and with respect to $R_0$; both are
$R$-modules in a natural way. Among other results, in this section
we will see that, under certain assumptions, these two Matlis duals
have the same set of associated prime ideals (over $R$, see 1.4.3
(ii)).
\bigskip Let $(R,\goth m)$ be a noetherian local equicharacteristic
complete ring with coefficient field $k$ and let $y_1,\dots ,y_i$ be
a sequence in $R$ such that $R_0:=k[[y_1,\dots ,y_i]]$ is regular
and of dimension $i$ (this is true, for example, if $\LCMo
^i_{(y_1,\dots ,y_i)R}(R)\neq 0$ holds, as this local cohomology
module agrees with $\LCMo^i_{(y_1,\dots,x_i)R_0}(R_0)\otimes
_{R_0}R$; also note that $R_0$ is, by definition, a subring of $R$).
Let $D_R$ denote the Matlis-dual functor with respect to $R$ and
$D_{R_0}$ the one with respect to $R_0$. By local duality (see e. g.
[BS, section 11] for a reference on local duality), we get
$$D_{R_0}(\LCMo ^i_{(y_1,\dots ,y_i)R}(R))=\Hom _{R_0}(R\otimes _{R_0}\LCMo
^i_{(y_1,\dots ,y_i)R_0}(R_0),\InjH _{R_0}(k))=\Hom _{R_0}(R,R_0)\ \
.$$ $\Hom _{R_0}(R,\InjH _{R_0}(k))$ is an injective $R$-module with
non-trivial socle; therefore, there exists an injective $R$-module
$E^\prime $ such that
$$\Hom _{R_0}(R,\InjH _{R_0}(k))=\InjH _R(k)\oplus E^\prime $$
holds. We set $E=\Gamma _{(y_1,\dots ,y_i)R}(E^\prime )$. We have
$$\eqalign {D_{R_0}(\LCMo ^i_{(y_1,\dots ,y_i)R}(R))&=\Hom
_{R_0}(\LCMo ^i_{(y_1,\dots ,y_i)R_0}(R_0)\otimes _{R_0}R,\InjH
_{R_0}(k))\cr &=\Hom _{R_0}(\LCMo ^i_{(y_1,\dots ,y_i)R_0}(R_0),\Hom
_{R_0}(R,\InjH _{R_0}(k)))\cr &=\Hom _R(R\otimes _{R_0}\LCMo
^i_{(y_1,\dots ,y_i)R_0}(R_0),\Hom _{R_0}(R,\InjH _{R_0}(k)))\cr
&=\Hom _R(\LCMo ^i_{(y_1,\dots ,y_i)R}(R),\InjH _R(k)\oplus E^\prime
)\cr &=D_R(\LCMo ^i_{(y_1,\dots ,y_i)R}(R))\oplus \Hom _R(\LCMo
^i_{(y_1,\dots ,y_i)R}(R),E^\prime )\cr &=D_R(\LCMo ^i_{(y_1,\dots
,y_i)R}(R))\oplus \Hom _R(\LCMo ^i_{(y_1,\dots ,y_i)R}(R),E)\cr }$$
and hence
$$\Ass _R(D_{R_0}(\LCMo ^i_{(y_1,\dots ,y_i)R}(R)))=\Ass _R(D_R(\LCMo ^i_{(y_1,\dots
,y_i)R}(R)))\cup \Ass _R(\Hom _R(\LCMo ^i_{(y_1,\dots
,y_i)R}(R),E))\ \ .\leqno{(1)}$$ It is natural to ask for relations
between $D_R(\LCMo ^i_{(y_1,\dots ,y_i)R}(R))$ and $D_{R_0}(\LCMo
^i_{(y_1,\dots ,y_i)R}(R))$; we will establish some in the sequel:
\bigskip
For every $\goth p\in Z:=\{ \goth p\in \Spec (R)\vert (y_1,\dots
,y_i)R\subseteq \goth p\subsetneq \goth m\} $ we choose a set
$\mu_\goth p$ such that
$$E=\bigoplus _{\goth p\in Z}\InjH _R(R/\goth p)^{(\mu_\goth p)}$$
holds.
\bigskip
{\bf 1.4.1 Remark}
\par
In the above situation, one has $\mu_\goth p\neq \emptyset $ for
every $ \goth p\in Z$.
\par
Proof:
\par
We have to show that $\goth p$ is associated to the $R$-module $\Hom
_{R_0}(R/\goth p,\InjH _{R_0}(k))$. The latter module is equal to
$\Hom _{R_0}(R/\goth p,k)$, because $\goth p$ is annihilated by
$y_1,\dots ,y_i$ (note that $k$ is the socle of $\InjH _{R_0}(k)$).
Thus we have to prove the following statement: If $(R,\goth m)$ is a
noetherian local equicharacteristic complete domain with coefficient
field $k$, then the zero ideal of $R$ is associated to the
$R$-module $\Hom _k(R,k)$:
\par
Let $x_1,\dots ,x_n\in R$ be a system of parameters for $R$,
$n:=\dim (R)$. Then $R_0:=k[[x_1,\dots ,x_n]]$ is a regular subring
of $R$, over which $R$ is module-finite. One has $\Hom _k(R,k)=\Hom
_{R_0}(R,\Hom _k(R_0,k))$ and, therefore, it is sufficient to prove
$\{ 0\} \in \Ass _{R_0}(\Hom _k(R_0,k))$, because in this case,
every $R_0$-injection
$$R_0\to \Hom _k(R_0,k)$$
induces an $R$-injection
$$\Hom _{R_0}(R,R_0)\to \Hom _k(R,k)$$
and $\{ 0\} \in \Supp _R(\Hom _{R_0}(R,R_0))$ holds, because $R$ is
finite over $R_0$. Thus we may assume $R=k[[x_1,\dots ,x_n]]$ from
now on:
\par
For $i=1,\dots ,n$ we set $R_i:=k[[x_1,\dots ,x_i]]$. Again we have
$$\Hom _k(R_i,k)=\Hom _{R_{i-1}}(R_i,\Hom _k(R_{i-1},k))$$
for $i=2,\dots ,n$. Using this and an obvious induction argument,
the statement follows from lemma 1.4.2 below.
\bigskip
{\bf 1.4.2 Lemma}
\par
Let $k$ be a field and let $R_0:=k[[X_1,\dots ,X_n]]$,
$R:=k[[X_1,\dots ,X_n,X]]=R_0[[X]]$ be power series rings in the
variables $X_1,\dots ,X_n,X$, respectively. Then
$$\{ 0\} \in \Ass _R(\Hom _{R_0}(R,R_0))\ \ .$$
Proof:
\par
By $\goth m_0$ we denote the maximal ideal of $R_0$. The canonical
short exact sequence
$$0\to R_0[X]\to R_0[[X]]\to R_0[[X]]/R_0[X]\to 0$$
induces an exact sequence
$$0\to \Hom _{R_0}(R_0[[X]]/R_0[X],R_0)\to \Hom _{R_0}(R_0[[X]],R_0)\buildrel
\alpha \over \to \Hom _{R_0}(R_0[X],R_0)\ \ .$$ The map $\alpha $ is
the Matlis dual (in the sense that
$$\Hom _{R_0}(\LCMo ^n_{\goth
m_0}(R_0[X]),\InjH _{R_0}(k))=\Hom _{R_0}(R_0[X]\otimes _{R_0}\LCMo
^n_{\goth m_0}(R_0),\InjH _{R_0}(k))=\Hom _{R_0}(R_0[X],R_0)$$ and
$$\Hom _{R_0}(\LCMo ^n_{\goth
m_0}(R_0[[X]]),\InjH _{R_0}(k))=\Hom _{R_0}(R_0[[X]]\otimes
_{R_0}\LCMo ^n_{\goth m_0}(R_0),\InjH _{R_0}(k))=\Hom
_{R_0}(R_0[[X]],R_0)$$ hold) of the canonical map
$$\LCMo ^n_{\goth m_0}(R_0[X])=\LCMo ^n_{\goth
m_0}(R_0)\otimes _{R_0}R_0[X]\to \LCMo ^n_{\goth m_0}(R_0)\otimes
R_0[[X]]=\LCMo ^n_{\goth m_0}(R_0[[X]])\ \ ,$$ which is obviously
injective. This means that $\alpha $ is surjective. The
$R_0[X]$-module $\Hom _{R_0}(R_0[X],R_0)$ can be written as
$R_0[[X^{-1}]]$ and in $R_0[[X^{-1}]]$ the element
$$h^\prime
:=1+X^{-1!}+X^{-2!}+\dots $$
has $R_0[X]$-annihilator zero
(essentially because the sequence of differences $2!-1!,3!-2!,\dots
$ becomes arbitrary large). Choose an element $h\in \Hom
_{R_0}(R_0[[X]],R_0)$ which is mapped to $h^\prime $ by $\alpha $
Then $\ann _{R_0[X]}(h)=\{ 0\} $ which implies $\ann
_{R_0[[X]]}(h)=\{ 0\} $, using a flatness argument.
\bigskip
{\bf 1.4.3 Remarks}
\par
Let $(R,\goth m)$ be a noetherian local ring and $y_1,\dots ,y_i$ a
sequence in $R$ and suppose that conjecture (*) holds.
\par
(i) For every fixed prime ideal $\goth p$ of $R$, one has
$$\Ass _R(\Hom _R(\LCMo ^i_{(y_1,\dots ,y_i)R}(R),\InjH _R(R/\goth p)))=\{
\goth q\in \Spec (R)\vert (\LCMo ^i_{(y_1,\dots ,y_i)R}(R/\goth
q))_\goth p\neq 0\} \ \ .$$ (ii) For every prime ideal $\goth p$ of
$R$, let $\nu _\goth p$ be a set. Then
$$\Ass _R(\Hom _R(\LCMo ^i_{(y_1,\dots ,y_i)R}(R),\bigoplus _{\goth
p\in \Spec (R)}\InjH _R(R/\goth p)^{(\nu _\goth p)}))=\bigcup _{\nu
_\goth p\neq \emptyset }\Ass _R(\Hom _R(\LCMo ^i_{(y_1,\dots
,y_i)R}(R),\InjH _R(R/\goth p)))\ \ .$$ As a consequence, in the
situation of (1) (note that then we had more assumptions: $R$ is
complete and equicharacteristic and $\underline y$ is such that
$\LCMo ^i_{(y_1,\dots ,y_i)R}(R)\neq 0$), one has
$$\Ass _R(D_{R_0}(\LCMo ^i_{(y_1,\dots ,y_i)R}(R)))=\Ass _R(D_R(\LCMo
^i_{(y_1,\dots ,y_i)R}(R)))\ \ .$$ Proof:
\par
(i) For every prime ideal $\goth p$ of $R$, $\InjH _R(R/\goth
p)=\InjH _{R_\goth p}(R_\goth p/\goth pR_\goth p)$ is naturally an
$R_\goth p$-module. This implies
$$\Hom _R(\LCMo ^i_{(y_1,\dots ,y_i)R}(R),\InjH _R(R/\goth p))=\Hom
_{R_\goth p}(\LCMo ^i_{(y_1,\dots ,y_i)R_\goth p}(R_\goth p),\InjH
_{R_\goth p}(R_\goth p/\goth pR_\goth p))$$ and, therefore and
because of (*),
$$\eqalign {\Ass _R(\Hom _R(\LCMo ^i_{(y_1,\dots ,y_i)R}(R),\InjH _R(R/\goth
p)))&=\{ \goth P\cap R\vert \goth P\in \Ass _{R_\goth p}(\Hom
_{R_\goth p}(\LCMo ^i_{(y_1,\dots ,y_i)R_\goth p}(R_\goth p),\InjH
_{R_\goth p}(R_\goth p/\goth pR_\goth p)))\} \cr &=\{ \goth q\in
\Spec(R)\vert \LCMo ^i_{(y_1,\dots ,y_i)R}(R/\goth q)_\goth p\neq
0\} \ \ .\cr }$$ (ii) We have natural inclusions
$$\eqalign {\bigoplus _{\goth p\in \Spec (R)}\Hom _R(\LCMo ^i_{(y_1,\dots
,y_i)R}(R),\InjH _R(R/\goth p))^{(\nu _\goth p)}&\subseteq \Hom
_R(\LCMo ^i_{(y_1,\dots ,y_i)R}(R),\bigoplus _{\goth p\in \Spec
(R)}\InjH _R(R/\goth p)^{(\nu _\goth p)})\cr &\subseteq \prod
_{\goth p\in \Spec (R)}\Hom _R(\LCMo ^i_{(y_1,\dots ,y_i)R}(R),\InjH
_R(R/\goth p))^{\nu _\goth p}\ \ .}$$
Every annihilator of a
non-trivial element of $\prod _{\goth p\in \Spec (R)}\Hom _R(\LCMo
^i_{(y_1,\dots ,y_i)R}(R),\InjH _R(R/\goth p))^{\nu_\goth p}$ is
contained in some associated prime ideal of some $\Hom _R(\LCMo
^i_{(y_1,\dots ,y_i)R}(R),\InjH _R(R/\goth p))$, where $\nu _\goth
p\neq \emptyset $. But the set
$$\Ass _R(\Hom _R(\LCMo ^i_{(y_1,\dots ,y_i)R}(R),\InjH
_R(R/\goth p)))$$ is stable under generalization because of the
conjecture (*). Therefore, we get $$\eqalign {\Ass _R(\Hom _R(\LCMo
^i_{(y_1,\dots ,y_i)R}(R),\bigoplus _{\goth p\in \Spec(R)}\InjH
_R(R/\goth p)^{(\nu _\goth p)})) &=\Ass _R(\bigoplus _{\goth p\in
\Spec (R)}\Hom _R(\LCMo ^i_{(y_1,\dots ,y_i)R}(R),\InjH _R(R/\goth
p))^{(\nu_\goth p)})\cr &=\bigcup _{\goth p\in Z}\Ass _R(\Hom
_R(\LCMo ^i_{(y_1,\dots ,y_i)R}(R),\InjH _R(R/\goth p)))\cr
&\subseteq \Ass _R(D_R(\LCMo ^i_{(y_1,\dots ,y_i)R}(R)))\ \ .\cr }$$
In particular, in the situation of (1), we have
$$\Ass _R(D_{R_0}(\LCMo ^i_{(y_1,\dots ,y_i)R}(R)))=\Ass _R(D_R(\LCMo
^i_{(y_1,\dots ,y_i)R}(R)))\ \ .$$ \vfil \eject {\bf 2 Associated
primes -- a constructive approach}
\bigskip
In this section we will prove results on the set
$$\Ass _R(D(\LCMo ^i_{(x_1,\dots ,x_i)R}(R)))\ \ ,$$
where $\underline x=x_1,\dots ,x_i$ is a sequence in a noetherian
local ring $R$. The proofs are based on the fact that, over the
formal power series ring $R=k[[X_1,\dots ,X_n]]$ ($k$ a field), the
$R$-module
$$E=k[X_1^{-1},\dots ,X_n^{-1}]$$
is an $R$-injective hull of $k$. The methods in this sections are
constructive to some extent, in fact, we construct certain elements
in $k[X_1^{-1},\dots ,X_n^{-1}]$. For the proofs, we will have to
distinguish between the equicharacteristic and the
mixed-characteristic case. One major result in this section is
(theorem 2.4, see also theorem 2.5 for the case of mixed
characteristic):
\par
If $\underline x=x_1,\dots x_i$ is a sequence in a noetherian local
equicharacteristic ring $(R,\goth m)$ and $\underline x$ is part of
a system of parameters of $R/\goth p$ for some fixed prime ideal
$\goth p$ of $R$, then one has
$$\goth p\in \Ass _R(D(\LCMo ^i_{(x_1,\dots ,x_i)R}(R)))\ \ .$$
We will also see that, in general, not all associated primed of
$D:=D(\LCMo ^i_{(x_1,\dots ,x_i)R}(R))$ are obtained in this way
(remark 2.7 (ii)). As a corollary, we are able to completely compute
the set $\Ass _R(D)$ in the case $i=1$ (corollary 2.6):
$$\Ass _R(D(\LCMo ^1_{xR}(R)))=\Spec (R)\setminus {\cal V}(x)$$
(note that ${\cal V}(x)=\{ \goth p\in \Spec(R)\vert x\in \goth p\}
$). In particular, this set is infinite (in general). The sections
ends with remarks on the questions of stableness under
generalization of the subsets
$$Z_1:=\{ \goth p\in \Spec (R)\vert \LCMo ^i_{(x_1,\dots ,x_i)R}(R/\goth
p)\neq 0\} $$ and
$$Z_2:=\{ \goth p\in \Spec (R)\vert x_1,\dots ,x_i\hbox { is part of a
system of parameters for }R/\goth p\}$$ of $\Spec(R)$. Note that we
have
$$Z_2\subseteq \Ass _R(D)\subseteq Z_1$$
by theorems 2.4, 2.5 and remark 1.1.2.
\bigskip
We start with a special case of the result mentioned above:
\bigskip
{\bf 2.1 Lemma}
\par
Let $k$ be a field, $n\geq 1$, $R=k[[X_1,\dots ,X_n]]$ and $i\in \{
1,\dots ,n\} $. We set $I:=(X_1,\dots ,X_i)R$ and $\goth
m:=(X_1,\dots ,X_n)R$. Then
$$\{ 0\} \in \Ass _R(D(\LCMo ^i_I(R)))$$
holds. \par
Proof:
\par
1. Case: $i=n$:
\par
Here $\LCMo ^i_I(R)=\InjH _R(R/\goth m)$ und also $D(\LCMo
^i_I(R))=R$ and the statement follows.
\par
2. Case: $i<n$: We have
$$\LCMo ^i_I(R)=\vtop{\baselineskip=1pt
\lineskiplimit=0pt \lineskip=1pt\hbox{lim} \hbox{$\longrightarrow $}
\hbox{$^{^{l\in \bf N \setminus \{ 0\} }}$}}(R/(X_1^l,\dots
,X_i^l)R)\ \ ,$$ the transition maps being induced by $R\to R$,
$r\mapsto (X_1\cdot \dots \cdot X_i)\cdot r$. So
$$D(\LCMo
^i_I(R))=\vtop{\baselineskip=1pt \lineskiplimit=0pt
\lineskip=1pt\hbox{lim} \hbox{$\longleftarrow $} \hbox{$^{^{l\in \bf
N\setminus \{ 0\} }}$}}(D(R/(X_1^l,\dots ,X_i^l)R))\ \ ;$$
here
$$D(R/(X_1^l,\dots ,X_i^l)R)=\Hom _R(R/(x_1^l,\dots ,x_i^l),
D(R))=\InjH _{R/(X_1^l,\dots ,X_i^l)R}(R/\goth m)(\subseteq \InjH
_R(R/\goth m))\ \ ,$$
the transition maps being induced by $\InjH
_R(R/\goth m)\to \InjH _R(R/\goth m)$, $e\mapsto (X_1\cdot \dots
\cdot X_i)\cdot e$ and we have $\InjH _R(R/\goth m)=k[X_1^{-1},\dots
,X_n^{-1}]$ (by definition, the last module is the $k$-vector space
with basis $(X_1^{i_1}\cdot \dots \cdot X_n^{i_n})_{i_1,\dots
,i_n\leq 0}$ and with an obvious $R$-module structure on it). We
define
$$\eqalign {\alpha :=&(1,X_1^{-1}\cdot \dots \cdot X_i^{-1}+X_{i+1}^{-1!}\cdot \dots
\cdot X_n^{-1!},\dots ,X_1^{-m}\cdot \dots \cdot X_i^{-m}+\cr
&+(X_{i+1}^{-1!}\cdot \dots \cdot X_n^{-1!})\cdot (X_1^{-(m-1)}\cdot
\dots \cdot X_i^{-(m-1)})+\dots +\cr &+(X_{i+1}^{-(m-1)!}\cdot \dots
\cdot X_n^{-(m-1)!})\cdot (X_1^{-1}\cdot \dots \cdot
X_i^{-1})+X_{i+1}^{-m!}\cdot \dots \cdot X_n^{-m!},\dots )\in
D(\LCMo ^i_I(R))\ \ .\cr }$$ Here we consider the projective limit
as a subset of a direct product. We state $\ann _R(\alpha )=\{ 0\}$:
Assume there is an $f\in \ann _R(\alpha )\setminus \{ 0\} $. We
choose $(a_1,\dots ,a_n)\in \Supp (f)$ such that $(a_1,\dots ,a_i)$
is minimal (using the ordering
$$(c_1,\dots c_i)\leq (c_1^\prime
,\dots ,c_i^\prime ):\iff c_1\leq c_1^\prime \wedge \dots \wedge
c_i\leq c_i^\prime \ \ )$$ in
$$\{ (a_1^\prime ,\dots ,a_i^\prime
)\vert \exists a_{i+1}^\prime ,\dots ,a_n^\prime :(a_1^\prime ,\dots
,a_n^\prime )\in \Supp (f)\} \ \ .$$
We may assume $a_1=\max \{
a_1,\dots ,a_i\} $. We replace $f$ by $X_2^{a_1-a_2}\cdot \dots
\cdot X_i^{a_1-a_i}\cdot f$; this means $a_1=\dots =a_i=:a$. Choose
$h_1,\dots ,h_i\in R$ and $g \in k[[X_{i+1},\dots ,X_n]]\setminus \{
0\} $ such that
$$f=X_1^{a+1}h_1+\dots +X_i^{a+1}h_i+(X_1^a\cdot \dots \cdot X_i^a)\cdot g$$
$f\cdot \alpha =0$ means: For every $m$ we have
$$\eqalign {0&=[ X_1^{a+1}h_1+\dots +X_i^{a+1}h_i+(X_1^a\cdot \dots \cdot X_i^a)\cdot
g] \cdot (X_1^{-m}\cdot \dots \cdot X_i^{-m}+\dots
+X_{i+1}^{-m!}\cdot \dots \cdot X_n^{-m!})\cr &=(X_1^{a+1}h_1+\dots
+X_i^{a+1}h_i)\cdot [X_1^{-m}\cdot \dots \cdot X_i^{-m}+\dots
+(X_{i+1}^{-(m-a-1)!}\cdot \dots \cdot X_n^{-(m-a-1)!})\cdot \cr
&\cdot (X_1^{-(a+1)}\cdot \dots \cdot X_i^{-(a+1)})]+g\cdot
(X_1^{-(m-a)}\cdot \dots \cdot X_i^{-(m-a)}+\dots
+X_{i+1}^{-(m-a)!}\cdot \dots \cdot X_n^{-(m-a)!})\ \ .\cr }$$
Choose $(b_{i+1},\dots ,b_n)$ minimal in $\Supp (g)$; then for all
$m>>0$ the following statements must hold:
$$\eqalign {(m-a)!-b_{i+1}&\leq (m-a-1)!\cr &\vdots \cr (m-a)!-b_n&\leq
(m-a-1)!\cr }$$ For $m>>0$ this leads to a contradiction, the
assumption is wrong and the lemma is proven.
\bigskip
{\bf 2.2 Lemma}
\par
Let $p$ be a prime number, $C$ a complete $p$-ring, $n\geq 1$,
$R=C[[X_1,\dots ,X_n]]$ and $i\in \{ 1,\dots ,n\} $. We set
$I:=(X_1,\dots ,X_i)R$ and $\goth m:=(p,X_1,\dots ,X_n)R$. Then
$$\{ 0\} \in \Ass _R(D(\LCMo ^i_I(R)))$$
holds.
\par
Proof:
\par
We have
$$\LCMo ^i_I(R)=\vtop{\baselineskip=1pt \lineskiplimit=0pt
\lineskip=1pt\hbox{lim} \hbox{$\longrightarrow $} \hbox{$^{^{l\in
\bf N \setminus \{ 0\} }}$}}(R/(X_1^l,\dots ,X_i^l)R)\ \ ,$$ the
transition maps being induced by $R\to R$, $r\mapsto (X_1\cdot \dots
\cdot X_i)\cdot r$. We deduce
$$D(\LCMo
^i_I(R))=\vtop{\baselineskip=1pt \lineskiplimit=0pt
\lineskip=1pt\hbox{lim} \hbox{$\longleftarrow $} \hbox{$^{^{l\in \bf
N\setminus \{ 0\} }}$}}(D(R/(X_1^l,\dots ,X_i^l)R))\ \ ;$$
we recall
$$D(R/(X_1,\dots ,X_i^l)R)=\InjH _{R/(X_1^l,\dots X_i^l)R}(R/\goth
m)(\subseteq \InjH _R(R/\goth m))\ \ ,$$
the transition maps being
induced by $\InjH _R(R/\goth m)\to \InjH _R(R/\goth m)$, $e\mapsto
(X_1\cdot \dots \cdot X_i)\cdot e$. Furthermore
$$\InjH _R(R/\goth m)=(C_p/C)[X_1^{-1},\dots ,X_n^{-1}]$$ holds (because of
$$\eqalign{\InjH _R(R/\goth m)&=\LCMo ^{n+1}_{(p,X_1,\dots ,X_n)R}(R)\cr &=\LCMo
^1_{pR}(R)\otimes _R\dots \otimes _R\LCMo ^1_{X_nR}(R)\cr
&=(C_p/C)\otimes _C((R_{X_1}/R)\otimes _R\dots \otimes
_R(R_{X_n}/R))\ \ .}$$
We define
$$\eqalign {\alpha :=&(p^{-1},p^{-1}X_1^{-1}\cdot \dots \cdot
X_i^{-1}+p^{-1!}X_{i+1}^{-1!}\cdot \dots \cdot X_n^{-1!},\dots
,p^{-1}X_1^{-m}\cdot \dots \cdot X_i^{-m}+\cr
&+(p^{-1!}X_{i+1}^{-1!}\cdot \dots \cdot X_n^{-1!})\cdot
(X_1^{-(m-1)}\cdot \dots \cdot X_i^{-(m-1)})+\dots +\cr
&+(p^{-(m-1)!}X_{i+1}^{-(m-1)!}\cdot \dots \cdot X_n^{-(m-1)!})\cdot
(X_1^{-1}\cdot \dots \cdot X_i^{-1})+p^{-m!}X_{i+1}^{-m!}\cdot \dots
\cdot X_n^{-m!},\dots )\in D(\LCMo ^i_I(R))\cr }$$ and, similar to
the proof of lemma 2.1, we show that $\ann _R(\alpha )=0$. Assume to
the contrary there is an $f\in \ann _R(\alpha )\setminus \{ 0\} $.
Choose $(a_1,\dots ,a_i)$ minimal in
$$\{ (a_1^\prime, \dots ,a_i^\prime )\vert \hbox{there exists }a_{i+1}^\prime
,\dots a_n^\prime \hbox{ such that }(a_1^\prime ,a_n^\prime )\in
\Supp (f)\}\ \ .$$ Like before we may assume $a_1=\dots =a_i=:a$.
Choose $h_1,\dots ,h_i\in R$ and $g\in C[[X_{i+1},\dots
,X_n]]\setminus \{ 0\} $ such that
$$f=X_1^{a+1}\cdot h_1+\dots +X_i^{a+1}\cdot h_i+X_1^a\cdot \dots \cdot
X_i^a\cdot g$$ $\alpha \cdot f=0$ implies, for all $m\in \bf N
\setminus \{ 0\} $,
$$0=(X_1^{a+1}h_1+\dots +X_i^{a+1}h_i+X_1^a\cdot \dots \cdot X_i^a\cdot
g)\cdot (p^{-1}X_1^{-m}\cdot \dots \cdot X_i^{-m}+\dots
+p^{-m!}X_{i+1}^{-m!}\cdot \dots \cdot X_n^{-m!})=$$
$$=(X_1^{a+1}h_1+\dots
+X_i^{a+1}h_i)\cdot [p^{-1}X_1^{-m}\cdot \dots \cdot X_i^{-m}+\dots
+(p^{-(m-a-1)!}\cdot X_{i+1}^{-(m-a-1)!}\cdot \dots \cdot
X_n^{-(m-a-1)!})\cdot $$
$$\cdot (X_1^{-(a+1)}\cdot \dots \cdot X_i^{-(a+1)})]+g\cdot
(p^{-1}X_1^{-(m-a)}\cdot \dots \cdot X_i^{-(m-a)}+\dots
+p^{-(m-a)!}X_{i+1}^{-(m-a)!}\cdot \dots \cdot X_n^{-(m-a)!})\ \ .$$
Now, let $(b_{i+1},\dots ,b_n)$ be minimal in $\Supp (g)$ and $c\in
C$ be the coefficient of $g$ in front of $X_{i+1}^{b_{i+1}}\cdot
\dots \cdot X_n^{b_n}$. In $C_p/C$ we have $c\cdot p^{-(m-a)!}\neq
0$ for all $m>>0$. So, like before, we must have
$$(m-a)!-b_{i+1}\leq (m-a-1)!$$
$$(m-a)!-b_n\leq (m-a-1)!$$
for all $m>>0$, which leads to a contradiction again.
\bigskip
{\bf 2.3 Lemma}
\par
Let $p$ be a prime number, $C$ a complete $p$-ring, $n\in \bf N$,
$R=C[[X_1,\dots ,X_n]]$, $i\in \{ 0,\dots ,n\} $, $I:=(p,X_1,\dots
,X_i)R$ and $\goth m:=(p,X_1,\dots ,X_n)R$. Then
$$\{ 0\} \in \Ass _R(D(\LCMo^{i+1}_I(R)))$$
holds. \par
Proof:
\par
1. Case: $i=n$: In this case we have $\LCMo ^{i+1}_I(R)=\InjH
_R(R/\goth m)$
   and hence $D(\LCMo ^{i+1}_I(R))=R$.
\par
2. Case: $i<n$: Similar to the situation in the proof of lemma we
have
$$D(\LCMo ^{i+1}_I(R))=\vtop{\baselineskip=1pt
\lineskiplimit=0pt \lineskip=1pt\hbox{lim} \hbox{$\longleftarrow
$}}(\InjH _{R/(p,X_1,\dots ,X_i)R}(R/\goth m)\buildrel p\cdot
X_1\cdot \dots \cdot X_i\over \longleftarrow \InjH
_{R/(p^2,X_1^2,\dots ,X_i^2)R}(R/\goth m)\buildrel p\cdot X_1\cdot
\dots \cdot X_i\over \longleftarrow \dots )$$
$$\InjH _R(R/\goth m)=(C_p/C)[X_1^{-1},\dots ,X_n^{-1}]$$
and we define
$$\eqalign {\alpha :=&(p^{-1},p^{-2}X_1^{-1}\cdot \dots \cdot
X_i^{-1}+p^{-2}X_{i+1}^{-1!}\cdot \dots \cdot X_n^{-1!},\dots
,p^{-(m+1)}X_1^{-m}\cdot \dots \cdot X_i^{-m}+\cr
&+p^{-(m+1)}X_{i+1}^{-1!}\cdot \dots \cdot X_n^{-1!}\cdot
X_1^{-(m+1)}\cdot \dots \cdot X_i^{-(m-1)}+\dots +\cr
&+p^{-(m+1)}X_{i+1}^{-(m-1)!}\cdot \dots \cdot X_n^{-(m-1)!}\cdot
X_1^{-1}\cdot \dots \cdot X_i^{-1}+p^{-(m+1)}X_{i+1}^{-m!}\cdot
\dots \cdot X_n^{-m!},\dots )\ \ .\cr }$$ Again we state $\ann
_R(\alpha )$ and assume, to the contrary, that there exists an $f\in
\ann _R(\alpha )\setminus \{ 0\} $, choose $(a_1,\dots ,a_i)$
minimal in
$$\{ (a_1^\prime ,\dots ,a_i^\prime )\vert \hbox{There exist }a_{i+1}^\prime
,\dots ,a_n^\prime \hbox{ such that }(a_1^\prime ,\dots ,a_n^\prime
)\in \Supp (f)\} $$ may assume $a_1=\dots =a_i=:a$ and choose
$h_1,\dots ,h_i\in R$, $g\in C[[X_{i+1},\dots ,X_n]]$ such that
$$f=X_1^{a+1}h_1+\dots +X_i^{a+1}h_i+(X_1^a\cdot \dots \cdot a_i^a)\cdot g\ \ .$$
This means, for all $m\in \bf N$,
$$0=(X_1^{a+1}h_1+\dots +X_i^{a+1}h_i)[p^{-(m+1)}X_1^{-m}\cdot \dots \cdot
X_i^{-m}+\dots +(p^{-(m+1)}X_{i+1}^{-(m-a-1)!}\cdot \dots \cdot
X_n^{-(m-a-1)!})\cdot $$
$$\cdot (X_1^{-(a+1)}\cdot \dots \cdot X_i^{-(a+1)})]+g\cdot
(p^{-(m+1)}X_1^{-(m-a)}\cdot \dots \cdot X_i^{-(m-a)}+\dots
+p^{-(m+1)}X_{i+1}^{-(m-a)!}\cdot \dots \cdot X_n^{-(m-a)!})\ \ .$$
Choose $(b_{i+1},\dots ,b_n)$ minimal in $\Supp (g)$ and let $c\in
C$ be the coefficient of $g$ in front of $X_{i+1}^{b_{i+1}}\cdot
\dots \cdot X_n^{b_n}$. In $C_p/C$ we have $g\cdot p^{-(m+1)}\neq 0$
for all $m>>0$, and so we must have for all $m>>0$
$$(m-a)!-b_{i+1}\leq (m-a-1)!$$
$$(m-a)!-b_n\leq (m-a-1)!$$
which leads to a contradiction, proving the lemma.
\bigskip
Now we are ready to prove that certain prime ideals are associated
to $D(\LCMo ^i_{(x_1,\dots ,x_i)R}(R))$ in a more general situation
($R$ does not have to be regular). This is done essentially by using
various base-change arguments and lemmas 2.1 -- 2.3:
\bigskip
{\bf 2.4 Theorem}
\par
Let $(R,\goth m)$ be a noetherian local ring, $i\geq 1$ and
$x_1,\dots ,x_i$ a sequence on $R$. Then
$$\Ass _R(D(\LCMo ^i_{(x_1,\dots ,x_i)R}(R)))\subseteq \{ \goth p\in \Spec (R)\vert \LCMo ^i_{(x_1,\dots
,x_i)R}(R/\goth p)\neq 0\} $$ holds. If $R$ is equicharacteristic,
$$\{ \goth p\in \Spec (R)\vert x_1,\dots ,x_i\hbox { is part of a
system of parameters for }R/\goth p\} \subseteq \Ass _R(D(\LCMo
^i_{(x_1,\dots ,x_i)R}(R)))$$ holds. \par
Proof:
\par
The first inclusion was shown in remark 1.2.1. For the second
inclusion let $\goth p\in \Spec (R)$ and $x_{i+1},\dots ,x_n\in R$
such that $x_1,\dots ,x_n$ (more precisely: their images in $R/\goth
p$) form a system of parameters for $R/\goth p$; then $n=\dim
(R/\goth p)$. $x_1,\dots ,x_n$ also form a system of parameters in
$\hat R/\goth p\hat R$. Choose $\goth q\in \Spec (\hat R)$ with
$\dim (\hat R/\goth q)=\dim (R/\goth p)$. This implies $\goth q\in
\Min (\hat R)$ and $\goth q\cap R=\goth p$. Because of $\dim (\hat
R/\goth q)=\dim (R/\goth p)$ the elements $x_1,\dots ,x_n$ form a
system of parameters of $\hat R/\goth q$. It is sufficient to show
$ \goth q\in \Ass _{\hat R}(D(\LCMo ^i_{(x_1,\dots ,x_i)\hat R}(\hat
R)))$. Namely, as
$$\eqalign {D(\LCMo ^i_{(x_1,\dots ,x_i)\hat R}(\hat R))&=\Hom _{\hat R}(\LCMo
^i_{(x_1,\dots ,x_i)\hat R}(\hat R),\InjH _{\hat R}(\hat R/\goth
m\hat R))\cr &=\Hom _{\hat R}(\LCMo ^i_{(x_1,\dots ,x_i)\hat R}(\hat
R),\InjH _R(R/\goth m))\cr &=\Hom _{\hat R}(\LCMo ^i_{(x_1,\dots
,x_i)R}(R)\otimes _R\hat R,\InjH _R(R/\goth m))\cr &=\Hom _R(\LCMo
^i_{(x_1,\dots ,x_i)R}(R),\InjH _R(R/\goth m))\cr &=D(\LCMo
^i_{(x_1,\dots ,x_i)R}(R))\ \ ,\cr }$$ every monomorphism $\hat
R/\goth q\to D(\LCMo ^i_{(x_1,\dots ,x_i)\hat R}(\hat R))$ induces a
monomorphism
$$R/\goth p\buildrel \hbox {kan.}\over \to \hat R/\goth q\to D(\LCMo ^i_{(x_1,\dots ,x_i)\hat
R}(\hat R))=D(\LCMo ^i_{(x_1,\dots ,x_i)R}(R))\ \ .$$ This means we
may assume that $R$ is complete.
\par
We have to show that the zero ideal of $R/\goth p$ is associated to
$$\Hom _R(R/\goth p,D(\LCMo^i_{x_1,\dots ,x_i)R}(R))= D(\LCMo
^i_{(x_1,\dots ,x_i)R/\goth p}(R/\goth p))$$ (this equality was
shown in the proof of the first inclusion). Replacing $R$ by
$R/\goth p$ we may assume that $R$ is a domain and $\goth p$ is the
zero ideal in $R$. Let $k\subseteq R$ denote a coefficient field.
$$R_0:=k[[x_1,\dots ,x_n]]\subseteq R$$
is an $n$-dimensional regular local subring of $R$, over which $R$
is module-finite. Let $\goth m_0$ denote the maximal ideal of $R_0$.
The $R$-Modul $\Hom _{R_0}(R,\InjH _{R_0}(R_0/\goth m_0))$ is
isomorphic to $\InjH _R(R/\goth m)$. We have
$$\eqalign {D(\LCMo ^i_{(x_1,\dots ,x_i)R}(R))&=
\Hom _R(\LCMo ^i_{(x_1,\dots ,x_i)R}(R),\InjH _R(R/\goth m))\cr
&=\Hom _R(\LCMo ^i_{(x_1,\dots ,x_i)R_0}(R_0)\otimes _{R_0}R,\InjH
_R(R/\goth m))\cr &=\Hom _{R_0}(\LCMo ^i_{(x_1,\dots ,x_i)R_0}(R_0),
\Hom _{R_0}(R,\InjH _{R_0}(R_0/\goth m_0)))\cr &=\Hom _{R_0}(R,\Hom
_{R_0}(\LCMo ^i_{(x_1,\dots ,x_i)R_0} (R_0),\InjH _{R_0}(R_0/\goth
m_0)))\cr &=\Hom _{R_0}(R,D(\LCMo ^i_{(x_1,\dots ,x_i)R_0}(R_0)))\ \
.\cr }$$ By lemma 2.1 there exists a monomorphism $R_0\to D(\LCMo
^i_{(x_1,\dots ,x_i)R_0}(R_0))$; so we get a monomorphism
$$\Hom _{R_0}(R,R_0)\to \Hom _{R_0}(R,D(\LCMo ^i
_{(x_1,\dots ,x_i)R_0}(R_0)))=D(\LCMo ^i_{(x_1,\dots ,x_i)R}(R))\ \
.$$ $R$ is a domain and module-finite over $R_0$, and thus $\{ 0\}
\in \Supp _R(\Hom _{R_0}(R,R_0))$; the statement now follows.
\bigskip
Again, there are versions for the case of mixed characteristic:
\bigskip
{\bf 2.5 Theorem}
\par
Let $(R,\goth m)$ be a noetherian local ring of mixed
characteristic, $p=\char (R/\goth m)$, $i\geq 0$ and $x_1,\dots
,x_i\in R$. Then
$$\{ \goth p\in \Spec (R)\vert p,x_1,\dots ,x_i\hbox { is part of a system of
parameters for }R/\goth p\} \subseteq \Ass _R(D(\LCMo
^{i+1}_{(p,x_1,\dots ,x_i)R}(R)))\ \ .$$ In case $i\geq 1$, we have
in addition
$$\{ \goth p\in \Spec (R)\vert p,x_1,\dots ,x_i\hbox { is part of a system of
parameters for }R/\goth p\} \subseteq \Ass _R(D(\LCMo ^i_{(x_1,\dots
,x_i)R}(R)))\ \ .$$ Theorem 2.5 is proved in a similar way like
Theorem 2.4, using lemmas 2.2 and 2.3 instead of lemma 2.1.
\bigskip
In the case $i=1$ the results proven so far are sufficient to
completely compute the set of associated primes:
\bigskip
{\bf 2.6 Corollary}
\par
Let $(R,\goth m)$ be a noetherian local equicharacteristic ring and
$x\in R$. Then
$$\Ass _R(D(\LCMo ^1_{xR}(R)))=\Spec (R)\setminus \goth V(x)$$
holds. In particular, this set is infinite in general.
\bigskip
{\bf 2.7 Remarks}
\par
(i) If one has $\Ass _R(D(\LCMo ^i _{(x_1,\dots
,x_i)R}(R)))=\emptyset $ in the situation of the theorem, it follows
that $\LCMo ^i_{(x_1,\dots ,x_i)R}(R)=0$ and also $\LCMo ^i
_{(x_1,\dots ,x_i)R}(R/\goth p)=0$ for every $\goth p \in \Spec (R)$
(by a well-known theorem), i. e. in this case conjecture (*) holds.
\smallskip
(ii) The second inclusion of theorem 2.4 is not an equality in
general: For a counterexample let $k$ be a field,
$R=k[[y_1,y_2,y_3,y_4]]$ and define $x_1=y_1y_3$, $x_2=y_2y_4$,
$x_3=y_1y_4+y_2y_3$. $x_1,x_2,x_3$ is not part of a system of
parameters for $R$, but we have
$$\sqrt
{(x_1,x_2,x_3)R}=(y_1,y_2)R\cap (y_3,y_4)R$$ and so a Mayer-Vietoris
sequence argument (see, e. g. [BS, 3.2.3] for a reference on the
Mayer-Vietoris sequence) shows
$$\InjH _R(k)=\LCMo ^3_{(y_1,y_2)R\cap
(y_3,y_4)R}(R)=\LCMo ^3_{(x_1,x_2,x_3)R}(R)$$
and so $D(\LCMo
^3_{(x_1,x_2,x_3)R}(R))=R$. Thus $\{ 0\} \in \Ass _R(D(\LCMo
^3_{(x_1,x_2,x_3)R}(R)))$.
\smallskip
(iii) In the situation of theorem 2.4, set
$$Z_1:=\{ \goth p\in \Spec (R)\vert \LCMo ^i_{(x_1,\dots ,x_i)R}(R/\goth
p)\neq 0\} $$ and
$$Z_2:=\{ \goth p\in \Spec (R)\vert x_1,\dots ,x_i\hbox { is part of a
system of parameters for }R/\goth p\}\ \ .$$ Then $Z_1$ is stable
under generalization (this follows e. g. from the following
well-known fact: If $I$ is an ideal of a noetherian domain $R$ such
that $0=\LCMo ^l_I(R)=\LCMo ^{l+1}_I(R)=\dots $ holds for some fixed
$l\in \Naturalsign $, then $0=\LCMo ^l_I(M)=\LCMo ^{l+1}_I(M)=\dots
$ holds for every $R$-module $M$). \par But note that, in general,
$Z_2\subseteq \Spec (R)$ is not stable under generalization, even
not if $R$ is regular; namely, for an example where $Z_2$ is not
stable under generalization, let $R=k[[x_1,x_2,x_3,x_4]]$ be a
formal power series algebra in four variables over a field $k$, set
$$\goth p_0=(x_1x_4+x_2x_3)R\hbox { and }\goth p=(x_3,x_4)R\ \ .$$
Then $x_1,x_2$ is a system of parameters for $R/\goth p$, but is not
a part for $R/\goth p_0$ (because $x_1x_4+x_2x_3$ is contained in
the ideal $(x_1,x_2)$), i. e. we have
$$\goth p_0\subseteq \goth p,\goth p\in Z_2, \goth
p_0\not\in Z_2\ \ .$$ Assume now that $R$ is regular; then, at
least, the following special form of stableness (of $Z_2$) under
generalization holds: Let $\goth p\in \Spec(R)$ such that $x_1,\dots
,x_i$ is part of a system of parameters for $R/\goth p$. Then
$x_1,\dots ,x_i$ is part of a system of parameters for $R$, i. e.
one has the implication
$$Z_2\neq \emptyset \Longrightarrow \{ 0\} \in Z_2\ \ .$$
This follows from the so-called height-formula which holds for
regular local rings and which says (we apply it to the ideal
$(x_1,\dots ,x_i)R+\goth p\subseteq R$):
$$\height ((x_1,\dots ,x_i)R+\goth p)\leq \height ((x_1,\dots
,x_i)R)+\height (\goth p)\leq i+\height (\goth p)\ \ .$$ But,
because of our assumption $\goth p\in Z_2$, we must have
$$\height ((x_1,\dots ,x_i)R+\goth p)= i+\height (\goth p)$$
and, therefore, $\height ((x_1,\dots ,x_i)R)=i$. \vfil \eject {\bf 3
Associated primes -- the characteristic-free approach}
\bigskip
In this section we investigate associated prime ideals of Matlis
duals $D(\LCMo ^i_I(M))$ of local cohomology modules ($R$ is local,
of course); there are two subsections: In the first one, we prove
characteristic-free versions of some results on the set of
associated primes of such a module; out methods here are different
to the ones used in section 2. Some results of this section can be
found in [HS1]. In the second part of this section, we concentrate
on the case $M=R$, $i=\dim (R)-1$, theorems 3.2.6 and 3.2.7 (where
we actually compute the set of associated primes of $D(\LCMo ^{\dim
(R)-1}_I(R))$) contain the main results of this second subsection.
\bigskip
\bigskip
{\bf 3.1 Characteristic-free versions of some results}
\bigskip
The following lemma is crucial for this subsection: \bigskip {\bf
3.1.1 Lemma}
\par
Let $R$ be a ring, $x,y\in R$ and $U$ an $R$-submodule of $R_x$ such
that $\im \iota _x\subseteq U$, where $\iota _x:R\to R_x$ is the
canonical map. Let $S:= \im \iota _y\subseteq R_y$. There exists an
$R$-epimorphism
$$R_x/U\to R_{xy}/(S_x+U_y).$$
Proof: \par
Let $V:=S_x+U_y\subseteq R_{xy}$ and let $(b_1,b_2,\dots
)\in R^{\Naturalsign ^+}$ be an infinite sequence. For $i\in
\Naturalsign $ we set
$$\rho _i:=\sum _{j=1}^i{b_j\over x^{i-j+1}y^j}+V\in R_{xy}/V\ \ (i\in
\Naturalsign ).$$ We calculate
$$\eqalign {x\rho _{i+1}-\rho _i&=(\sum _{j=1}^{i+1}{xb_j\over
x^{i-j+2}y^j}+V)-(\sum _{j=1}^i{b_j\over x^{i-j+1}y^j}+V)\cr
&={b_{i+1}\over y^{i+1}}+V\cr &=0,\cr }$$ because
$${b_{i+1}\over
y^{i+1}}\in (\im \iota _x)_y\subseteq U_y\subseteq V\ \ .$$
Thus we
have $x\rho _{i+1}=\rho _i$ for all $i\in \Naturalsign $ and so we
get a map $\varphi :R_x\to R_{xy}/V$ given by
$${r\over x^i}\mapsto r\rho _i\ \ (r\in R,i\in \Naturalsign ).$$
It is easy to see that $\varphi $ is $R$-linear. Let $u\in U$ be
arbitrary. There are $r\in R$ and $i\in \Naturalsign $ such that
$u={r\over x^i}$. We have
$$\varphi (u)=r\rho _i=\sum _{j=1}^i{rb_j\over x^{i-j+1}}+V=u\sum
_{j=1}^i{x^{j-1}b_j\over y^j}+V=0,$$ because $$u\sum
_{j=1}^i{x^{j-1}b_j\over y^j}\in U_y\subseteq V\ \ .$$
This implies
$U\subseteq \ker (\varphi )$ and hence we get an induced
$R$-homomorphism $f:R_x/U\to R_{xy}/V$. The set $\{ {1\over
x^i}+U\vert i\in \Naturalsign ^+\} $ is a generating set for $R_x/U$
and so we have
$$f\hbox { is surjective }\iff \varphi \hbox{ is surjective }\iff \{ \rho _1,\rho
_2, \dots \} \hbox { generates } R_{xy}/V.$$ The set $\{ {1\over
x^iy^j}+V\vert i,j\in \Naturalsign ^+\} $ generates $R_{xy}/V$. For
$i\in \Naturalsign ^+$ we set
$${\cal B}_i:=\pmatrix {&b_1&b_2&b_3&\ldots &b_i\cr &b_2&b_3&b_4&\ldots
&b_{i+1}\cr &\vdots &\vdots &\vdots & &\vdots \cr
&b_i&b_{i+1}&b_{i+2}&\ldots &b_{2i-1}\cr }$$ Then we have for $i\in
\Naturalsign ^+$:
$$(\rho _i,y\rho _{i+1},\dots ,y^{i-1}\rho _{2i-1})^T={\cal B}_i({1\over
x^iy}+V,{1\over x^{i-1}y^2}+V,\dots ,{1\over xy^i}+V)^T.$$ If we
choose $b_1,b_2,\dots \in R$ in such a way that $\det {\cal B}_i\in
R^*$ for all $i\in \Naturalsign ^+$ (which is possible, $\cal B$
consists only of ones and zeroes), then $\{ \rho _1,\rho _2, \dots
\}$ generates $R_{xy}/V$.
\bigskip
From now on we assume that $R$ is noetherian, we can use \v{C}ech
cohomology to compute local cohomology. Thus, lemma 3.1.1 implies:
\bigskip
{\bf 3.1.2 Theorem}
\par
Let $R$ be a noetherian ring, $M$ an $R$-module, $m\in \Naturalsign
^+, n\in \Naturalsign ,x_1,\dots ,x_m,y_1,\dots ,y_n\in R$. Then
there exists an $R$-epimorphism
$$\LCMo ^m_{(x_1,\dots ,x_m)R}(M)\to \LCMo ^{m+n}_{(x_1,\dots ,x_m,y_1,\dots
,y_n)R}(M).$$ Proof: \par
Obviously it suffices to prove the
statement for the case $M=R$. Using \v{C}ech cohomology to compute
both local cohomology modules the statement follows immediately from
lemma by induction on $n$.
\bigskip
By dualizing the surjection from the preceding theorem we get an
injection. But, then, the set of associated prime ideals of the
right-hand side is contained in the set of associated prime ideals
of the left-hand side. This is the basic idea in the proof of
statement (ii) in the following theorem (the same is true for (iii),
(iv) and (v), as these statements follow from (ii), see the proof
below for details):
\bigskip
{\bf 3.1.3 Theorem}
\par
Let $(R,\goth m)$ be a noetherian local ring, $m\in \Naturalsign ^+,
x_1,\dots ,x_m\in R$ and $M$ a finitely generated $R$-module. Then
the following statements hold:
\par
(i) $\dim (M/\goth pM)\geq m$ for every $\goth p\in \Ass _R(D(\LCMo
^m_{(x_1,\dots ,x_m)R}(M)))$.
\par
(ii) $\{ \goth p\in \Supp _R(M)\vert x_1,\dots ,x_m\hbox { is part
of a system of parameters of }R/\goth p\} \subseteq \Ass _R(D(\LCMo
^m_{(x_1,\dots ,x_m)R}(M)))$.
\par
(iii) $\Ass _R(D(\LCMo ^1_{xR}(R)))=\Spec (R)\setminus {\frak V}(x)$
for every $x\in R$.
\par
(iv) If $x_1,\dots ,x_m$ is part of a system of parameters of $M$,
we have $\Assh (M)\subseteq \Ass _R(D(\LCMo ^m_{(x_1,\dots
,x_m)R}(M)))$; furthermore, if $m=\dim (M)$, equality holds: $\Assh
(M)=\Ass _R(D(\LCMo ^{\dim (M)}_\goth m(M)))$ (note that, by
definition, $\Assh (M)$ consists of the associated prime ideals of
$M$ of highest dimension).
\par
(v) If $R$ is complete, $\goth p\in \Supp _R(M)$ and $\dim (R/\goth
p)=m$, the equivalence
$$\goth p\in \Ass _R(D(\LCMo ^m_{(x_1,\dots ,x_m)R}(M)))\iff x_1,\dots
,x_m\hbox { is a system of parameters of }R/\goth p$$ holds.
\par
Proof: \par
We set $I:=(x_1,\dots ,x_m)R$.
\par
(i) Let $\goth p\in \Ass _R(D(\LCMo ^m_I(M)))$. We conclude
$$0\neq \Hom _R(R/\goth p,D(\LCMo ^m_I(M)))=D(\LCMo ^m_I(M)\otimes _R(R/\goth
p))=D(\LCMo ^m_I(M/\goth pM)).$$ Thus we have $\LCMo ^m_I(M/\goth
pM)\neq 0$ and statement (i) follows (note that it follows also from
remark 1.2.1).
\par
(ii) Let $\goth p\in \Supp _R(M)$ such that $x_1,\dots ,x_m$ is part
of a system of parameters of $R/\goth p$. By completing $x_1,\dots
,x_m$ to a system of parameters of $M/\goth pM$ and using theorem
3.1.2, we may assume that $x_1,\dots ,x_m$ is a system of parameters
of $M/\goth pM$. So we have $\dim M/\goth pM$=$\dim (R/\goth p)=m$.
Therefore we get
$$\eqalign {\Hom _R(R/\goth p,D(\LCMo ^m_I(M)))&=D(\LCMo ^m_I(M/\goth pM))\cr
&=D(\LCMo ^m_\goth m(M/\goth pM)))\cr &\neq 0.\cr }$$ On the other
hand we have $\Hom _R(R/\goth q,D(\LCMo ^m_I(M)))=0$ for every prime
ideal $\goth q$ of $R$ containing $\goth p$ properly, by (i);
statement (ii) follows.
\par
(iii) Using (ii), it remains to show that $x\not\in \goth p$ holds
for every $\goth p\in \Ass _R(D(\LCMo ^1_{xR}(R)))$. As we have seen
above, our hypothesis implies $\LCMo ^1_{xR}(R/\goth p)\neq 0$. So
we must have $x\not\in \goth p$.
\par
(iv) The first statement follows from (ii) (note that, for every
$\goth p\in \Assh (M)$, $x_1,\dots ,x_m$ is part of a system of
parameters of $R/\goth p$, too) and then the second statement from
(i).
\par
(v) Let $\goth p\in \Supp _R(M)$ such that $\goth p\in \Ass
_R(D(\LCMo ^m_I(M)))$. We have to show that $x_1,\dots ,x_m$ is a
system of parameters of $M/\goth pM$: $\LCMo ^m_I(M/\goth pM)\neq 0$
implies $\LCMo ^m_I(R/\goth p)\neq 0$. As $R$ and hence $R/\goth p$
are complete we may conclude from Hartshorne-Lichtenbaum vanishing
(see, e. g. , [BS, 8.2.1] or theorem 6.1.4 for a reference on
Hartshorne-Lichtenbaum vanishing) that $\dim (R/(I+\goth p))=0$, i.
e. $x_1,\dots ,x_m$ is a system of parameters of $R/\goth p$.
\bigskip
\bigskip
{\bf 3.2 On the set $\Ass _R(D(\LCMo ^{\dim (R)-1}_I(R)))$}
\bigskip
We prove a series of lemmas which we will need for the main results
3.2.6 and 3.2.7.
\bigskip
{\bf 3.2.1 Lemma}
\par
Let $(S,\goth m)$ be a noetherian local complete Gorenstein ring of
dimension $n+1$ ($\geq 1$) and $\goth P\subseteq S$ a prime ideal of
height $n$. Then
$$D(\LCMo ^n_\goth P(S))=\widehat {S_\goth P}/S$$
holds canonically.
\par
Proof: \par Local duality over the Gorenstein ring $S$ shows that
there are natural isomorphisms
$$D(\LCMo ^n_\goth P(S))=D(\vtop{\baselineskip=1pt \lineskiplimit=0pt \lineskip=1pt\hbox{lim}
\hbox{$\longrightarrow $} \hbox{$^{^{l\in \bf N}}$}}\Ext
^n_S(S/\goth P^l,S))=\vtop{\baselineskip=1pt \lineskiplimit=0pt
\lineskip=1pt\hbox{lim} \hbox{$\longleftarrow $} \hbox{$^{^{l\in \bf
N}}$}}\LCMo ^1_\goth m(S/\goth P^l)\ \ .$$ Take $y\in \goth
m\setminus \goth P$. Now, $\sqrt {y(S/\goth P^l)}=\goth m (S/\goth
P^l)$ implies
$$\LCMo ^1_\goth m(S/\goth
P^l)=\LCMo ^1_{y(S/\goth P^l)}(S/\goth P^l)=(S_y/\goth
P^lS_y)/(S/\goth P^l)=(S_\goth P/\goth P^lS_\goth P)/(S/\goth P^l)$$
and the statement follows by observing that the maps
$$(S_\goth
P/\goth P^{l+1}S_\goth P)/(S/\goth P^{l+1})\to (S_\goth P/\goth
P^lS_\goth P)/(S/\goth P^l)$$
which we get from this, are the natural
ones.
\bigskip
{\bf 3.2.2 Lemma}
\par
Let $(R,\goth m)$ be a noetherian local complete domain and
$I\subseteq R$ a prime ideal such that $\dim (R/I)=1$. Then there
exist a noetherian local complete regular ring $S$, a local
homomorphism $S\buildrel \rho \over \to R$ and a prime ideal $\goth
Q\subseteq S$ such that $R$ is finite as an $S$-module and such that
$$\height (\ker (\rho ))=1, \dim (S/\goth Q)=1, \sqrt {\goth QR}=I, \ker
(\rho )\subseteq \goth Q$$ hold.
\par
Proof: \par Either $R$ contains a field $k$ or, if not, a
coefficient ring $(V,tV)$; choose $y_1,\dots ,y_n,y\in R$ such
that
$$I=\sqrt {(y_1,\dots,y_n)R}$$
and $y_1,\dots ,y_{n-1},y$ is a
system of parameters of $R$ ($\dim (R)=n$); in the case of mixed
characteristic we may take $y_1:=t$ if $t\in I$ and $y:=t$ if
$t\not\in I$. Now we define a subring of $R$:
$$R_0:=k[[y,y_1,\dots ,y_n]]$$
(if $R$ contains a field) resp.
$$R_0:=V[[y,y_2,\dots ,y_n]]$$
(if $R$ contains no field and $t\in I$) resp.
$$R_0:=V[[y_1,\dots ,y_n]]$$
(if $R$ contains no field and $t\not\in I$). Furthermore we define a
power series ring
$$S:=k[[Y,Y_1,\dots ,Y_n]]\hbox { resp. } V[[Y,Y_2,\dots ,Y_n]]\hbox { resp. }
V[[Y_1,\dots ,Y_n]]$$ (all capital letters denote variables) and it
is clear how to define a surjection $S\buildrel \rho \over \to
R_0(\subseteq R)$. We set
$$\goth Q:=(Y_1,\dots,Y_n)S\hbox { resp. } (t,Y_2,\dots ,Y_n)S\hbox { resp. }
(Y_1,\dots ,Y_n)S,$$ where in all three cases we have $\ker (\rho
)\subseteq \goth Q$ because of $y\not\in I$. All other statements
are obvious now.
\bigskip
{\bf 3.2.3 Lemma}
\par
Let $R$ be a noetherian ring.
\par
(i) Let $\goth P$ be a prime ideal of $R$ which is not maximal. Then
the equivalence
$$R_\goth P=\widehat {R_\goth P}\iff \goth P\hbox { is minimal in }\Spec(R)$$
holds.
\par
(ii) Assume that $R$ is local (and noetherian) and that all prime
ideals associated to $R$ are minimal in $\Spec (R)$. Then $\Ass
_R(\hat R/R)\subseteq \Ass (R)$ holds. In particular if $R$ is a
non-complete (local) domain (i. e. if $R\subsetneq \hat R$),
$$\Ass _R(\hat R/R)=\{ 0\} $$
holds.
\par
Proof: \par
(i) The implication $\Leftarrow $ is clear as every
zero-dimensional local noetherian ring is complete. We assume there
exists a prime ideal $P$ of $R$ which is neither minimal nor maximal
in $\Spec(R)$ and such that $R_P=\widehat {R_P}$. $PR_P$ is not
minimal in $\Spec (R)$. Choose $Q,Q^\prime \in \Spec (R)$ such that
$Q^\prime \subsetneq P\subsetneq Q$ and such that $\dim
(R_Q/PR_Q)=1$. We set ${\frak R}:=R_Q/Q^\prime R_Q$ and ${\frak
P}:=PR_Q/Q^\prime R_Q\in \Spec ({\frak R})$ and we get
$${\frak R_P}=R_P/Q^\prime R_P=\widehat {R_P}/Q^\prime \widehat {R_P}=\widehat
{R_P/Q^\prime R_P}=\widehat {\frak R_P}.$$ So we may assume that $R$
is a local domain  and $\dim (R/P)=1$.
\par
Take $y\in \goth m\setminus P$. Assume that for some $n\in
\Naturalsign $
$$P^{(n)}\subseteq P^{(n+1)}+yR$$
holds ($P^{(n)}:=P^nR_P\cap R$ is a $P$-primary ideal of $R$ such
that $P^{(n)}R_P=P^nR_P$, $P^{(n)}$ is the so-called $n$-th symbolic
power of $P$). It would follow that
$$P^{(n)}=P^{(n)}\cap (P^{(n+1)}+yR)=P^{(n+1)}+(P^{(n)}\cap
yR)=P^{(n+1)}+yP^{(n)}$$ (the last equality follows, because $\goth
P^{(n)}$ is $\goth P$-primary) and then $P^{(n)}=P^{(n+1)}$, by the
lemma of Nakayama (see, e. g., [Ma, Theorem 2.2] for the lemma of
Nakayama). Again by the lemma of Nakayama, this would imply
$P^nR_P=0$ and so $PR_P$ would be minimal in $\Spec (R_P)$. We
conclude that for every $n\in \Naturalsign $
$$P^{(n)}\not\subseteq P^{(n+1)}+yR$$
holds. For every $n\in \Naturalsign $ we choose $x_n\in
P^{(n)}\setminus (P^{(n+1)}+yR)$ and define (for every $n\in
\Naturalsign ^+$)
$$\xi _n:=\sum _{i=0}^{n-1}{x_i\over y^{(i+1)^2}}\in R_P.$$
Because of
$$\xi  _{n+1}-\xi _n={x_n\over y^{(n+1)^2}}\in P^nR_P$$
(for
every $n$), we have
$$(\xi _n+P^nR_P)_{n\in \Naturalsign ^+}\in \widehat {R_P}=R_P.$$
Therefore, there exists $\xi \in R_P$ such that
$$(\xi +P^nR_P)_{n\in \Naturalsign ^+}=(\xi _n+P^nR_P)_{n\in \Naturalsign
^+},$$ i. e.
$$\xi -\xi _n\in P^nR_P$$ holds for every $n\in \Naturalsign ^+$.
\par
Write $\xi ={a\over b}$, where $a\in R,b\in R\setminus P$. The ideal
$P+bR$ of $R$ is either $R$ or $\goth m$-primary, so there exist
$p\in \Naturalsign ^+$ and $c\in R$ such that $y^p-bc\in P$; it
follows that
$$y^{pn}-bc_n\in P^n,$$
where
$$c_n:=b^{-1}(y^{pn}-(y^p-bc)^n)\in R$$
(note that $y^{pn}-(y^p-bc)^n$ is divisible by $b$ in $R$) and we
conclude that
$$\xi -{ac_n\over y^{pn}}={ay^{pn}-abc_n\over by^{pn}}={a(y^p-bc)^n\over
by^{pn}}\in P^nR_P$$ for every $n\in \Naturalsign ^+$. We get
$$\xi _n-{ac_n\over y^{pn}}=\xi -{ac_n\over y^{pn}}-(\xi -\xi _n)\in P^nR_P$$
for every $n\in \Naturalsign ^+$. From this we get (for $n>p$) after
multiplication by $y^{n^2}$ that
$$\sum _{i=0}^{n-1}x_iy^{n^2-(i+1)^2}-ac_ny^{n(n-p)}\in P^{(n)}$$
and in particular $x_{n-1}\in P^{(n)}+yR$ which is a contradiction.
\par
(ii) We only have to prove the first statement, the second one
follows from it immediately; Let $P$ be an arbitrary element of
$\Spec (R)\setminus \Ass (R)$; We conclude $\Hom _R(R/P,R)=0$ and
hence also $\Hom _R(R/P,\hat R)=0$ (because $P$ contains an element
which operates injectively on $R$ and $\hat R$ is flat over $R$).
Thus the short exact sequence
$$0\to R\buildrel \subseteq \over \to \hat R\to \hat R/R\to 0$$
induces an exact sequence
$$0\to \Hom _R(R/P,\hat R/R)\to \Ext ^1_R(R/P,R)\buildrel \varphi \over \to
\Ext ^1_R(R/P,\hat R).$$ By our hypothesis there exists $x\in P$
such that $x\not\in Q$ for all $Q\in \Ass (R)$. We get short exact
sequences
$$0\to R\buildrel x\over \to R\to R/xR\to 0$$
and
$$0\to \hat R\buildrel x\over \to \hat R\to \hat R/x\hat R\to 0.$$
Because of $x\in P$ a commutative diagram with exact rows is
induced:
$$\matrix{&0&\to &\Hom _R(R/P,R/xR)&\to &\Ext ^1_R(R/P,R)&\to &0\cr
&&&\downarrow \psi &&\downarrow \varphi &&\cr &0&\to &\Hom
_R(R/P,\hat R/x\hat R)&\to &\Ext ^1_R(R/P,\hat R)&\to &0.\cr }$$
$\psi $ is injective as $R/xR\subseteq \widehat {R/xR}=\hat R/x\hat
R$. Therefore, $\varphi $ is injective which implies that $\Hom
_R(R/P,\hat R/R)=0$, i. e. $P\not\in \Ass _R(\hat R/R)$.
\bigskip
The following result is a special case of both 3.2.6 and 3.2.7.
\bigskip {\bf
3.2.4 Theorem}
\par
Let $(R,\goth m)$ be a $d$-dimensional local noetherian complete
domain, where $d\geq 2$; let $P$ be a prime ideal of $R$ such that
$\dim (R/P)=1$. Then
$$\{ 0\} \in \Ass _R(D(\LCMo ^{d-1}_P(R)))$$
holds.
\par
Proof: \par We apply lemma 3.2.2, set $R_0:=\im (\rho )$ and
consider the ideal $\goth Q$ from lemma 3.2.2 as an ideal of $R_0$.
Because of lemma 3.2.2, $R_0$ is a complete intersection, in
particular it is Gorenstein. By $\goth m_0$ we denote the maximal
ideal of $R_0$. $R$ is finite over $R_0$ and so we have
$D_{R_0}(R)=\Hom _{R_0}(R,\InjH _{R_0}(R_0/\goth m_0)))=\InjH
_R(R/\goth m)=D_R(R)$, which implies $D_{R_0}(M)=D_R(M)$ for every
$R$-module $M$. On the other hand the functor $\LCMo ^{d-1}_\goth
Q(\_ )$ is right exact by Hartshorne-Lichtenbaum vanishing; in
particular we have
$$D_R(\LCMo ^{d-1}_P(R))=D_{R_0}(\LCMo ^{d-1}_\goth Q(R_0)\otimes
_{R_0}R)=\Hom _{R_0}(R,D_{R_0}(\LCMo ^{d-1}_\goth Q(R_0)))$$ and so
every $R_0$-monomorphism $R_0\to D_{R_0}(\LCMo ^{d-1}_\goth Q(R_0))$
induces an $R$-monomorphism
$$\Hom _{R_0}(R,R_0)\to D_R(\LCMo
^{d-1}_\goth P(R))\ \ .$$
It is easy to see that $\{ 0\} \in \Ass
_{R_0}(\Hom _{R_0}(R,R_0))$ holds (e. g. by localizing) and thus it
suffices to show $\{ 0\} \in \Ass _{R_0}(D_{R_0}(\LCMo ^{d-1}_\goth
Q(R_0)))$, i. e. we may assume $R_0$ is Gorenstein. Now, by lemma
3.2.1, we have a commutative diagram with exact rows:
$$\matrix {&0&\to &R&\buildrel \subseteq \over \to &\widehat {R_P}&\to &D(\LCMo
^{d-1}_P(R))&\to &0\cr &&&\downarrow \subseteq &&\downarrow =
&&&&\cr &0&\to &R_P&\buildrel \subseteq \over \to &\widehat
{R_P}&\to &\widehat {R_P}/R_p&\to &0.\cr }$$ This diagram induces an
epimorphism
$$D(\LCMo ^{d-1}_P(R))\to \widehat {R_P}/R_P.$$
By lemma 3.2.3 (i) we have $\widehat {R_P}/R_P\neq 0$ and it follows
from lemma 3.2.3 (ii) that $(\widehat {R_P}/R_P)\otimes _R(Q(R))\neq
0$. Thus we have $D(\LCMo ^{d-1}_P(R))\otimes _RQ(R)\neq 0$ by the
above epimorphism, which is equivalent to the statement of the
theorem.
\bigskip
{\bf 3.2.5 Lemma}
\par
Let $(R,\goth m)$ be a noetherian local complete domain, $d:=\dim
(R)\geq 1$ and $J\subseteq R$ an ideal of $R$ such that $\dim
(R/J)=1$. Then
$$\Assh (R)\cap \Ass _R(D(\LCMo ^{d-1}_J(R)))=\{ Q\in \Assh
(R)\vert \dim (R/(J+Q))\geq 1\} $$ holds. In particular, if $\LCMo
^d_J(R)=0$,
$$\Assh (D(\LCMo ^{d-1}_J(R)))=\Assh (R)$$
holds. \par Proof: \par The second statement follows from the first
one by Hartshorne-Lichtenbaum vanishing. First, we prove the
statement in the special case where $\dim (R/(J+Q))\geq 1$ for all
$Q\in \Assh (R)$; then we will reduce the general to this special
situation. Now, in the special case, it suffices to show the
inclusion ``$\supseteq $'': By Hartshorne-Lichtenbaum vanishing we
have $\LCMo ^d_J(R)=0$, i. e. $\LCMo ^{d-1}_J$ is right exact. Let
$Q\in \Assh (R)$ be arbitrary. The canonical epimorphism $R\to
R/Q=:\overline R$ induces a monomorphism
$$D_{\overline R}(\LCMo ^{d-1}_{J\overline R}(\overline R))=D_R(\LCMo ^{d-1}_J(\overline
R))\to D(\LCMo ^{d-1}_J(R)).$$ If $\{ 0\} \in \Ass _{\overline
R}(D_{\overline R}(\LCMo ^{d-1}_{J\overline R}(\overline R)))$ then
$Q\in \Ass _R(D_R(\LCMo ^{d-1}_J(R)))$, and so we may assume that
$R$ is a domain. If we can write $J=J_1\cap J_2$ with non-$\goth
m$-primary ideals $J_1,J_2$ of $R$ such that $J_1+J_2$ is $\goth
m$-primary then, because of $\LCMo ^d_{J_1}(R)=\LCMo ^d_{J_2}(R)=0$
(Hartshorne-Lichtenbaum vanishing), a Mayer-Vietoris sequence
argument shows the existence of an epimorphism
$$\LCMo ^{d-1}_J(R)\to \LCMo ^d_{J_1+J_2}(R)=\LCMo ^d_\goth m(R).$$
But then theorem 3.1.3 (iv) implies that
$$\{ 0\} =\Assh (R)=\Ass _R(D(\LCMo ^d_\goth m(R)))\subseteq \Ass _R(D(\LCMo
^{d-1}_J(R))).$$ If there is no such decomposition $J=J_1\cap J_2$
of $J$ we may assume that $J$ is a prime ideal; but then the
statement follows from theorem 3.2.4. Now we turn to the general
case, i. e. we assume there is a $Q\in \Assh (R)$ such that $\dim
(R/(J+Q))=0$. We define $U(R)$ to be the intersection of all
$Q^\prime $-primary components of a primary decomposition of the
zero ideal in $R$ for all $Q^\prime \in \Assh (R)$. Apparently we
have $\Ass _R(R/U(R))=\Assh (R)$ and $\dim (U(R))<d$. Because of the
latter fact the short exact sequence $0\to U(R)\buildrel \subseteq
\over \to R\to R/U(R)\to 0$ induces an exact sequence
$$0\to D(\LCMo ^{d-1}_J(R/U(R)))\to D(\LCMo ^{d-1}_J(R))\to D(\LCMo
^{d-1}_J(U(R))).$$ Trivially $\dim _R(\Supp _R(\LCMo
^{d-1}_J(U(R))))\leq d-1$ holds. By considering $R/U(R)$ rather then
$R$ we may assume that $\Ass _R(R)=\Assh (R)$. We write $0=I^\prime
\cap I^{\prime \prime }$ with ideals $I^\prime ,I^{\prime \prime }$
of $R$ such that $\Ass _R(R)=\Ass _R(R/I^\prime )\cup \Ass
_R(R/I^{\prime \prime })$ and $\dim (R/(J+Q^\prime ))\geq 1$ for all
$Q^\prime \in \Ass _R(R/I^\prime )$ and $\dim (R/(J+Q^{\prime \prime
}))=0$ for all $Q^{\prime \prime }\in \Ass _R(R/I^{\prime \prime
}))$. It follows that $\dim (R/(J+I^{\prime \prime }))=0$. By using
a Mayer-Vietoris argument and the facts that $\LCMo ^d_J(R/I^\prime
)=0$ (Hartshorne-Lichtenbaum vanishing) and $\LCMo ^i_J(R/I^{\prime
\prime })=\LCMo ^i_\goth m(R/I^{\prime \prime })$ for all $i\in
\Naturalsign $ we get a short exact sequence
$$D(\LCMo ^{d-1}_\goth m(R/I^\prime +I^{\prime \prime }))\to D(\LCMo
^{d-1}_J(R/I^\prime ))\oplus D(\LCMo ^{d-1}_\goth m(R/I^{\prime
\prime }))\to
$$
$$\to D(\LCMo ^{d-1}_J(R))\to D(\LCMo ^{d-2}_\goth m(R/(I^\prime +I^{\prime
\prime }))).$$ It is clear that we have
$$\dim _R(\Supp _R(D(\LCMo ^{d-1}_\goth m(R/(I^\prime +I^{\prime \prime
})))))\leq d-1$$ and $$\dim _R(\Supp _R(D(\LCMo ^{d-2}_\goth
m(R/(I^\prime +I^{\prime \prime })))))\leq d-1.$$ $R$ is complete
and so we can use local duality to conclude that
$$\dim _R(\Supp _R(D(\LCMo ^{d-1}_\goth m(R/I^{\prime \prime }))))\leq d-1.$$
Thus we get, by what is already shown,
$$\Assh (R)\cap \Ass _R(D(\LCMo ^{d-1}_J(R)))=\{ Q\in \Ass
_R(D(\LCMo ^{d-1}_J(R/I^\prime )))\vert \dim (R/Q)=d\} =\Assh
(R/I^\prime ).$$
\bigskip
The following theorems 3.2.6 (where $R$ is not necessarily complete)
and 3.2.7 (where $R$ will be complete) contain the main results.
\bigskip
{\bf 3.2.6 Theorem}
\par
Let $(R,\goth m)$ be a $d$-dimensional noetherian local ring and
$J\subseteq R$ an ideal such that $\dim (R/J)=1$ and $\LCMo
^d_J(R)=0$. Then
$$\Assh (R)=\Assh (D(\LCMo ^{d-1}_J(R)))$$
holds. \par Proof: \par
One has $\LCMo ^d_{J\hat R}(\hat R)=\LCMo
^d_J(R)\otimes _R\hat R=0$ and
$$\eqalign {D_{\hat R}(\LCMo ^{d-1}_{J\hat R}(\hat R))&=D_{\hat R}(\LCMo
^{d-1}_J(R)\otimes \hat R)\cr &=\Hom _R(\LCMo ^{d-1}_J(R),D_{\hat
R}(\hat R))\cr &=D_R(\LCMo ^{d-1}_J(R))\cr }$$ Therefore, every
$\hat R$-monomorphism $\varphi :\hat R/\goth P\to D_{\hat R}(\LCMo
^{d-1}_{J\hat R}(\hat R))$, where $\goth P$ is a prime ideal of
$\hat R$, induces an $R$-monomorphism $\hat R/\goth P\to D_R(\LCMo
^{d-1}_J(R))$. On the other hand we have a $R$-monomorphism $R/\goth
p\to \hat R/\goth P$, where $\goth p:=\goth P\cap R$. Composition of
these monomorphisms gives us a monomorphism
$$R/\goth p\to D_R(\LCMo ^{d-1}_J(R)).$$
Because of $\Assh (R)=\{ \goth P\cap R\vert \goth P\in \Assh (\hat
R)\}$ we may assume that $R$ is complete. But then the statement
follows from lemma 3.2.5.
\bigskip
{\bf 3.2.7 Theorem}
\par
Let $R$ be a $d$-dimensional local complete ring and $J\subseteq R$
an ideal such that $\dim (R/J)=1$ and $\LCMo ^d_J(R)=0$. Then
$$\Ass _R(D(\LCMo ^{d-1}_J(R))=\{ P\in \Spec(R)\vert \dim (R/P)=d-1,\dim
(R/(P+J))=0\} \cup \Assh (R)$$ holds. \par Proof: \par Let $P\in
\Spec (R)$. If $\dim (R/P)\leq d-2$ we have
$$\Hom _R(R/P,D(\LCMo
^{d-1}_J(R)))=D(\LCMo ^{d-1}_J(R/P))=0$$
and hence $P\not\in \Ass
_R(D(\LCMo^ {d-1}_J(R)))$. If $\dim (R/P)=1$ then (set $\overline
R:=R/P$):
$$\Hom _R(R/P,D(\LCMo ^{d-1}_J(R)))=D(\LCMo ^{d-1}_J(R/P))=D(\LCMo
^{d-1}_{J\overline R}(\overline R)).$$ $R$ is complete and so, by
Hartshorne-Lichtenbaum vanishing, the equivalence
$$\LCMo ^{d-1}_{J\overline R}(\overline R)\neq 0\iff \dim (\overline R/J\overline R)=0$$
holds. On the other hand we have $\overline R/J\overline R=R/(P+J)$
and, therefore
$$\{ P\in \Ass _R(D(\LCMo ^{d-1}_J(R)))\vert \dim (R/P)=d-1\}=$$
$$=\{ P\in \Spec (R)\vert \dim (R/P)=d-1,\dim (R/(P+J))=0\} \ \ .$$
Now the statement follows from lemma 3.2.5. \vfil \eject
\bigskip {\bf 4 The regular case and how to reduce to it}
\bigskip
\bigskip
By "regular case" we mean the following situation: Let $k$ be a
field, $R=k[[X_1,\dots ,X_n]]$ a power series algebra over $k$ in
$n$ variables and $I$ the ideal $(X_1,\dots ,X_h)R$ of $R$ ($1\leq
h\leq n$). We are interested in the associated prime ideals of
$$D:=D(\LCMo ^h_I(R))\ \ .$$
In the first subsection we will demonstrate how one can reduce
conjecture (*) to the regular case, in subsection 4.2 we present
results on $\Ass _R(D)$ for general $h$; subsection 4.3 concentrates
on the case $h=n-2$, which is in some sense the "first" interesting
case.
\bigskip
{\bf 4.1 Reductions to the regular case}
\bigskip
Suppose that $(R,\goth m)$ is a noetherian local ring. After
completing $R$, we can write $R$ as a quotient of a regular local
ring $S$; on the other hand we can find a regular local subring $S$
of $R$ such that $R$ is module-finite over $S$. We will use both
methods to reduce to the regular case, i. e. to make facts about
$\Ass _S(D_S)$ into facts about $\Ass _R(D_R)$ ($D_S$ and $D_R$
stand for the Matlis duals of local cohomology modules of $S$ resp.
$R$), see remark 4.1.1 and theorem 4.1.2 for details.
\bigskip
{\bf 4.1.1 Remark}
\par
Suppose that $(R,\goth m)$ is a noetherian local equicharacteristic
domain, $i\geq 1$ and $x_1,\dots ,x_i$ is a sequence in $R$ such
that $\LCMo ^i_{(x_1,\dots ,x_i)R}(R)\neq 0$. Suppose furthermore,
that one wants to show $\{ 0\} \in \Ass _R(D(\LCMo ^i_{(x_1,\dots
,x_i)R}(R)))$ (that is conjecture (*) in the equicharacteristic
case). W. l. o. g. one can replace $R$ by $\hat R/\goth p_0$, where
$\hat R$ is the $\goth m$-adic completion of $R$ and $\goth p_0\in
\Spec (\hat R)$ is lying over the zero ideal of $R$ (because then
$$\eqalign {D_{\hat R/\goth p_0}(\LCMo ^i_{(x_1,\dots ,x_i)(\hat
R/\goth p_0)}(\hat R/\goth p_0))&=\Hom _{\hat R/\goth p_0}(\LCMo
^i_{(x_1,\dots ,x_i)(\hat R)}(\hat R)\otimes _{\hat R}(\hat R/\goth
p_0),\Hom _{\hat R}(\hat R/\goth p_0,\InjH _{\hat R}(k))))\cr &=\Hom
_{\hat R}(\LCMo ^i_{(x_1,\dots ,x_i)(\hat R)},\Hom _{\hat R}(\hat
R/\goth p_0,\InjH _{\hat R}(k)))\cr &=\Hom _{\hat R}(\hat R/\goth
p_0,\Hom _{\hat R}(\LCMo ^i_{(x_1,\dots ,x_i)\hat R}(\hat R),\InjH
_{\hat R}(k)))\cr &=\Hom _{\hat R}(\hat R/\goth p_0,\Hom _{\hat
R}(\LCMo ^i_{(x_1,\dots ,x_i)\hat R}(\hat R),\InjH _R(k)))\cr &=\Hom
_{\hat R}(\hat R/\goth p_0,D_R(\LCMo ^i_R(\LCMo ^i_{(x_1,\dots
,x_i)R}(R))))\cr }$$ contains en element $d$ with $\hat
R$-annihilator $\goth p_0$, i. e. with $R$-annihilator $\goth
p_0\cap R=0$; but $d$ is naturally an element of $D_R(\LCMo
^i_{(x_1,\dots ,x_i)R}(R))$), and so we may assume that $(R,\goth
m)$ is a noetherian local equicharacteristic complete domain. Let
$k$ be a coefficient field of $R$. Now if we use a surjective
$k$-algebra homomorphism $k[[X_1,\dots ,X_n]]\to R$ ($k[[X_1,\dots
,X_n]]$ is a power series algebra over $k$ in $n$ variables
$X_1,\dots ,X_n$) mapping $X_1,\dots ,X_n$ to $x_1,\dots x_n$,
respectively, we can reduce to the following problem (note that,
below, $\goth p$ corresponds to the zero ideal of $R$):
\par
If $R=k[[X_1,\dots ,X_n]]$ is a power series ring over a field $k$
in $n$ variables $X_1,\dots ,X_n$, $1\leq i\leq n$, $\goth q\in \Ass
_R(D(\LCMo ^i_{(X_1,\dots ,X_i)R}(R)))$, $\goth p\in \Spec (R)$,
$\goth p\subseteq \goth q$, then $\goth p\in \Ass _R(D(\LCMo
^i_{(X_1,\dots ,X_i)R}(R)))$ holds (that is: The set $\Ass
_R(D(\LCMo ^i_{(X_1,\dots ,X_i)R}(R)))$ is stable under
generalization).
\par
Thus we have reduced conjecture (*) (in the equicharacteristic case)
to the preceding statement, a similar reduction is possible in the
case of mixed characteristic. \vfil \eject {\bf 4.1.2 Theorem}
\par
Let $(R,\goth m)$ be a noetherian local complete ring with
coefficient field $k\subseteq R$, $l\in \Naturalsign ^+$ and
$x_1,\dots ,x_l\in R$ a part of a system of parameters of $R$. Set
$I:=(x_1,\dots ,x_l)R$. Let $x_{l+1},\dots ,x_d\in R$ be such that
$x_1,\dots ,x_d$ is a system of parameters of $R$. Denote by $R_0$
the (regular) subring $k[[x_1,\dots ,x_d]]$ of $R$. Then if $\Ass
_{R_0}(D(\LCMo ^l_{(x_1,\dots ,x_l)R_0}(R_0)))$ is stable under
generalization, $\Ass _R(D(\LCMo ^l_I(R)))$ is also stable under
generalization.
\par
Proof: \par Set $X:=\Ass _R(D(\LCMo ^m_I(R)))$ and let $\goth p_1\in
\Spec (R),\goth p\in X,\goth p_1\subseteq \goth p$. We have to show
$\goth p_1\in X$. The hypothesis on $\goth p$ implies
$$0\neq \LCMo ^l_I(R/\goth p)=\LCMo ^l_{(x_1,\dots ,x_m)R_0}(R_0/\goth p\cap
R_0)\otimes _ {R_0}R.$$ But $\Ass _{R_0}(D(\LCMo ^l_{(x_1,\dots
,x_l)R_0}(R_0)))$ is stable under generalization and so by using
Matlis duality we first conclude $\goth p\cap R_0\in \Ass
_{R_0}(D(\LCMo ^l_{(x_1,\dots ,x_l)R_0}(R_0)))$ and then, by using
the stableness hypothesis again, $\goth p_1\cap R_0\in \Ass
_{R_0}(D(\LCMo ^l_{(x_1,\dots ,x_l)R_0}(R_0)))$. Now the existence
of an $R_0$-linear injection $R_0/\goth p_1\cap R_0\to D(\LCMo
^l_{(x_1,\dots ,x_l)R_0}(R_0))$ implies the existence of an
$R$-linear injection
$$\eqalign {\Hom _{R_0}(R,R_0/\goth p_1\cap R_0)&\to \Hom _{R_0}(R,D(\LCMo
^l_{(x_1,\dots ,x_l)}(R_0)))\cr &=\Hom _{R_0}(\LCMo ^l_{(x_1,\dots
,x_l)R_0}(R_0),(\Hom _{R_0}(R,\InjH _{R_0}(k))))\cr &=D(\LCMo
^l_I(R)),\cr }$$ where the last equality follows from the fact $\Hom
_{R_0}(R,\InjH _{R_0}(k))=\InjH _R(k)$. Thus it is sufficient to
show
$$\goth p_1\in \Ass _R(\Hom _{R_0}(R,R_0/\goth p_1\cap R_0))=\Hom
_{R_0/\goth p_1\cap R_0}(R/(\goth p_1\cap R_0)R,R_0/\goth p_1\cap
R_0)\ \ .$$
But $R$ is finite as $R_0$-module and so $\Hom
_{R_0/\goth p_1\cap R_0}(R/\goth p_1,R_0/\goth p_1\cap R_0)\neq 0$;
on the other hand $\goth p_1$ is minimal in the support of $R/(\goth
p_1\cap R_0)R$ and so, combining these facts, the statement of
theorem 4.1.2 follows.
\bigskip
\bigskip
{\bf 4.2 Results in the general case, i. e. $h$ is arbitrary}
\bigskip
We collect some properties of $\Ass _R(D)$ in the regular case; note
that $R$ does not have to contain a field:
\bigskip
{\bf 4.2.1 Theorem} \par Let $(R,\goth m)$ be a noetherian local
complete regular ring. Let $X_1,\dots ,X_n$ be a regular system of
parameters of $R$, $n=\dim (R)$. Set $I:=(X_1,\dots ,X_h)R$ for some
$h\in \{ 1, \dots ,n\} $. Set $D:=D(\LCMo ^h_I(R))$.
\par
(i) For $h=n$ one has
$$\Ass _R(D)=\{ \{ 0\} \}\ \ .$$
(ii) For $h=n-1$ one has
$$\Ass _R(D)=\{ \{ 0\} \} \cup \{ pR\vert p\in R\hbox { prime
element}, p\not\in I\} \ \ .$$ (iii) For general $h$ the following
statements hold: \par ($\alpha $) For every $\goth p\in \Spec (R)$
the implication
$$\goth p\in \Ass _R(D)\Longrightarrow \height (\goth p)\leq n-h$$
holds.
\par
($\beta $) For every $\goth p\in \Spec(R)$ such that $\height (\goth
p)=n-h$ one has the equivalence
$$\goth p\in \Ass _R(D)\iff I+\goth p\hbox { is }\goth m\hbox
{-primary}\ \ .$$ ($\gamma $) Every $f\in I\setminus \goth mI$ is
contained in no $\goth p\in \Ass _R(D)$; in particular, if $f=p$ is
a prime element of $R$ (such that $p\in I\setminus \goth mI$), one
has
$$pR\not\in \Ass _R(D)\ \ .$$
($\delta $) If $p\in R$ is a prime element such that $p\not\in I$
then
$$pR\in \Ass _R(D)$$
holds. \par
Proof:
\par
(i) follows from $D(\LCMo ^n_\goth m(R))=R$, (ii) from theorem
3.2.7. (iii) ($\alpha )$ and ($\beta $) follow from theorem 3.1.3
(i) resp. (v). (iii) ($\gamma $): In the given situation one has
$$\Hom _R(R/fR,D)=D(\LCMo ^h_I(R)/f\LCMo ^h_I(R))=D(\LCMo
^h_{I(R/fR)}(R/fR))$$ and $\LCMo ^h_{I(R/fR)}(R/fR)=0$, because $f$
is a minimal generator of $I$, i. e. $I(R/fR)$ can be generated by
$h-1$ elements; thus multiplication by $f$ is injective on $D$, the
statement follows. (iii) ($\delta $) follows from theorem 3.1.3
(ii).
\bigskip
{\bf 4.2.2 Remark}
\par
In the situation of theorem 4.2.1, the largest $h$, for which we
cannot completely determine $\Ass _R(D)$, is $h=n-2$; the theorem
leaves open the question which prime ideals $\goth p=pR$, $p\in R$ a
prime element, $p\in I$, are associated to $R$. In the next
subsection we will concentrate on the case $h=n-2$. We will give a
partial answer to this open question (see corollary 4.3.1) and we
will see (remarks 4.3.2 (i) and (ii)) that both $pR\in \Ass _R(D)$
and $pR\not\in \Ass _R(D)$ can occur (both in the special case $p\in
I$, $h=n-2$).
\bigskip
{\bf 4.2.3 Theorem}
\par
Let $(R_0,\goth m_0)$ be a noetherian local complete
equicharacteristic ring, let $\dim (R_0)=n-1$, $k\subseteq R_0$ a
coefficient field of $R_0$, $1\leq h\leq n$. Let $x_1,\dots ,x_n$ be
elements of $R_0$ such that $\sqrt {(x_1,\dots ,x_n)R_0}=\goth m_0$.
Set $I_0:=(x_1,\dots ,x_h)R_0$. Let $R:=k[[X_1,\dots ,X_n]]$ be a
power series algebra over $k$ in the variables $X_1,\dots ,X_n$,
$I:=(X_1,\dots ,X_h)R$. Then the $k$-algebra homomorphism $R\to R_0$
determined by $X_i\mapsto x_i$ $(i=1,\dots n$) induces a
module-finite ring map $\iota :R/fR\to R_0$ for some prime element
$f\in R$. We set
$$D:=D(\LCMo ^h_I(R))\ \ .$$
Then \smallskip (i) $D$ has an associated prime which contains $f$
if and only if $\LCMo ^h_{I_0}(R_0)\neq 0$. \smallskip Furthermore
if $R_0$ is regular and $\height (I_0)<h$, the following statements
hold:
\smallskip
(ii) $D$ has no associated prime ideal which contains $f$ and has
height $n-h$.
\par
(iii) If $\LCMo ^h_{I_0}(R_0)\neq 0$, ($f$ is contained in an
associated prime of $D$ and) every maximal element $\goth q$ of
$\Ass _R(D)$ containing $f$ has $\dim (R/\goth q)>h$; we will see
below (remark 4.3.2 (ii)) that this situation really occurs and,
therefore, it is in general not true that all maximal elements of
$\Ass _R(D(\LCMo ^h_I(R)))$ have dimension $h$; note that his was
conjecture (+) in [HS1, section 0] (see also remark 1.2.4 for more
details on this).
\par
Proof:
\par
(i) follows from
$$\eqalign {\exists _{\goth p\in \Ass _R(D)}f\in \goth p&\iff \Hom _R(R/fR,D)\neq 0\cr &\iff D(\LCMo ^h_I(R)/f\LCMo ^h_I(R))\neq 0\cr &\iff D(\LCMo ^h_I(R/fR))\neq 0
\cr &\iff \LCMo ^h_I(R/fR)\neq 0\cr &\iff \LCMo ^h_{I_0}(R_0)\neq 0\
\ .}$$ Note that, for the last equivalence, we use the fact that via
$\iota $ $R_0$ is a finite and torsion-free $R/fR$-module.
\par
From now on we assume, in addition, that $R_0$ is regular and that
$\height (I_0)<h$. \par
(ii) We assume, to the contrary, that there
is a prime ideal $\goth p\in \Ass _R(D)$ such that $\height (\goth
p)=n-h$ and such that $f\in \goth p$: $R_0$ is module-finite over
$R/fR$, and so there exists $\goth q\in \Spec (R_0)$ such that
$\goth q\cap (R/fR)=\goth p/fR$. But now $\goth q\cap R=\goth p$
implies
$$\height (\goth q)=(n-1)-\dim (R_0/\goth q)=(n-1)-\dim (R/\goth p)=\height (\goth p)-1=n-h-1$$
and therefore, from $\height (I_0)<h$, we conclude
$$\height (I_0+\goth q)<n-1\ \ ,$$
which means that $(I_0+\goth q)/\goth q$ is not $\goth m_0/\goth q$-primary in $R_0/\goth q$. Hence, by Hartshorne-Lichtenbaum
vanishing,
$$\LCMo ^h_{I_0}(R_0/\goth q)=0\ \ .$$
But $R_0/\goth q$ is a torsion-free finite $R/\goth p$-module, and
so the last vanishing result implies
$$\LCMo ^h_I(R/\goth p)=0\ \ ,$$
which contradicts the assumption $\goth p\in \Ass _R(D)$.
\par
(iii) The first statement implies that there is an associated prime
of $D$ which contains $f$ and (ii) shows that every such prime ideal
$\goth p$ has height smaller than $h$.
\bigskip
\bigskip
{\bf 4.3 The case $h=\dim (R)-2$, i. e. the set $\Ass _R(D(\LCMo
^{n-2}_{(X_1,\dots ,X_{n-2})R}(k[[X_1,\dots ,X_n]])))$}
\bigskip
We can give a partial answer to the question which height one prime
ideals contained in $I$ are associated to $D$:
\bigskip
{\bf 4.3.1 Corollary}
\par
If we are in the special case where $h=n-2$, $R_0$ is regular and
$\height (I_0)<h$ in the situation of theorem 4.2.3, we clearly have
(because of theorem 4.2.3 (ii))
$$fR\in \Ass _R(D)\iff \LCMo ^{n-2}_{I_0}(R_0)\neq 0\ \ .$$
In this case, $fR$ is a maximal element of $\Ass _R(D)$. By [HL,
Theorem 2.9], the latter holds if and only if $\dim (R_0/I_0)\geq 2$
and $\Spec (\overline {R_0}/I_0\overline {R_0})\setminus \{ \goth
m_0(\overline {R_0}/I_0\overline {R_0})\} $ is connected, where
$\overline {R_0}$ is defined as the completion of the strict
henselization of $R_0$; this means that $\overline {R_0}$ is
obtained from $R_0$ by replacing the coefficient field $k$ by its
separable closure in any fixed algebraic closure of $k$.
\bigskip
{\bf 4.3.2 Remarks}
\par
(i) In the situation of the statement (i) of theorem 4.2.3, it can
easily happen that both $f\in I$ and $\LCMo ^h_{I_0}(R_0)=0$ hold;
then we have, in particular, $fR\not\in \Ass _R(D)$. Hence, in
general, not all height one prime ideals contained in $I$ are
associated to $D$. In fact, if $\height (I_0)<h$, then $f$ is
necessarily contained in $I$. Hence, if $\ara (I_0)<h$, then both
$f\in I$ and $\LCMo ^h_{I_0}(R_0)=0$ hold and therefore one has
$fR\not\in \Ass _R(D)$.
\par
(ii) In the situation of corollary 4.3.1, it can happen that $fR\in
\Ass _R(D)$. For example, we can take
$$n=5, h=3, k={\bf Q}\hbox{ (the
rationals), }R_0={\bf Q}[[y_1,y_2,y_3,y_4]]\ \ ,$$ a power series
algebra over $\bf Q$ in the variables $y_1,y_2,y_3,y_4$. We set
$$x_1=y_1y_3,
x_2=y_2y_4,x_3=y_1y_4+y_2y_3,x_4=y_1+y_3,x_5=y_2+y_4\ \ .$$
Then
$\height (I_0)=2$ and $\LCMo ^3_{I_0}(R_0)\neq 0$. Furthermore
$$f:=-X_2X_4^2+X_3X_4X_5-X_1X_5^2+4X_1X_2-X_3^2\in R$$
generates the
kernel of the $\bf Q$-algebra homomorphism $R\to R_0$, which is
determined by $X_i\mapsto x_i$ ($i=1,\dots ,5$), where $R$ is
defined as the power series ring ${\bf Q}[[X_1,X_2,X_3,X_4,X_5]]$
over $\bf Q$ in the variables $X_1,X_2,X_3,X_4,X_5$. Now, by
corollary 4.3.1, $fR$ is a maximal element of $\Ass _R(D)$. In
particular, this example clearly provides a counterexample to
conjecture (+) from [HS1, section 0] (see also remark 1.2.4 for
details on this).
\par
Proof:
\par
(i) We assume that $\height (I_0)<h$ and prove $f\in I$: Let $\goth
p_0$ be a prime ideal minimal over $I_0$ and such that $\height
(\goth p_0)<h-1$; the inclusion $\goth p_0\cap R\supseteq I+fR$
induces a surjection $R/(I+fR)\to R/\goth p_0\cap R$; on the other
hand, $R_0/\goth p_0$ is finite over $R/\goth p_0\cap R$. Therefore
we have
$$\dim (R_0/\goth p_0)=\dim (R/\goth p_0\cap R)\leq \dim (R/(I+fR))\ \
.$$ Now, if $f$ was not contained in $I$, one would conclude $\dim
(R_0/\goth p_0)\leq n-h-1$ and, hence, $\height (\goth p_0)\geq h$.
(ii) It is easy to see that
$$\sqrt {I_0}=(y_1,y_2)R_0\cap (y_3,y_4)R_0$$
and so $\height (I_0)=2$ and a Mayer-Vietoris sequence argument,
applied to the ideals $(y_1,y_2)R_0$ and $(y_3,y_4)R_0$, shows that
$\LCMo ^3 _{I_0}(R_0)\neq 0$. $f$ generates the kernel of the $\bf
Q$-algebra homomorphism $R\to R_0$; this can be seen e. g. by
observing that $f$, as an element of ${\bf Q}[X_1,X_2,X_3,X_4,X_5]$,
generates the kernel of the associated map over polynomial instead
of formal power series rings, which in turn is true, because first
of all an easy calculation shows that $f$ is in this kernel and,
secondly, as a polynomial, $f$ is irreducible, which can either be
seen by a direct calculation or by using a computer algebra system
like, for example, Macaulay 2. The rest follows from corollary
4.3.1.
\bigskip
Assume that $p\in I$ is a prime element. The next example and, more
generally, theorem 4.3.4 show that under certain conditions, $p$ is
contained in infinitely many associated height two prime ideals of
$D$. This is useful for two reasons: It will lead to a
generalization of an example of Hartshorne of a non-artinian (but
zero-dimensional) local cohomology module (see theorem 6.2.3); and
secondly, it will show that either conjecture (*) holds (for
$h=n-2$) or, if not, $D$ satisfies a remarkable property (see remark
4.3.6 for details on this property).
\bigskip
{\bf 4.3.3 Example} \par Still in the above situation consider
$p:=X_1X_n+X_2X_{n-1}\in I\cap (X_{n-1},X_n)R$. For every $\lambda
\in k$ set $\goth p_\lambda :=(X_{n-1}+\lambda X_1,X_n-\lambda
X_2)R$. Then $\goth p_\lambda \in \Ass _R(D(\LCMo ^{n-2}_I(R)))$
holds (this follows from theorem 3.1.3 (v)) and
$$p=X_1(X_n-\lambda X_2)+X_2(X_{n-1}+\lambda X_1)$$
is contained in every $\goth p_\lambda $.
\bigskip
The same idea works more general:
\bigskip
{\bf 4.3.4 Theorem}
\par
Let $R=k[[X_1,\dots ,X_n]]$ be a power series ring in the variables
$X_1,\dots ,X_n$ ($n\geq 4$) over a field $k$ and let $I$ be the
ideal $(X_1,\dots ,X_{n-2})R$ (i. e. $h=n-2$ in the above notation).
Furthermore, let $p\in R$ be a prime element such that $p\in I\cap
(X_{n-1},X_n)R$.
\par
The set
$$\{ \goth p\in \Spec (R)\vert \goth p\in \Ass _R(D(\LCMo ^{n-2}_I(R))),p\in
\goth p, \height(\goth p)=2\} $$ is infinite.
\par
Proof: \par It is easy to see that there exist $f,g\in I, f\not\in
(X_{n-1},X_n)R$ and $l\geq 1$ such that
$$p=X_n^lf+X_{n-1}g$$
holds. Let $m\in \Naturalsign ^+$ be arbitrary. We have
$$p=(X_n^l+X_1^mg)f+(X_{n-1}-X_1^mf)g$$
and so
$$p\in I_m:=(X_n^l+X_1^mg,X_{n-1}-X_1^mf)R\ \ .$$
The elements $$X_1,\dots ,X_{n-2},X_n^l+X_1^mg,X_{n-1}-X_1^mf$$
form
a system of parameters of $R$ and so there exists a $\goth p_m\in
\Ass _R(D(\LCMo ^{n-2}_I(R)))$ containing $I_m$. For $m,m^\prime \in
\Naturalsign ^+,m\neq m^\prime $
$$\sqrt {I_m+I_{m^\prime }}=(X_1,X_n,X_{n-1})R\cap \sqrt
{(X_{n-1},X_n,f,g)R}$$ holds; in particular, all primes containing
$I_m+I_{m^\prime }$ have height at least three. The statement
follows.
\bigskip
{\bf 4.3.5 Remark}
\par
In the situation of theorem 4.3.4, conjecture (*) would clearly
imply $pR\in \Ass _R(D)$. Now, if $pR$ was not associated to $D$,
there would be a remarkable consequence, which is somewhat
counterintuitive (note that, in the situation below, the way we
choose $\goth p_l$ has nothing to do with the way how we choose
$d_{l+1},d_{l+2},\dots $) and which is explained in the next remark.
\bigskip
{\bf 4.3.6 Remark}
\par
Let $R=k[[X_1,\dots ,X_n]]$ be a power series algebra over a field
$k$ in the variables $X_1,\dots ,X_n$; set $I=(X_1,\dots ,X_{n-2})R,
D:=D(\LCMo ^{n-2}_I(R))$ and $Y:=X_1\cdot \dots \cdot X_{n-2}$;
furthermore, assume that $p\in I$ is a prime element of $R$ such
that there are infinitely many height two prime ideals associated to
$D$ and containing $p$ (by theorem 4.3.4, this is true for example,
if $p\in (X_{n-1},X_n)R$ holds) and such that $pR\not\in \Ass
_R(D)$.
\smallskip Then for any sequence $(\goth p_i)_{i\in \Naturalsign} $
of pairwise different elements of $\Ass _R(\Hom_R(R/pR,D))$ and for
any sequence $(d_i)_{i\in \Naturalsign }$ in $D$ such that $\Ann
_R(d_i)=\goth p_i$ for every $i\in \Naturalsign $, there exists a
number $N$ such that
$$\Ann _R(d_{l+1}Y^{l+1}+d_{l+2}Y^{l+2}+\dots )\subseteq \goth p_l$$
holds for every $l>N$ (see the proof below for remarks on our
notation).
\par
Proof:
\par
It is well-known that $\LCMo^{n-2}_I(R)$ is the cohomology in the
$n-2$-th degree of the \v{C}ech-complex
$$0\to R\to \oplus _{i_1=1}^{n-2}R_{X_{i_1}}\to \oplus _{1\leq
i_1<i_2\leq n-2}R_{X_{i_1}X_{i_2}}\to \dots \to R_{X_1\dots
X_{n-2}}\to 0\ \ ;$$ Therefore we can write $\LCMo ^{n-2}_I(R)$ as
$$k[[X_{n-1}X_n]][X_1^{-1},\dots ,X_{n-2}^{-1}]\ \ ;$$
by definition, this expression shall stand for $$\oplus _{i_1,\dots
,i_{n-2}\leq 0}k[[X_{n-1},X_n]]\cdot X_1^{i_1}\cdot \dots \cdot
X_{n-2}^{i_{n-2}}$$ with the obvious $R$-module structure on it.
Using this, a straight-forward calculation shows
$$D=k[X_{n-1}^{-1},X_n^{-1}][[X_1,\dots ,X_{n-2}]]\ \ ,$$
where we use similar notation like above, i. e. we write the
elements of $D$ as formal power series in $X_1,\dots ,X_{n-2}$ and
coefficients in
$$k[X_{n-1}^{-1},X_n^{-1}]=\oplus _{i_{n-1},i_n\leq 0}k\cdot X_{n-1}^{i_{n-1}}\cdot
X_n^{i_n}\ \ .$$ Using this description of $D$ it is clear that
$d_{l+1}Y^{l+1}+d_{l+2}Y^{l+2}+\dots $ is an element of $D$ for
every $l\in \Naturalsign $. In the same way it is clear that the
element
$$d:=d_0+Y\cdot d_1+Y^2\cdot d_2+\dots \in D$$
is well-defined. By construction $p$ annihilates $d$ and so, because
of $pR\not\in \Ass _R(D)$, there exists $r\in \Ann _R(d)\setminus
pR$. We conclude
$$0=rd=rd_0+rYd_1+rY^2d_2+\dots $$
The height of every prime ideal associated to $D$ is at most two and
thus for every $l\in \Naturalsign $ either $\Ann _R(rd_l)=\goth p_l$
or $rd_l=0$ holds. The latter condition is equivalent to $r\in \goth
p_l$, which holds if and only if $\goth p_l$ contains the ideal
$(r,p)R$. Hence there are only finitely many $l\in \Naturalsign $
such that $rd_l=0$, the set
$$M:=\{ l\in \Naturalsign \vert rd_l=0\} $$
is finite. For every $l\in \Naturalsign \setminus M$ we have
$$\Ann _R(d_{l+1}Y^{l+1}+d_{l+2}Y^{l+2}+\dots )\subseteq \Ann
_R(rd_{l+1}Y^{l+2}+rd_{l+2}Y^{l+2}+\dots )=\Ann_R(rd_0+\dots
+rd_lY^l)\subseteq \goth p_l\ \ ;$$ note that the last inclusion
follows from lemma 4.3.7 below. In particular, for every $l>\max M$
we have
$$\Ann _R(d_{l+1}Y^{l+1}+d_{l+2}Y^{l+2}+\dots )\subseteq \goth
p_l\ \ ,$$ we can take $N:=\max M$.
\bigskip
{\bf 4.3.7 Lemma}
\par
In the situation of theorem 4.3.4, assume that $d_1,\dots ,d_n$ are
elements of $D$ such that for every $i=1,\dots ,n$ the ideal $\Ann
_R(d_i)=:\goth p_i$ is a height-two prime ideal of $D$ and such that
the $\goth p_i$ are pairwise different. Then
$$\Ann _R(d_1+\dots +d_n)=\goth p_1\cap \dots \cap \goth p_n$$
holds. \par
Proof:
\par
By induction on $n$, the case $n=1$ being trivial. We assume that
$n>1$ and that the lemma is true for smaller $n$. The inclusion
$\supseteq $ is trivial. Now, if there was an element $r\in \Ann
_R(d_1+\dots +d_n)\setminus \goth p_i$ for some $i\in \{ 1,\dots
,n\} $, we would have
$$-rd_1=rd_2+\dots +rd_n$$
and
$$\goth p_2\cap \dots \cap \goth p_n=\Ann _R(d_2+\dots
+d_n)\subseteq \Ann _R(rd_2+\dots +rd_n)=\goth p_1\ \ ,$$ which
would be a contradiction. \vfil \eject {\bf 5 On the meaning of a
small arithmetic rank of a given ideal}
\bigskip
\bigskip
We investigate the condition that the arithmetic rank of a given
ideal is small in the sense that it is one or two. We start with an
example of an ideal whose cohomological dimension is one but whose
arithmetic rank is two (example 5.1); this makes the question when
$\ara (I)\leq 1$ holds more difficult; we present criteria for this
condition and also for $\ara (I)\leq 2$ (theorem 5.2.5 and corollary
5.2.6). While this works equally well in the local and in the graded
case, we distinguish some subtle differences between these two cases
in the third subsection 5.3.
\bigskip
\bigskip
{\bf 5.1 An Example} \bigskip We start with an example of an ideal
$I$ of a noetherian ring $R$ where both $0=\LCMo ^2_I(R)=\LCMo
^3_I(R)=\dots $ and $\ara (I)\geq 2$ hold: Let $k$ be any field and
$R=k[[x,y,z,w]]$ a power series ring over $k$ in four variables. Set
$$f=xw-yz\ \ ,$$
$$g_1=y^3-x^2z, g_2=z^3-w^2y$$
and $$I=\sqrt {(f,g_1,g_2)R}\ \ .$$ The ideal $I\subseteq R$ is the
complete version of the vanishing ideal of a rational quartic curve
in projective three-space over $k$; it is well-known that
$I\subseteq R$ is a height-two prime ideal of $R$. We claim that
both $\LCMo ^s_I(R/fR)=0$ for every $s\geq 2$ and $\ara
(I(R/fR))\geq 2$ hold (the last statement may be known, we include a
proof for lack of a reference):
\par
Let $y_0,\dots ,y_3$ be new variables and set
$S:=k[[y_0,y_1,y_2,y_3]]$. Denote by $R_1$ the three-dimensional
subring $R_1:=k[[y_0y_1,y_0y_2,y_1y_3,y_2y_3]]$ of $S$. The ring
homomorphism
$$R\to R_1,x\mapsto y_0y_1,y\mapsto y_0y_2,z\mapsto y_1y_3,w\mapsto
y_2y_3$$ clearly induces an isomorphism
$$R/fR\cong R_1(\subseteq S)\ \ .$$
Now consider the $k$-linear map
$$k[y_0,y_1,y_2,y_3]\buildrel \varphi \over \to R_1$$
that sends  a term $y_0^{\alpha _0}y_1^{\alpha _1}y_2^{\alpha
_2}y_3^{\alpha _3}$ to $y_0^{\alpha _0}y_1^{\alpha _1}y_2^{\alpha
_2}y_3^{\alpha _3}\in R_1$ if $\alpha _0+\alpha _3=\alpha _1+\alpha
_2$ holds, and to zero otherwise. Note that $\varphi $ is
well-defined by construction and naturally induces a map
$$S=k[[y_0,y_1,y_2,y_3]]\buildrel \tilde \varphi \over \to R_1\ \ .$$
Now it is easy to see that $\tilde \varphi $ is $R_1$-linear and
makes $R_1$ into a direct summand in $S$ (as an $R_1$-submodule).
Thus $\LCMo ^2_I(R/fR)$ is isomorphic to a direct summand of $\LCMo
^2_{IS}(S)$. We have
$$IS=(g_1,g_2)S=((y_0y_2^3-y_1^3y_3)\cdot
y_0^2,(y_0y_2^3-y_1^3y_3)\cdot (-y_3^2))S$$ and
$$\sqrt {IS}=(y_0y_2^3-y_1^3y_3)S\ \ .$$
This implies $\LCMo ^2_{IS}(S)=0$ and thus, by what we have seen
above, $\LCMo ^2_I(R/fR)=0$. \smallskip Now we show $\ara
(I(R/fR))=2$: We assume $\ara (I(R/fR))\neq 2$; then we clearly have
$\ara (I(R/fR))=1$. Let $h\in R$ be such that
$$I(R/fR)=\sqrt {h(R/fR)}$$
holds. This implies
$$\sqrt {IS}=\sqrt {hS}\ \ .$$
We have seen before that
$$\sqrt {IS}=(y_0y_2^3-y_1^3y_3)S$$
holds. $S$ is a unique factorization domain and so there exist
$N\geq 1$ and $s\in S$ such that
$$h=(y_0y_2^3-y_1^3y_3)^N\cdot s\hbox { and } (y_0y_2^3-y_1^3y_3)\not \vert
s$$ hold. From $h\in R_1\subseteq S$ it follows that all terms
$y_0^{\alpha _0}y_1^{\alpha _1}y_2^{\alpha _2}y_3^{\alpha _3}$ in
$h\in S$ have the property $\alpha _0+\alpha _3=\alpha _1+\alpha
_2$; on the other hand, all terms $y_0^{\alpha _0}y_1^{\alpha
_1}y_2^{\alpha _2}y_3^{\alpha _3}$ of $(y_0y_2^3-y_1^3y_3)^N$ have
the property $(\alpha _0+\alpha _3)-(\alpha _1+\alpha _2)=-2N$. So
we can assume that all terms $y_0^{\alpha _0}y_1^{\alpha
_1}y_2^{\alpha _2}y_3^{\alpha _3}$ of $s$ have the property $(\alpha
_0+\alpha _3)-(\alpha _1+\alpha _2)=2N$. But then $s$ cannot be a
unit in $S$ and so
$$(y_0y_2^3-y_1^3y_3)S=\sqrt {hS}=(y_0y_2^3-y_1^3y_3)S\cap \sqrt
{sS}$$ clearly leads to a contradiction.
\bigskip
\bigskip
{\bf 5.2 Criteria for $\ara (I)\leq 1$ and $\ara (I)\leq 2$}
\bigskip
{\bf 5.2.1 Remark} \par Let $(R,\goth m)$ be a noetherian local
ring. By $\InjH :=\InjH _R(R/\goth m)$ we denote an $R$-injective
hull of $R/\goth m$. Let $I$ be an ideal of $R$. Then the following
statements are equivalent:
\smallskip
(i) $\ara (I)\leq 1$.
\par
(ii) $\LCMo ^i_I(R)=0$ for $i\geq 2$ and $\exists f\in I: f $
operates surjectively on $\LCMo ^1_I(R)$.
\par
(iii) $\LCMo ^i_I(R)=0$ for $i\geq 2$ and $\exists f\in I: f $
operates injectively on $D(\LCMo ^1_I(R))$.
\par
(iv) $\LCMo ^i_I(R)=0$ for $i\geq 2$ and $I\not\subseteq \bigcup
_{\goth p\in \Ass _R(D(\LCMo ^1_I(R)))}\goth p$.
\smallskip
Furthermore, if conditions (ii) or (iii) hold, we have $\sqrt
I=\sqrt {fR}$.
\par
Proof:
\par
(ii) -- (iv) are obviously equivalent, we show (i) $\iff $ (ii):
\par
(i) $\Rightarrow $ (ii): Assume $\sqrt I=\sqrt {fR}$ for some $f\in
R$. $f$ clearly operates surjectively on $\LCMo ^1_{fR}(R)$; but
$\sqrt I=\sqrt {fR}$ implies $\LCMo ^1_I=\LCMo ^1_{fR}$.
\par (ii) $\Rightarrow $(i): $\LCMo ^1_I(\_)$ is right-exact on the
category of $R$-modules. Therefore we have an exact sequence
$$\LCMo ^1_I(R)\buildrel f\over \to \LCMo ^1_I(R)\to \LCMo ^1_I(R/fR)\to 0\ \
.$$ Thus $\LCMo ^1_I(R/fR)=0$ holds, implying $\LCMo ^1_I(R/\goth
p)=0$ for all $\goth p\in \Spec (R)$ containing $f$. But because of
our hypothesis, this means that we have $\LCMo ^i_I(R/\goth p)=0$
for all $\goth p\in \Spec (R)$ containing $I$ and for all $i\geq 1$.
Thus $f$ must be contained in every prime ideal of $R$ containing
$I$. $\sqrt I=\sqrt {fR}$ follows.
\bigskip
{\bf 5.2.2 Remark} \par Now we consider the following situation
(referred to from now on as graded situation): Let $K$ be a field,
$l\in \Naturalsign ,$
$$R=K[X_0,\dots ,X_N]/J$$ ($N\in \Naturalsign
,J\subseteq K[X_0,\dots ,X_N]$ a homogenous ideal, where every $X_i$
has a multidegree in $\Naturalsign ^l$),
$$I\subseteq R\hbox { a homogenous
ideal, }$$ $\goth m$ the maximal homogenous ideal $(X_0,\dots
,X_N)R$ of $R$, $\InjH :=\InjH _R(R/\goth m)$ an $R$-injective hull
of $R/\goth m$; $\InjH $ has a natural grading and serves also as a
*-$R$-injective hull of $R/\goth m$. Here we follow the use of
*-notation from [BS, in particular sections 12 and 13]. $*D$ shall denote the functor from the category of
graded $R$-modules to itself defined by $$(*D)(M):=*\Hom _R(M,\InjH
)$$ for a graded $R$-module $M$. We define the homogenous arithmetic
rank of $I$ to be
$$\ara ^h (I):=\min \{ l\in \Naturalsign \vert \exists r_1,\dots ,r_l\in R\hbox {
homogenous}:\sqrt I=\sqrt {(r_1,\dots ,r_l)R}\} $$ and we set
$$I^h:=\{ x\in I\vert x\hbox { is homogenous}\} \ \ .$$ Now, just
like in the local case, one can show that the following statements
are equivalent:
\smallskip
(i) $\ara ^h(I)\leq 1$.
\par
(ii) $\LCMo ^i_I(R)=0$ for $i\geq 2$ and $\exists $ homogenous $f\in
I: f $ operates surjectively on $\LCMo ^1_I(R)$.
\par
(iii) $\LCMo ^i_I(R)=0$ for $i\geq 2$ and $\exists $ homogenous
$f\in I: f $ operates injectively on $(*D)(\LCMo ^1_I(R))$.
\par
(iv) $\LCMo ^i_I(R)=0$ for $i\geq 2$ and $I^h \not\subseteq \bigcup
_{\goth p\in \Ass _R((*D)(\LCMo ^1_I(R)))}\goth p$.
\smallskip
Furthermore, if conditions (ii) or (iii) hold, we have $\sqrt
I=\sqrt {fR}$.
\bigskip {\bf 5.2.3 Definition}
\par
Let $(R,\goth m)$ be a noetherian local ring and $X\subseteq \Spec
(R)$ a subset. We say that $X$ satisfies prime avoidance if, for
every ideal $J$ of $R$,
$$J\subseteq \bigcup _{\goth p\in X}\goth p$$
implies
$$\exists \goth p_0\in X:J\subseteq \goth p_0\ \ .$$
\bigskip
{\bf 5.2.4 Definition} \par In the graded situation, let $X\subseteq
\Spec ^h(R):=\{ \goth p\in \Spec (R)\vert \goth p\hbox {
homogenous}\} $ be any subset. We say that $X$ satisfies homogenous
prime avoidance if, for every homogenous ideal $J$ of $R$,
$$J^h\subseteq \bigcup _{\goth p\in X}\goth p$$
implies
$$\exists \goth p_0\in X: J\subseteq \goth p_0\ \ .$$
{\bf 5.2.5 Theorem}
\smallskip
(i) Let $(R,\goth m)$ be a noetherian local ring and $I$ an ideal of
$R$ such that
$$0=\LCMo ^2_I(R)=\LCMo ^3_I(R)=\dots \eqno{(1)}$$
holds. Then
$$\ara (I)\leq 1\iff \Ass _R(D(\LCMo ^1_I(R)))\hbox { satisfies prime
avoidance}\ \ .$$ (ii) Let $R$ be graded and $I\subseteq R$ an
homogenous ideal such that $0=\LCMo ^2_I(R)=\LCMo ^3_I(R)=\dots $.
Then
$$\ara ^h(I)\leq 1\iff \Ass _R((*D)(\LCMo ^1_I(R)))\hbox { satisfies homgenous
prime avoidance}\ \ .$$
\par
Proof: \par
(i) We set
$$D:=D(\LCMo ^1_I(R))\ \ .$$
\par
$\Rightarrow $: Let $J\subseteq R$ be an ideal such that
$$J\subseteq \bigcup _{\goth p\in \Ass _R(D)}\goth p\ \ .$$
We claim that $\Hom _R(R/J,D)\neq 0$. Assumption: $\Hom
_R(R/J,D)=0$: It is a general fact that for every ideal $K\subseteq
R$ and every $R$-module $M$ one has
$$\Hom _R(R/K,D(M))=D(M/KM)\ \ .$$
(Proof of this general fact: If $K=(k_1,\dots k_l)R$ for some
$k_1,\dots ,k_l\in R$, the exact sequence
$$R^l\buildrel (k_1,\dots ,k_l)\over \to R\buildrel \hbox
{can.}\over \to R/K\to 0$$ induces an exact sequence
$$M^l\buildrel (k_1,\dots ,k_l)\over \to M\buildrel \hbox
{can.}\over \to M/KM\to 0\ \ ;$$ The functor $D$ is exact and so we
get an exact sequence
$$0\to D(M/KM)\buildrel \hbox {can.}\over \to D(M)\buildrel \pmatrix
{k_1\cr \vdots \cr k_l\cr }\over \to D(M)^l\ \ ,$$ from which the
statement $\Hom _R(R/K,D(M))=D(M/KM)$ follows.)
\par
We apply this general fact in the case $K=J$, $M=\LCMo ^1_I(R)$ and
conclude that
$$0=\Hom _R(R/J,D(\LCMo ^1_I(R)))=D(\LCMo ^1_I(R)/J\LCMo
^1_I(R))=D(\LCMo ^1_I(R/J)).$$ For the last equality we use the fact
that the functor $\LCMo ^1_I$ is right-exact (because of hypothesis
(1): $0=\LCMo ^2_I(R)=\LCMo ^3_I(R)=\dots $). But $D(\LCMo
^1_I(R/J))=0$ implies that $\LCMo ^1_I(R/J)=0$. Again, because of
hypothesis $(1)$, it follows that $\LCMo ^1_I(R/\goth p)=0$ for all
prime ideals $\goth p$ of $R$ containing $J$. Clearly, the last
condition implies $I\subseteq \goth p$ for all $\goth p$ containing
$J$, that is $I\subseteq \sqrt J$. There is an $x\in R$ such that
$\sqrt I=\sqrt {xR}$. Hence $x^l\in J$ for $l>>0$. So there is a
$\goth p\in \Ass _R(D)$ containing $x$. Now we have
$$0=\LCMo ^1_{xR}(R/\goth p)=\LCMo ^1_I(R/\goth p)$$
and thus
$$0=D(\LCMo ^1_I(R/\goth p))=\Hom _R(R/\goth p,D(\LCMo ^1_I(R)))$$
contradicting $\goth p\in \Ass _R(D)$. Thus the assumption $\Hom
_R(R/J,D)=0$ is false and the claim $\Hom _R(R/J,D)\neq 0$ is
proven; so there exists a $d\in D\setminus \{ 0\} $ such that
$J\subseteq \ann _R(d)$.
\par
$\Leftarrow $: We have to show the existence of an $x\in I$
operating surjectively on $\LCMo^1_I(R)$. Assume to the contrary
$$I\subseteq \bigcup _{\goth p\in \Ass_R(D)}\goth p\ \ .$$
From the hypothesis we get a $\goth p_0\in \Ass _R(D)$ such that
$I\subseteq \goth p_0$. But this $\goth p_0$ would satisfy
$$0\neq \LCMo ^1_I(R/\goth p_0)=0\ \ .$$
\par
(ii) The proof consists mainly of a graded version of the proof of
(i):
\par
$\Rightarrow $: Let $J\subseteq R$ be an homogenous ideal such that
$$J^h\subseteq \bigcup _{\goth p\in \Ass _R((*D)(\LCMo
^1_I(R)))}\goth p$$ and $x\in R^h$ an element such that $\sqrt
I=\sqrt {xR}$. We assume
$$\Hom _R(R/J,*\Hom _R(\LCMo ^1_I(R),\InjH ))=0$$
and remark that for the first $\Hom $ (in the preceding formula) it
would not make any difference if we replaced $\Hom $ by $*\Hom $.
This implies
$$*\Hom _R((R/J)\otimes _R\LCMo ^1_I(R),\InjH )=0$$
and hence $\LCMo ^1_I(R/J)=0$. Thus $I\subseteq \goth q$ for all
prime ideals $\goth q$ of $R$ containing $J$. This implies the
existence of a $\goth p_0\in \Ass _R((*D)(\LCMo ^1_I(R)))$ such that
$x\in \goth p_0$ contradicting $\LCMo ^1_I(R/\goth p_0)\neq 0$.
\par
$\Leftarrow $: We assume that for every $x\in I^h$ there exists a
$\goth p\in \Ass _R((*D)(\LCMo ^1_I(R)))$ such that $x\in \goth p$,
i. e.
$$I^h\subseteq \bigcup _{\goth p\in \Ass _R(*\Hom _R(\LCMo ^1_I(R),\InjH
))}\goth p\ \ .$$ There is a $\goth p_0\in \Ass _R(*\Hom _R(\LCMo
^1_I(R),\InjH ))$ containing $I$, contradicting $\LCMo ^1_I(R/\goth
p_0)\neq 0$.
\bigskip
Theorem 5.2.5 implies criteria for $\ara (I)\leq 2$ resp. for $\ara
^h(I)\leq 2$:
\bigskip
{\bf 5.2.6 Corollary}
\smallskip
(i) Let $(R,\goth m)$ be a noetherian local ring and $I$ an ideal of
$R$. Then $\ara (I)\leq 2$ if and only if there exists $g\in I $
such that $0=\LCMo ^2_I(R/gR)=\LCMo ^3_I(R/gR)=\dots $ and such that
$\Ass _R(D(\LCMo ^1_I(R/gR))$ satisfies prime avoidance.
\par
(ii) Let $R$ be a graded ring and $I$ an ideal of $R$. Then $\ara
^h(I)\leq 2 $ if and only if there exists a homogenous $g\in I $
such that $0=\LCMo ^2_I(R/gR)=\LCMo ^3_I(R/gR)=\dots $ and such that
$\Ass _R(D(\LCMo ^1_I(R/gR))$ satisfies homogenous prime avoidance.
\par
Proof: \par $\Rightarrow $ follows immediately from theorem 5.2.5
(i) resp. (ii); for the other implication observe that the
conditions on the right side imply $\ara _{R/gR}(I/(gR))=1$ resp.
$\ara ^h _{R/gR}(I/(gR))=1$ again by theorem 5.2.5 (i) resp. (ii).
\bigskip
\bigskip
{\bf 5.3 Differences between the local and the graded case}
\bigskip
{\bf 5.3.1 Lemma}
\smallskip
Let $R$ be a graded domain and $f\in R\setminus \{ 0\} $. Then the
ideal $\sqrt {fR}$ is homogenous if and only if $f$ is homogenous.
In particular, for any homogenous ideal $I$ of $R$ we have
$$\ara (I)\leq 1\iff \ara ^h(I)\leq 1\ \
.$$
\par
Proof:
\par
$\Leftarrow $ is clear. $\Rightarrow $: $R$ is $\Naturalsign
^l$-graded. This given grading may be seen as $l$ given
$\Naturalsign $-gradings on $R$ and so we may assume $l=1$. Let
$\delta :=\deg (f)$. Then $f_\delta (=$degree-$\delta $-part of $f)
\in \sqrt {fR}$, i. e. $\exists n\in \Naturalsign ^+$ and $\exists
g\in R: f_\delta ^n=fg$. $R$ is a domain and so $f$ (as well as $g$)
must be homogenous.
\bigskip
{\bf 5.3.2 Remark} \par In the graded situation, given graded
$R$-modules $M$ and $N$,
$$*\Hom _R(M,N)\subseteq \Hom _R(M,N)$$ holds. For finite $M$ one has
equality here, but for arbitrary $M$ equality does not hold in
general. In fact one has
$$\Ass _R(*\Hom _R(M,N))\subsetneq \Ass _R(\Hom _R(M,N))$$ in general as we
will see below in the case $M=\LCMo ^1_I(R)$, $N=\InjH :=\InjH
_R(R/\goth m)$; then we will also see that, in some sense, $\Ass
_R(\Hom _R(\LCMo ^1_I(R),\InjH))$ is much larger than $\Ass _R(*\Hom
_R(\LCMo ^1_I(R),\InjH))$.
\bigskip
{\bf 5.3.3 Definition and remark} \par Let $R$ be a graded ring and
$I\subseteq R$ homogenous ideal such that $0=\LCMo ^2_I(R)=\LCMo
^3_I(R)=\dots $. Let $f\in I$ be an element, not necessarily
homogenous. Now we define two conditions on $I$ and $f$:
$$\forall _{\goth p\in \Ass _R(\Hom _R(\LCMo ^1_I(R),\InjH ))} f\not \in
\goth p\leqno {(C_1)}$$
$$\forall _{\goth p\in \Ass _R(*\Hom _R(\LCMo ^1_I(R),\InjH ))} f\not \in
\goth p\leqno {(C_2)}$$ Condition $(C_1)$ is just a reformulation of
$\sqrt {(I)}=\sqrt {fR}$ (see the proof of remark 5.2.1). In
contrary to $(C_1)$, all objects in $(C_2)$ are graded and so
condition $(C_2)$ may be seen as a graded version of the condition
"f generates $I$ up to radical"; furthermore $(C_1)$ clearly implies
$(C_2)$.
\smallskip
Terminology: For a given homogenous ideal $I$ of $R$ and a given
element $f\in I$ we say that condition $(C_i)(I;f)$ holds if $(C_i)$
holds for $I$ and $f$ ($i=1,2$). \bigskip In the next section we
will investigate to what extent condition $(C_1)$ differs from
condition $(C_2)$. Theorem 5.3.5 will show that there are (in fact
many) non-homogenous $f\in I$ such that $(C_2)(I;f)$ holds, but
there are no non-homogenous $f\in I$ such that $(C_1)(I;f)$ holds.
\bigskip
{\bf 5.3.4 Remark} \par It is easy to see that for every homogenous
element $g\in I$ the conditions $(C_1)(I;f)$ and $(C_2)(I;f)$ are
equivalent.
\bigskip
{\bf 5.3.5 Theorem}
\smallskip
(i) Let $I$ be a homogenous ideal of a graded ring $R$ such that
$\ara ^h(I)\leq 1$. Let $g_1,\dots ,g_n\in I\setminus \{ 0\} $ be
homogenous of pairwise different degrees (in $\Naturalsign ^l$) and
such that
$$\sqrt I=\sqrt {(g_1,\dots ,g_n)R}\ \ .$$
Then
$$(C_2)(I;g_1+\dots +g_n)\hbox { holds.}$$
\par
(ii) Let $I$ be a homogenous ideal of a graded ring $R$ (and such
that $0=\LCMo ^2_I(R)=\LCMo ^3_I(R)=\dots $). Let $g\in I$ be a
non-homogenous element. Then
$$(C_1)(I;g) \hbox { does not hold.}$$
Proof: \par (i) We have $\LCMo ^1_I(R/(g_1,\dots ,g_n)R)=0$ and
hence
$$(g_1,\dots ,g_n)R\not\subseteq \goth p$$
for all $\goth p\in \Ass _R(*\Hom _R(\LCMo ^1_I(R),\InjH ))$.
Theorem 5.2.5 (ii) implies
$$((g_1,\dots ,g_n)R)^h\not \subseteq \bigcup _{\goth p\in \Ass _R(*\Hom
_R(\LCMo ^1_I(R),\InjH ))} \goth p\ \ .$$ Because of the different
degrees of the $g_i$ we conclude
$$(g_1+\dots +g_n)R\not \subseteq \bigcup _{\goth p\in \Ass _R(*\Hom
_R(\LCMo ^1_I(R),\InjH ))} \goth p$$ and the statement follows.
\par
(ii) The first statement of lemma 5.3.1 implies that if $R$ is a
graded domain and $I\subseteq R$ is a homogenous ideal such that
$\ara (I)\leq 1$ ($\iff \ara ^h(I)\leq 1$), every non-homogenous
$g\in I$ does not operate injectively on $\Hom _R(\LCMo
^1_I(R),\InjH )$. Furthermore, if $\ara (I)>1$ ($\iff \ara
^h(I)>1$), it is clear (use the ideas of section 5.2) that no $g\in
I$ operates injectively on $\Hom _R(\LCMo ^1_I(R),\InjH )$.
\bigskip
{\bf 5.3.6 Remark} \par While in the situation of theorem 5.3.5
statement (i) says there are (many) non-homogenous $f\in I$
operating injectively on $(*D)(\LCMo ^1_I(R))$, (ii) says there are
no non-homogenous $f\in I$ operating injectively on $D(\LCMo
^1_I(R))$. \vfil \eject {\bf 6 Applications}
\bigskip
\bigskip
{\bf 6.1 Hartshorne-Lichtenbaum vanishing}
\bigskip
The (more difficult) part of Hartshorne-Lichtenbaum vanishing
theorem (for another reference, see, e. g., [BS, 8.2.1]) says that
for an ideal $I$ in a noetherian local complete domain $(R,\goth m)$
there is the implication
$$\LCMo ^{\dim(R)}_I(R)\neq 0\Longrightarrow \sqrt I=\goth m\ \ .$$
We present two new proofs for it: theorem 6.1.2 works with the
normalization of $R$ and the Matlis dual of the local cohomology
module in question, while theorem 6.1.4 uses the fact that, over a
noetherian local complete Gorenstein ring $(S,\goth m)$ of dimension
$n+1$ and every height $n$ prime ideal $\goth P$ in $S$, one has
$D(\LCMo ^n_\goth P(S))=\widehat {S_\goth P}/S$ (this is lemma
3.2.1); it is remarkable that the proof of theorem 6.1.4 uses (this
is hidden in the proof of lemma 6.1.3) the ring structure on
$\widehat {S_\goth P}$.
\bigskip
{\bf 6.1.1 Theorem}
\par
Let $(R,\goth m)$ be a noetherian local ring and $M$ a finitely
generated $R$-module. Then $\LCMo ^{\dim _R(M)}_\goth m(M)\neq 0$.
\par
Proof:
\par
It is well-known and not difficult to see that for every $n\in {\bf
N} $, every ideal $I\subseteq R$ and every finitely generated
$R$-module $N$ the following statements are equivalent:
\par
(i) $H^i_I(N)=0$ for all $i\geq n$.
\par
(ii) $H^i_I(R/\goth p)=0$ for all $i\geq n$ and all $ \goth p\in
\Supp _R(N)$
\par
This fact implies (setting $d:=\dim _R(M)=\dim (R/\ann _R(M))$)
$\LCMo ^d_\goth m(M)\neq 0\iff \LCMo ^d_\goth m(R/\ann _R(M))\neq
0$. Thus we may assume that $M=R$ and $R$ is a domain. Again, we set
$d:=\dim (R)$ and choose a system of parameters $x_1,\dots ,x_d\in
R$ for $R$. Theorem 3.1.3 (ii) implies $\{ 0\} \in \Ass _R(D(\LCMo
^d_\goth m(R)))$; in particular, $\LCMo ^d_\goth m(R)\neq
0$.\bigskip {\bf 6.1.2 Theorem}
\par
Let $(R,\goth m)$ be a noetherian local complete equicharacteristic
domain, $n:=\dim (R)\geq 1$ and $I\subsetneq R$ an ideal. Then
$$\LCMo ^n_I(R)\neq 0\iff \sqrt I=\goth m$$
holds. \par
Proof:
\par
$\Longleftarrow $ follows from theorem 6.1.1. $\Longrightarrow $: By
induction on $n$, the case $n=1$ being trivial; we assume that $n>1$
and that the theorem is true for smaller $n$. Let $\tilde R$ be the
normalization of $R$. $\tilde R$ is a noetherian local (as $R$ is a
domain) complete equicharacteristic domain and is module-finite over
$R$, i. e. $\dim (\tilde R)=\dim (R)$; we denote the maximal ideal
of $\tilde R$ by $\goth m_{\tilde R}$. One has $\LCMo ^n_{I\tilde
R}(\tilde R)=\LCMo ^n_I(\tilde R)\neq 0$, because of $\Supp
_R(\tilde R)=\Spec (R)$. It suffices to show $\sqrt {I\tilde
R}=\goth m_{\tilde R}$ and so we may assume that $R$ is normal.
\smallskip
We choose $x_1,\dots ,x_n\in I$ such that $\sqrt {(x_1,\dots
,x_n)R}=\sqrt I$ and define the subring $R_0:=k[[x_1,\dots ,x_n]]$
of $R$, where $k$ is any fixed coefficient field of $R$; by $\goth
m_0$ we denote the maximal ideal of $R_0$. Because of $\LCMo
^n_{\goth m_0}(R)\neq 0$ we may conclude $\dim (R_0)=n$, i. e. $R_0$
is a formal power series ring over $k$ in the $n$ variables
$x_1,\dots ,x_n$. We set
$$x:=x_1,I_x:=\{ r\in R\vert \forall _{\varphi \in \Hom
_{R_0}(R,R_0)} \varphi (r)\in xR_0\} \ \ .$$ $I_x$ is an ideal of
$R$ such that $R\cdot x\subseteq I_x$; furthermore, we have
$I_x\subsetneq R$ because of
$$0\neq \Hom _{R_0}(\LCMo ^n_I(R),\InjH _{R_0}(k))\ \ .$$
For every $r\in I_x$ and every $\varphi \in \Hom _{R_0}(R,R_0)$ we
have $\im (r\cdot \varphi )\subseteq R_0\cdot x$ and thus there
exists $\varphi _0\in \Hom _{R_0}(R,R_0)$ such that $r\cdot \varphi
=x\cdot \varphi _0$. Therefore, $r\over x$ operates in a canonical
way on the finite $R$-module $\Hom _{R_0}(R,R_0)$ (note that $\LCMo
^n_I(R)$ is artinian as surjective image of the artinian $R$-module
$\LCMo ^n_\goth m(R)$; we conclude that $r\over x$ (as an element of
$Q(R)$, the quotient field of $R$) is integral over $R$; But $R$ is
normal and so we have $r\in R\cdot x$; this implies $I_x=R\cdot x$.
We have
$$\Ann _R(\Hom _{R_0/xR_0}(R/I_x,R_0/xR_0))=I_x=R\cdot x$$
and for every $\goth P\in \Ass _R(\Hom _{R_0/xR_0}(R/I_x,R_0/xR_0))$
there exists a non-trial $R_0/xR_0$-linear map $R/\goth P\to
R_0/xR_0$; by an easy Matlis duality argument, we conclude $\LCMo
^{n-1}_I(R/\goth P)\neq 0$ and, therefore, $\height (\goth P)=1$ and
$\sqrt {I+\goth P}=\goth m$ (induction hypothesis). We have shown
$$\Ass _R(\Hom _{R_0/xR_0}(R/I_x,R_0/xR_0))=\Min _R(R/xR)$$
and $\sqrt {I+\goth P}=\goth m$ for every $\goth P\in \Min
_R(R/xR)$. Now, because of
$$\bigcap _{\goth P\in \Min _R(R/xR)}\goth
P=\sqrt {xR}\subseteq \sqrt I$$
it follows that $\sqrt I=\goth m$.
\bigskip
{\bf 6.1.3 Lemma}
\par
Let $(S,\goth m)$ be a noetherian local regular ring, $\goth
P\subseteq S$ a prime ideal and $f\in \goth P$ an irreducible
element. Then $f$ operates injectively on $\widehat {S_\goth P}/S$.
\par
Proof: \par
We take a primary decomposition
$$f\widehat {S_\goth P}=\goth q_1\cap \dots \cap \goth q_n$$
of $f\widehat {S_\goth P}$ such that $\height(\goth q_i)=1$
($i=1,\dots ,n$). $f\in \goth P$ implies $n\geq 1$. Clearly we have
$f\in \goth q_1\cap S$ and $\goth q_1\cap S\subseteq S$ is a primary
ideal of height one (the canonical map $S\to \widehat {S_\goth P}$
is flat and thus going-down holds). On the other hand $fS\subseteq
S$ is a prime ideal of height one. Therefore we have
$$fS\supseteq
\goth q_1\cap S\supseteq f\widehat {S_\goth P}\cap S$$
and hence
$fS=f\widehat {S_\goth P}\cap S$. The statement follows.
\bigskip
{\bf 6.1.4 Theorem}
\par
Let $(R,\goth m)$ be a noetherian local complete domain and
$I\subseteq R$ an ideal such that $\sqrt I\subsetneq \goth m$. Then
$$\LCMo ^{\dim (R)}_I(R)=0$$
holds.
\par
Proof: \par Set $n:=\dim (R)$. For proper ideals $I_1\subseteq I_2$
of $R$ the canonical map $\LCMo ^n_{I_2}(R)\to \LCMo ^n_{I_1}(R)$ is
surjective and so we may assume $I$ is a prime ideal of $R$ of
height $n-1$. Now we choose $S,\goth Q$ and $\rho $ as in lemma
3.2.2. Then we have
$$\LCMo ^{\dim (R)}_I(R)=\LCMo ^{\dim (R)}_\goth
Q(S/fS)\otimes _{S/fS}R$$ ($f$ is a generator of the height one
prime ideal $\ker (\rho )$. Lemmas 3.2.1 and 6.3.1 imply that $f$
operates injectively on $D_S(\LCMo ^{\dim (R)}_\goth Q(S))$ and thus
the statement follows.
\bigskip
\bigskip
{\bf 6.2 Generalization of an example of Hartshorne}
\bigskip
The idea of this subsection is that, by Theorem 4.3.4, the Matlis
duals of certain local cohomology modules have infinitely many
associated prime ideals; but then this local cohomology module can
not be artinian. It turns out that this leads to a generalization of
an example of Hartshorne ([Ha1, section 3]); more details on this
generalization can be found in [HS2, section 1]. The author thanks
Gennady Lyubeznik for drawing his attention to this example.
\bigskip {\bf
6.2.1 Example}
\par
Let $k$ be a field, $R=k[[X_1,X_2,X_3,X_4]]$ a power series algebra
over $k$ in four variables, $I=(X_1,X_2)R$ and, for every $\lambda
\in k$, define
$$\goth p_\lambda :=(X_3+\lambda X_1,X_4+\lambda X_2)R\ \ .$$
Clearly, every $\goth p_\lambda $ is a height two prime ideal of $R$
and, by theorem 3.1.3 (v), is associated to $D=D(\LCMo ^2_I(R))$. On
the other hand, for every $\lambda \in k$, one has
$$p:=X_1X_4+X_2X_3\in \goth p_\lambda $$
(because of $p=X_1(X_4-\lambda X_2)+X_2(X_3+\lambda X_1)$).
Therefore, at least if $k$ is infinite, $D$ has infinitely many
associated primes containing $p$. This implies that
$$\Hom _R(R/pR,D)$$
cannot be finitely generated. But $\Hom _R(R/pR,D)$ is the Matlis
dual of
$$\LCMo ^2_I(R/pR)$$
and so $\LCMo ^2_I(R/pR)$ cannot be artinian.
\bigskip
{\bf 6.2.2 Remark}
\par
This is essentially Hartshorne's example ([Ha1, section 3]), the
main difference is that Hartshorne works over a ring of the form
$k[X_3,X_4][[X_1,X_2]]$, while we work over a ring of the form
$k[[X_1,X_2,X_3,X_4]]$; but the two versions are essentially the
same, because the module
$$\LCMo ^2_{(X_1,X_2)}(k[X_3,X_4][[X_1,X_2]]/(X_1X_4+X_2X_3))$$
is naturally a module over $k[[X_1,X_2,X_3,X_4]]$, because its
support is $\{ (X_1,X_2,X_3,X_4)\} $. This is true, because for
every prime ideal $\goth p\neq (X_1,X_2,X_3,X_4)$ of
$k[X_3,X_4][[X_1,X_2]]$ containing $X_1X_4+X_2X_3$ the ring
$(k[X_3,X_4][[X_1,X_2]]/(X_1X_4+X_2X_3))_\goth p$ is regular, and so
Hartshorne-Lichtenbaum vanishing shows that $\LCMo
^2_{(X_1,X_2)}(k[X_3,X_4][[X_1,X_2]]/(X_1X_4+X_2X_3))_\goth p=0$.
\bigskip A similar technique like
in the example above works to show that $\LCMo ^{n-2}_I(R/pR)$ is
not artinian in the general situation $R=k[[X_1,\dots ,X_n]]$,
$n\geq 4$, $I=(X_1,\dots ,X_{n-2})R$ and $p\in R$ a prime element
such that $p\in (X_{n-1},X_n)R$, even if the field $k$ is finite:
\bigskip
{\bf 6.2.3 Theorem}
\par
Let $k$ be a field, $n\geq 4$, $R=k[[X_1,\dots ,X_n]]$,
$I=(X_1,\dots ,X_{n-2})R$ and $p\in R$ a prime element such that
$p\in (X_{n-1},X_n)R$. Then $\LCMo ^{n-2}_I(R/pR)$ is not artinian.
\par
Proof: \par Set $D:=D(\LCMo ^{n-2}_I(R))$. If $p\not\in I$, it is
easy to see that $$\Supp _R(\LCMo ^{n-2}_I(R/pR))={\cal V}(I+pR)\ \
,$$ the set of prime ideals of $R$ containing $I+pR$, and so $\LCMo
^{n-2}_I(R/pR)$ is not artinian (it is not zero-dimensional). We
assume $p\in I$: If $\LCMo ^{n-2}_I(R/pR)$ was artinian, $D(\LCMo
^{n-2}_I(R/pR))$ would be finitely generated; but we have seen
before that, because of the exactness of $D$ and the right-exactness
of $\LCMo ^{n-2}_I$,
$$D(\LCMo ^{n-2}_I(R/pR))=\Hom _R(R/pR,D)\ \ ,$$
and from Theorem 4.3.4 we
know that the latter module is not finitely generated (it has
infinitely many associated prime ideals).
\bigskip
{\bf 6.2.4 Remark} Marley and Vassilev have shown \smallskip {\bf
Theorem} ([MV, theorem 2.3])
\par
Let $(T,\goth m)$ be a noetherian local ring of dimension at least
two. Let $R=T[x_1,\dots ,x_n]$ be a polynomial ring in $n$ variables
over $T$, $I=(x_1,\dots ,x_n)$, and $f\in R$ a homogenous polynomial
whose coefficients form a system of parameters for $T$. Then the
*socle of $\LCMo ^n_I(R/fR)$ is infinite dimensional.
\smallskip
In their paper [MV], Marley and Vassilev say (in section 1) that
Hartshorne's example is obtained by letting $T=k[[u,v]],n=2$ and
$f=ux+vy$; there is a slight difference between the two situations
that comes from the fact that Hartshorne works over a ring of the
form $k[x,y][[u,v]]$ while Marley and Vassilev work over a ring of
the form $k[[u,v]][x,y]$. The two rings are not the same. But, as
$$\Supp _R(\LCMo ^2_{(u,v)}(R/(uy+vx)))=\{ (x,y,u,v)\} $$
(both for $R=k[x,y][[u,v]]$ and for $R=k[[u,v]][x,y]$), the local
cohomology module in question is (in both cases) naturally a module
over $k[[x,y,u,v]]$ and, therefore, both versions are equivalent, i.
e. the result of Marley and Vassilev is a generalization of
Hartshorne's example.
\bigskip
{\bf 6.2.5 Remark}
\par
[MV, theorem 2.3] and theorem 6.2.3 are both generalizations of
Hartshorne's example, but, due to different hypotheses, they can
only be compared in the following special case: $k$ a field, $n\geq
4$,
$$R_0=k[[X_{n-1},X_n]][X_1,\dots ,X_{n-2}]\ \ ,$$
$$R=k[[X_1,\dots ,X_n]]\ \ ,$$
$I=(X_1,\dots ,X_{n-2})R$, $p\in R_0$ a homogenous element such that
$p$ is prime as an element of $R$. Then [MV, theorem 2.3] says
(implicitly) that
$$\LCMo ^{n-2}_I(R/pR)$$
is not artinian, if the coefficients of $p\in R_0$ in
$k[[X_{n-1},X_n]]$ form a system of parameters in
$k[[X_{n-1},X_n]]$, while theorem 6.2.3 says that the same module is
not artinian if none of these coefficients of $p$ is a unit in
$k[[X_{n-1},X_n]]$.
\bigskip
\bigskip
{\bf 6.3 A necessary condition for set-theoretic complete
intersections}
\bigskip
Let $(R,\goth m)$ be a noetherian local ring and $I=(x_1,\dots
,x_i)R=I\subseteq R$ a set-theoretic complete intersection ideal (in
the sense that its height is $i$). Then $\LCMo ^i_I(R)\neq 0$ (this
can be seen by localizing at a height-$i$ prime ideal of $R$
containing $I$). On the other hand, statement (ii) from theorem
6.3.1 below presents a necessary condition for $\LCMo ^i_I(R)\neq
0$.
\bigskip {\bf 6.3.1 Theorem}
\par
Let $(R,\goth m)$ be a noetherian local complete domain containing a
field $k$ and $x_1,\dots ,x_i$ a sequence in $R$ ($i\geq 1$). Define
$R_0:=k[[x_1,\dots ,x_i]]$ as a subring of $R$.
\par
(i) The following two statements are equivalent:
\par
($\alpha $) $\LCMo ^i_{(x_1,\dots ,x_i)R}(R)\neq 0$.
\par
($\beta $) $\Hom _{R_0}(R,R_0)\neq 0$ and $\dim (R_0)=i$.
\par
(ii) If the equivalent statements of (i) hold, one has
$$R\cap Q(R_0)=R_0\ \ ,$$
where $Q(R_0)$ denotes the quotient field of $R_0$ and where the
intersection is taken in the quotient field $Q(R)$ of $R$.
\par
Proof:
\par
We denote the maximal ideal of $R_0$ by $\goth m_0$. \par (i) We
have
$$\LCMo ^i_{(x_1,\dots ,x_i)R}(R)=R\otimes _{R_0} \LCMo ^i_{(x_1,\dots
,x_i)R_0}(R_0)$$ and, thus, $\LCMo ^i_{(x_1,\dots ,x_i)R_0}(R_0)\neq
0$ implies $\dim (R_0)=i$; therefore we may assume that $R_0$ is a
formal power series ring in $x_1,\dots ,x_i$. Therefore, we may
assume that $R_0$ is $i$-dimensional. Let $\InjH _{R_0}(k)$ be a
fixed $R_0$-injective hull of $R_0$. We have
$$\eqalign {\Hom _{R_0}(\LCMo ^i_{(x_1,\dots ,x_i)R}(R),\InjH
_{R_0}(k))&=\Hom _{R_0}(R\otimes _{R_0} \LCMo ^i_{(x_1,\dots
,x_i)R_0}(R_0),\InjH _{R_0}(k))\cr &=\Hom _{R_0}(R,\Hom _{R_0}(\LCMo
^i_{(x_1,\dots ,x_i)R_0}(R_0),\InjH _{R_0}(k)))\cr &=\Hom
_{R_0}(R,R_0)\ \ ,\cr }$$ where we have used the fact that $R_0$ is
a formal power series ring in $x_1,\dots ,x_i$ over $k$. This
identity shows the stated equivalence.
\par
(ii) Under the given assumptions, we have $\Hom _{R_0}(R,R_0)\neq
0$. Let
$$\varphi \in \Hom _{R_0}(R,R_0)$$
be any non-zero element and let $r_0\in R_0\setminus \{ 0\} ,r\in R$
such that $r_0\cdot r\in R_0$. We have to show $r\in R_0$: We set
$$r_0^\prime :=r_0r$$
and conclude
$$r_0\varphi (r)=\varphi (r_0^\prime )=r_0^\prime \varphi (1)\ \ .$$
This shows
$$\varphi (1)r=\varphi (1){r_0^\prime \over r_0}=\varphi (r)\in R_0\
\ .$$ On the other hand, we have
$$r_0^{\prime 2}=r_0^2r^2$$
and thus
$$r_0^2\varphi (r^2)=r_0^{\prime 2}\varphi (1)$$
and
$$\varphi (1)r^2=\varphi (1){r_0^{\prime 2}\over r_0^2}=\varphi
(r^2)\in R_0\ \ .$$ Continuing in the same way, one sees that, for
every $l\geq 1$, one has
$$\varphi (1)r^l\in R_0\ \ .$$
But this implies that the $R_0$-module
$$\varphi (1)\cdot <1,r,r^2,\dots >_{R_0}$$
is finitely generated ($<1,r,r^2,\dots >_{R_0}$ stands for the
$R_0$-submodule of $R$ generated by $1,r,r^2,\dots $). But, as $R$
is a domain,
$$<1,r,r^2,\dots >{R_0}$$
is then finitely generated, too, i. e. $r$ is necessarily contained
in $R_0$.
\bigskip
\bigskip
{\bf 6.4 A generalization of local duality}
\bigskip
Over some rings (e. g. over complete Cohen-Macaulay rings), there is
a correspondence between certain $\Ext $-modules on the one hand and
certain local cohomology modules on the other hand; this
correspondence is given (in both directions) by taking the Matlis
dual and is called local duality. This result can e. g. be found in
[BS, section 11]. In the form in which it is usually presented,
local duality works only if the support ideal is $\goth m$, i. e. if
one takes local cohomology with support in $\goth m$. But, below we
generalize this result to a large class of support ideals $I$.
\bigskip
{\bf 6.4.1 Theorem}
\bigskip
Let $(R,\goth m)$ be a noetherian local ring, $I\subseteq R$ an
ideal, $h\in \Naturalsign $ such that
$$\LCMo ^l_I(R)\neq 0\iff l=h$$
holds and $M$ an $R$-module. Then, for every $i\in \{ 0,\dots ,h\}
$, one has
$$\Ext ^i_R(M,D(\LCMo ^h_I(R)))=D(\LCMo ^{h-i}_I(M))\ \ .$$
Proof:
\par
We take the sequence of functors $(D\circ \LCMo ^{h-i}_I)_{i\in
\Naturalsign }$ from the category of $R$-modules to itself; of
course, $\LCMo ^M_I=0$ for every $M<0$. Given a short exact sequence
of $R$-modules
$$0\to M^\prime \to M\to M^{\prime \prime }\to 0\ \ ,$$
we clearly get an exact sequence of the form
$$\eqalign {0&\to D(\LCMo ^h_I(M^{\prime \prime }))\to D(\LCMo
^h_I(M))\to D(\LCMo ^h_I(M^\prime ))\to \cr &\to D(\LCMo
^{h-1}_I(M^{\prime \prime }))\to \dots \cr }$$ (note that, by our
hypothesis, $\LCMo ^{h+1}_I(M^\prime )=M^\prime \otimes _R\LCMo
^{h+1}_I(R)=0$). In the case $i=0$ we get, for any $R$-module $M$,
$$\eqalign {D(\LCMo ^h_I(M))&=\Hom _R(\LCMo ^h_I(M),\InjH _R(R/\goth
m))\cr &=\Hom _R(M\otimes _R\LCMo ^h_I(R),\InjH _R(R/\goth m))\cr
&=\Hom _R(M,\Hom _R(\LCMo ^h_I(R),\InjH _R(R/\goth m)))\cr &=\Hom
_R(M,D(\LCMo ^h_I(R)))\ \ .\cr }$$ Finally, for every $i>0$ and
every $m\in \Naturalsign $, we have $\LCMo ^{h-i}_I(R^m)=0$ and
hence $\LCMo ^{h-i}_I(F)=0$ for every free $R$-module $F$; we get
$$D(\LCMo ^{h-i}_I(F))=0$$
for every $i>0$ and every free $R$-module $F$. By some well-known
homology theory, the last three properties imply our statement.
\bigskip
{\bf 6.4.2 Remark}
\par
If, in the situation of theorem 6.4.1, $R$ is complete and
Cohen-Macaulay and $I=\goth m$ (then $h=\dim (R)$ necessarily), the
statement takes the form
$$\Ext ^i_R(M,D(\LCMo ^{\dim (R)}_\goth m(R)))=D(\LCMo ^{\dim
(R)-i}_\goth m(M))$$ for every $R$-module $M$ and every $i\in \{
0,\dots ,\dim (R)\} $. If we assume furthermore that $M$ is finitely
generated, then $\LCMo ^{\dim (R)-i}_\goth m(M)$ is artinian and the
above statement implies (in fact, is equivalent to)
$$D(\Ext ^i_R(M,D(\LCMo ^{\dim (R)}_\goth m(R))))=\LCMo ^{\dim
(R)-i}_\goth m(M)\ \ .$$ We study the $R$-module $D(\LCMo ^{\dim
(R)}_\goth m(R))$: By Matlis duality, it is finitely generated. We
calculate its type, which is defined as the following $R/\goth
m$-vector space dimension:
$$\dim _{R/\goth m}(\Ext ^{\dim (R)}_R(R/\goth m,D(^{\dim (R)}_\goth
m(R))))=\dim _{R/\goth m} (D(\LCMo ^0_\goth m(R/\goth m)))=1$$ (note
that, for the first equality, we use theorem 6.4.1 again). The
$\goth m$-depth of $D(\LCMo ^{\dim (R)}_\goth m(R))$ is $\dim (R)$
(this follows from theorem 1.1.2, take any parameter sequence
$\underline x$ of $R$, it will be a regular sequence on $D(\LCMo
^{\dim (R)}_\goth m(R))$), i. e. $D(\LCMo ^{\dim (R)}_\goth m(R))$
is a maximal Cohen-Macaulay module (over $R$). By definition, these
properties show that
$$D(\LCMo ^{\dim (R)}_\goth m(R))=:\omega _R$$
is a canonical module for $R$, and the statement of theorem 6.4.1
becomes
$$D(\Ext ^i_R(M,\omega _R))=\LCMo ^{\dim (R)-i}_\goth m(M)\ \ .$$
Clearly, this is a version of the local duality theorem (see, e. g.,
[BS, section 11] for more details on local duality). \vfil \eject
{\bf 7 Further Topics}
\bigskip
\bigskip
{\bf 7.1 Local Cohomology of formal schemes} \bigskip In some cases
we can consider the Matlis duals of local cohomology modules as
certain local cohomology modules of the structure sheaf of some
formal scheme (see [Og, in particular section 2]), here are the
details:
\bigskip {\bf 7.1.1 Remark}
\par
Let $(R,\goth m)$ be a noetherian local Gorenstein ring and $I$ an
ideal of $R$. We define
$$X:=\Spec (R), Y:={\cal V}(I)\ \ ,$$
i. e. $Y$ is the closed subscheme of $X$ defined by $I$. We denote
by $\cal X$ the formal completion of $X$ along $Y$ and by $p$ the
closed point of the topological space underlying $\cal X$ (note that
as topological spaces $\cal X$ and $X$ are the same). Then, for
every $i\in \Naturalsign $, there is a canonical isomorphism
$$\LCMo ^i_p({\cal X},{\cal O_X})=D(\LCMo ^{\dim (R)-i}_I(R))$$
(see [Og, 2.2.3]). This follows essentially from local duality, the
fact that $\LCMo ^i_\goth m(R/I^v)$ is artinian for every $v\geq 1$
and the existence of a short exact sequence
$$0\to R^1\invlim _v \LCMo ^{i-1}_\goth m(R/I^v)\to \LCMo ^i_p({\cal
X},{\cal O_X})\to \invlim _v\LCMo ^i_\goth m(R/I^v)\to 0\ \ .$$ Note
that local cohomology of a (formal) coherent sheaf on a formal
scheme is defined in the sense of Grothendieck, i. e. as (local)
cohomology of an abelian sheaf on a topological space.
\bigskip
\bigskip
{\bf 7.2 $D(\LCMo ^i_I(R))$ has a natural $D$-module structure}
\bigskip
Let $k$ be a field and $R=k[[X_1,\dots X_n]]$ a power series ring
over $k$ in $n$ variables. Let
$$D(R,k)\subseteq \End _k(R)$$
be the (non-commutative) subring defined by the multiplication maps
by $r\in R$ (for all $r\in R$) and by all $k$-linear derivation maps
from $R$ to $R$. $D:=D(R,k)$ is the so-called ring of $k$-linear
differential operators on $R$. [Bj] contains material on the ring
$D(R,k)$ and on similar rings; $D$-modules in relation with local
cohomology modules have been studied in [Ly1]. For $i=1,\dots ,n$
let $\partial _i$ denote the partial derivation map from $R$ to $R$
with respect to $X_i$. Then, as an $R$-module, one has
$$D(R,k)=\oplus _{i_1,\dots
,i_n\in \Naturalsign } R\cdot \partial _1^{i_1}\dots \partial
_n^{i_n}\ \ .\leqno{(1)}$$ Now, let $I\subseteq R$ be an ideal and
$i\in \Naturalsign $. We will demonstrate that there is a canonical
left-$D$-module structure on $D(\LCMo ^i_I(R))$ (the following idea
was inspired by Gennady Lyubeznik). To do so, by identity (1), it is
sufficient to determine the action of an arbitrary $k$-linear
derivation $\delta :R\to R$ on $D(\LCMo ^i_I(R))$, to extend it to
an action of $D(R,k)$ on $D(\LCMo ^i_I(R))$ and to show that this
action is well-defined and satisfies all axioms of a
left-$D$-module. The derivation $\delta $ induces a $k$-linear map
$$R/I^v\to R/I^{v-1}\ \ (v\geq 1)$$
and, in a canonical way, a map of complexes from the \v{C}ech
complex of $R/I^v$ with respect to $X_1,\dots ,X_n$ to the \v{C}ech
complex of $R/I^{v-1}$ with respect to $X_1,\dots ,X_n$ ($v\geq 1$).
By taking cohomology, we get a map
$$\LCMo ^{n-i}_\goth m(R/I^v)\to \LCMo ^{n-i}_\goth m(R/I^{v-1})\ \
(v\geq 1)\ ,$$ where $\goth m$ stands for the maximal ideal of $R$.
These maps induce a map
$$\invlim _{v\in \Naturalsign }(\LCMo ^{n-i}_\goth m(R/I^v))\to \invlim _{v\in \Naturalsign }(\LCMo ^{n-i}_\goth
m(R/I^v))$$ (note that the maps of the above inverse system are
induced by the canonical epimorphisms $R/I^v\to R/I^{v-1}$). But, by
local duality and $\LCMo ^i_I(R)=\dirlim _{v\in \Naturalsign }(Ext
^i_R(R/I^v,R))$, one has
$$\invlim _{v\in \Naturalsign }(\LCMo ^{n-i}_\goth m(R/I^v))=D(\dirlim
_{v\in \Naturalsign }(Ext ^i_R(R/I^v,R)))=D(\LCMo ^i_I(R))\ \ .$$
Now, having determined the action of the element $\delta $ on
$D(\LCMo ^i_I(R))$, by (1) it is clear how to extend this to an
action of $D(R,k)$ on $D(\LCMo ^i_I(R))$ such that $D(\LCMo
^i_I(R))$ becomes a left-$D$-module (note that, for every $k$-linear
derivation $\delta :R\to R$ and every $r\in R$, we have $\delta
(r\cdot d)=\delta (r)\cdot d+r\cdot \delta (d)$, i. e. the action of
$D$ on $D(\LCMo ^i_I(R)$ makes it a left-$D$-module). We have seen
(in various situations) in sections 2, 3 and 4, that, in general,
$$D(\LCMo ^i_I(R))$$
has infinitely many associated primes. On the other hand, one knows
from [Ly1, Theorem 2.4 (c)] (at least if $\char (k)=0$), that every
finitely generated left-$D$-module has only finitely many associated
prime ideals (as $R$-module, of course). This shows that, in
general, $D(\LCMo ^i_I(R))$ is an example of a non-finitely
generated left-$D$-module. In particular, $D(\LCMo ^i_I(R))$ is not
holonomic in general (see [Bj] for the notion of holonomic modules).
\bigskip
\bigskip
{\bf 7.3 The zeroth Bass number of $D(\LCMo ^i_I(R))$ (w. r. t. the
zero ideal) is not finite in general}
\bigskip
Let $(R,\goth m)$ be a noetherian local domain, $i\geq 1$ and
$x_1,\dots ,x_i\in R$. Then, as we have seen in theorem 3.1.3 (ii),
one has
$$\{ 0\} \in \Ass _R(D(\LCMo ^i_{(x_1,\dots ,x_i)R}(R)))$$
in some situations; actually, if conjecture (*) holds, this is true
provided $\LCMo ^i_{(x_1,\dots ,x_i)R}(R)\neq 0$ holds. It is
natural to ask for the associated Bass number of $D(\LCMo
^i_{(x_1,\dots ,x_i)R}(R))$, i. e. the $Q(R)$-vector space dimension
of
$$D(\LCMo ^i_{(x_1,\dots ,x_i)R}(R))\otimes _RQ(R)\ \ ,$$
where $Q(R)$ stands for the quotient field of $R$. As we will see
below, this number is not finite in general; more precisely, we
consider the following case: Let $k$ be a field, $R=k[[X_1,\dots
,X_n]]$ a power series algebra over $k$ in $n\geq 2$ variables,
$1\leq i<n$ and $I$ the ideal $(X_1,\dots ,X_i)R$ of $R$; in this
situation
$$\dim _{Q(R)}(D(\LCMo ^i_I(R))\otimes _RQ(R))=\infty $$
holds, see theorem 7.3.2 below for a proof. \vfil \eject {\bf 7.3.1
Remark}
\par
Note that in section 4.3 (see, in particular, the proof of remark
4.3.6) we introduced some notation on polynomials in "inverse
variables" and we explained and proved (note that the situation here
is more general then in remark 4.3.6, where $i$ was $n-2$, but the
proof of 4.3.6 works in this more general situation too) the
following formulas:
$$\LCMo ^i_I(R)=k[[X_{i+1},\dots ,X_n]][X_1^{-1},\dots ,X_i^{-1}]\ \ ,$$
$$\InjH _R(k)=k[X_1^{-1},\dots ,X_n^{-1}]$$
and
$$D(\LCMo ^i_I(R))=k[X_{i+1}^{-1},\dots ,X_n^{-1}][[X_1,\dots
,X_i]]\ \ .$$ Also note that the latter module is different from and
larger than the module
$$k[[X_1,\dots ,X_i]][X_{i+1}^{-1},\dots ,X_n^{-1}]\ \ .$$
\bigskip
The following proof is technical; its basic idea is the following
one: Let $k$ be a field, $R=k[[X,Y]]$ a power series algebra over
$k$ in two variables; then we have
$$\LCMo ^1_{XR}(R)=k[[Y]][X^{-1}]$$
and $$D:=D(\LCMo ^1_{XR}(R))=k[Y^{-1}][[X]]\ \ .$$ Set
$$\eqalign {d_2&:=\sum _{l\in \Naturalsign }Y^{-l^2}X^l\cr
&=1+Y^{-1}X+Y^{-4}X^2+Y^{-9}X^3+\dots \in D\cr }$$ and let $r\in
R\setminus \{ 0\} $ be arbitrary. Because of $r\neq 0$ we can write
$$r=X^{a+1}\cdot h+X^a\cdot g$$
with some $h\in R, g\in k[[Y]]\setminus \{ 0\} $. Then, at least for
$l>>0$, the coefficient of $r\cdot d_2$ in front of $X^l$ is
$$h^*\cdot Y^{-(l-a-1)^2}+g\cdot Y^{-(l-a)^2}$$
for some $h^*\in k[[Y]]$. Now, if we write
$$g=c_bY^b+c_{b+1}Y^{b+1}+\dots $$
for some $b\in \Naturalsign , c_b\neq 0$ and observe the fact
$$-(l-a)^2+b<-(l-a-1)^2 \ \ (l>>0)\ \ ,$$
it follows that the term
$$c_b\cdot Y^{-(l-a)^2+b}$$
(coming from $h^*\cdot Y^{-(l-a-1)^2}+g\cdot Y^{-(l-a)^2}$) cannot
be canceled out by any other term. In fact, for $l>>0$, the lowest
non-vanishing $Y$-exponent of the coefficient in front of $X^l$, is
$-(l-a)^2+b$. The crucial point is that the sequences $-(l-a)^2+b$
and $-l^2$ agree up to the two shifts given by $a$ and $b$. This
means that some information about $d_2$ is stored in $rd_2$.
\bigskip
{\bf 7.3.2 Theorem}
\par
Let $k$ be a field, $R=k[[X_1,\dots ,X_n]]$ a power series algebra
over $k$ in $n\geq 2$ variables, $1\leq i<n$ and $I$ the ideal
$(X_1,\dots ,X_i)R$ of $R$. Then
$$\dim _{Q(R)}(D(\LCMo ^i_I(R))\otimes _RQ(R))=\infty $$
holds.
\par
Proof:
\par
As the proof is technical we will first show the case $n=2,i=1$; in
the remark after this proof we will explain how one can reduce the
general to this special case. Set $X=X_1,Y=X_2$ and
$$D:=D(\LCMo ^2_I(R))=k[Y^{-1}][[X]]\ \ .$$
For every $n\in \Naturalsign \setminus \{ 0\} $, set
$$d_n:=\sum _{l\in \Naturalsign }Y^{-l^n}\cdot X^l\in D\ \ .$$
It is sufficient to show the following statement: The elements
$(d_n\otimes 1)_{n\in \Naturalsign \setminus \{ 0\} }$ in $D\otimes
_RQ(R)$ are $Q(R)$-linear independent:
\bigskip
We define an equivalence relation on ${\bf Z}^\Naturalsign $ (the
set of all maps from $\Naturalsign $ to $\bf Z$, i. e. infinite
sequences of integers) by saying that $(a_n),(b_n)\in {\bf
Z}^\Naturalsign $ are equivalent (short form: $(a_n)\sim (b_n)$) iff
there exist $N,M\in \Naturalsign $ and $p\in {\bf Z}$ such that
$$a_{N+1}=b_{M+1}+p,a_{N+2}=b_{M+2}+p,\dots $$
hold. It is easy to see that $\sim $ is an equivalence relation on
${\bf Z}^\Naturalsign $. For every $d\in D$, we define $\delta
(d)\in {\bf Z}^\Naturalsign $ in the following way: Let $f_l\in
k[Y^{-1}]$ be the coefficient of $d$ in front of $X^l$; we set
$$(\delta (d))(l):=0$$
if $f_l=0$ and
$$(\delta (d))(l):=s$$
if $s$ is the smallest $Y$-exponent of $f_l$, i. e.
$$f_l=c_sY^s+c_{s+1}Y^{s+1}+\dots +c_0\cdot 1$$
for some $c_s\neq 0$.
\bigskip
Now suppose that $r_1,\dots ,r_{n_0}\in R$ are given such that
$r_{n_0}\neq 0$. We claim that
$$\delta (r_1d_1+\cdot +r_{n_0}d_{n_0})\sim \delta (d_{n_0})$$
holds. Note that if we prove this statement we are done, essentially
because then $r_1d_1+\dots +r_{n_0}d_{n_0}$ can not be zero.\bigskip
It is obvious that one has $\delta (d+d^\prime )\sim \delta
(d_{N_2})$ for given $d,d^\prime \in D$ such that
$$\delta (d)\sim \delta (d_{N_1}), \delta (d^\prime )\sim \delta
(d_{N_2}), N_2>N_1\ \ .$$ For this reason it is even sufficient to
prove the following statement: For a fixed $n\in \Naturalsign
\setminus \{ 0\} $ and for any $r\in R\setminus \{ 0\} $ one has
$$\delta (rd_n)\sim \delta (d_n)\ \ .$$
We can write
$$r=X^{a+1}\cdot h+X^a\cdot g$$
with $a\in \Naturalsign , h\in k[[X,Y]]$ and $g\in k[[Y]]\setminus
\{ 0\} $. We get
$$\delta (r\cdot d_n)\sim \delta (\sum _{l\geq a+1}
(hY^{-(l-a-1)^n}+gY^{-(l-a)^n})X^l)$$ and we write
$$g=c_bY^b+c_{b+1}Y^{b+1}+\dots $$
with $c_b\in k^*$. Now, because of
$$-(l-a)^n+b<-(l-a-1)^n\ \ (l>>0)$$
it is clear that, for $l>>0$, the smallest $Y$-exponent in front of
$X^l$ (of the power series $r\cdot d_n$) is $-(l-a)^n+b$. Therefore,
one has
$$\delta (r\cdot d_n)\sim (-l^n)\sim \delta (d_n)$$
and we are done.
\bigskip
{\bf 7.3.3 Remark}
\par
A proof of the general case of theorem can be obtained e. g. in the
following way: First, we use theorem 3.1.2 repeatedly to get a
surjection
$$\LCMo ^i_{(X_1,\dots ,X_i)R}(R)\to \LCMo ^{n-1}_{(X_1,\dots
,X_{n-1})R}(R)$$ and hence an injection
$$D(\LCMo ^{n-1}_{(X_1,\dots
,X_{n-1})R}(R))\to D(\LCMo ^i_{(X_1,\dots ,X_i)R}(R))\ \ ,$$ which
allows us to reduce to the case $i=n-1$; then it is possible to
adapt our proof of theorem 7.3.2 with some minor changes: Instead of
working with maps $\Naturalsign \to {\bf Z}$, one works with maps
$$\Naturalsign ^{n-1}\to {\bf Z}$$
and also with multi-indices instead of indices.
\bigskip
\bigskip
{\bf 7.4 On the module $\LCMo ^h_I(D(\LCMo ^h_I(R)))$}
\bigskip
In the previous section we were interested in modules of the form
$$D:=D(\LCMo ^i_I(R))\ \ ,$$
where $I$ is an ideal in a local ring $R$. In this section we
compute the local cohomology module
$$\LCMo ^i_I(D)\ \ .$$
Our results say (essentially) that this module is $\InjH _R(R/\goth
m)$ if $I$ is a set-theoretic complete intersection and it is either
$\InjH _R(R/\goth m)$ or zero in general (see theorems 7.4.1 and
7.4.2 for precise formulations and proofs).
\bigskip
{\bf 7.4.1 Theorem}
\par
Let $(R,\goth m)$ be a noetherian local complete Cohen-Macaulay ring
with coefficient field $k$ and $x_1,\dots ,x_i\in R$ ($i\geq 1$) a
regular sequence in $R$. Set $I:=(x_1,\dots ,x_i)R)$ ($I$ is a
set-theoretic complete intersection ideal of $R$). Then one has
$$\LCMo ^i_I(D(\LCMo ^i_I(R)))=\InjH _R(k)\ \ .$$
Proof:
\par
First we show a special case: Assume that $R=k[[X_1,\dots ,X_n]]$ is
a formal power series algebra over $k$ in $n$ variables and
$x_1=X_1,\dots ,x_i=X_i$. Then, as we have seen in the proof of
remark 4.3.6 (note that the situation here is more general then in
remark 4.3.6, where $i$ was $n-2$, but the proof of 4.3.6 works in
this more general situation too), we have
$$\LCMo ^i_I(R)=k[[X_{i+1},\dots ,X_n]][X_1^{-1},\dots ,X_i^{-1}]$$
and
$$D(\LCMo ^i_I(R))=k[X_{i+1}^{-1},\dots ,X_n^{-1}][[X_1,\dots
,X_i]]$$ (again, see section 4.3, in particular the proof of remark
4.3.6 for the notation). As the functor $\LCMo ^i_I$ is right-exact,
we have
$$\eqalign {\LCMo ^i_I(D(\LCMo ^i_I(R)))&=\LCMo ^i_I(R)\otimes _RD(\LCMo
^i_I(R))\cr &=k[[X_{i+1},\dots ,X_n]][X_1^{-1},\dots
,X_i^{-1}]\otimes _Rk[X_{i+1}^{-1},\dots X_n^{-1}][[X_1,\dots
,X_i]]\cr &\buildrel (*)\over =k[X_1^{-1},\dots ,X_n^{-1}]\cr
&=\InjH _R(k)\cr }$$ Proof of equality (*): The map
$$k[X_{i+1},\dots ,X_n][X_1^{-1},\dots
,X_i^{-1}]\otimes k[X_{i+1}^{-1},\dots X_n^{-1}][X_1,\dots ,X_i]\to
k[X_1^{-1},\dots ,X_n^{-1}]$$
$$X_{i+1}^{r_{i+1}}\cdot \dots \cdot X_n^{r_n}\cdot X_1^{-s_1}\cdot
\dots X_i^{-s_i}\otimes X_{i+1}^{-t_{i+1}}\cdot \dots \cdot
X_n^{-t_n}\cdot X_1^{u_1}\cdot \dots \cdot X_i^{u_i}\mapsto $$
$$\mapsto X_{i+1}^{r_{i+1}-t_{i+1}}\cdot X_n^{r_n-t_n}\cdot
X_1^{u_1-s_1}\cdot \dots \cdot X_i^{u_i-s_i}\hbox { if
}r_{i+1}-t_{i+1},\dots ,r_n-t_n,u_1-s_1,\dots ,u_i-s_i\leq 0$$ and
to zero otherwise, induces an $R$-linear map
$$k[[X_{i+1},\dots ,X_n]][X_1^{-1},\dots
,X_i^{-1}]\otimes _Rk[X_{i+1}^{-1},\dots X_n^{-1}][[X_1,\dots
,X_i]]\to k[X_1^{-1},\dots ,X_n^{-1}]\ \  ,\leqno {(1)}$$ which is
surjective and maps the $k$-vector space generating system
$$\{ X_1^{-s_1}\cdot \dots \cdot X_i^{-s_i}\otimes
X_{i+1}^{-t_{i+1}}\cdot \dots \cdot X_n^{-t_n}\vert s_1,\dots
,s_i,t_{i+1},\dots ,t_n\geq 0\} $$ of the vector space on the left
side of (1) to the $k$-basis
$$\{ X_1^{-s_1}\cdot \dots \cdot X_i^{-s_i}\cdot
X_{i+1}^{-t_{i+1}}\cdot \dots \cdot X_n^{-t_n}\vert s_1,\dots
,s_i,t_{i+1},\dots ,t_n\geq 0\} $$ of the vector space on the right
side of (1), and therefore provides us with the desired isomorphism
in our special case.
\par
We come to the general case: Choose $x_{i+1},\dots ,x_n\in R$ such
that $\sqrt {(x_1,\dots ,x_n)R}=\goth m$ ($x_1,\dots x_n$ is a s. o.
p. of $R$). Define
$$R_0:=k[[x_1,\dots ,x_n]]\subseteq R$$
$R_0$ is regular of dimension $n$ and $R$ is a finite-rank free
$R_0$-module. Define $I_0:=(x_1,\dots ,x_i)R_0$. We have
$$\LCMo ^i_I(D(\LCMo ^i_I(R)))=\LCMo ^i_I(R)\otimes _RD(\LCMo
^i_I(R))$$ and
$$\LCMo ^i_I(R)=\LCMo ^i_{I_0}(R_0)\otimes _{R_0}R$$
and
$$\eqalign {D(\LCMo ^i_I(R)))&=\Hom _R(\LCMo ^i_{I_0}(R_0)\otimes
_{R_0}R,\InjH _R(k))\cr &=\Hom _{R_0}(\LCMo ^i_{I_0}(R_0),\InjH
_R(k))\cr &\buildrel (2)\over =\Hom _{R_0}(\LCMo ^i_{I_0}(R_0),\Hom
_{R_0}(R,\InjH _{R_0}(k)))\cr &=\Hom _{R_0}(R,D_{R_0}(\LCMo
^i_{I_0}(R_0)))\cr }$$ For (2) we use the fact
$$\InjH _R(k)=\Hom _{R_0}(R,\InjH _{R_0}(k))$$
We get
$$\eqalign {\LCMo ^i_I(D(\LCMo ^i_I(R)))&=\LCMo ^i_{I_0}(R_0)\otimes _{R_0}\Hom _{R_0}(R,D_{R_0}(\LCMo
^i_{I_0}(R_0)))\cr &\buildrel (3)\over =\Hom _{R_0}(R,\LCMo
^i_{I_0}(R_0)\otimes _{R_0}D_{R_0}(\LCMo ^i_{I_0}(R_0)))\cr &=\Hom
_{R_0}(R,\InjH _{R_0}(k))\cr &\buildrel (2)\over =\InjH _R(k)\cr }$$
For (3) we use the fact that $R$ is a finite-rank free $R_0$-module.
\bigskip
\bigskip
{\bf 7.4.2 Theorem}
\bigskip
Let $R$ be a noetherian local complete regular ring of
equicharacteristic zero, $I\subseteq R$ an ideal of height $h\geq
1$, $x_1,\dots ,x_h\in I$ an $R$-regular sequence and assume that
$$\LCMo ^l_I(R)=0\hbox { for every }l>h\ \ .$$
Then $\LCMo ^h_I(D(\LCMo ^h_I(R)))$ is either $\InjH _R(k)$ or zero.
\par
Proof:
\par
We set
$$D:=D(\LCMo ^h_{(x_1,\dots ,x_h)R}(R))\ \ $$
By theorem 1.1.2, we know that $x_1,\dots ,x_h$ is a $D$-regular
sequence and, therefore, we have
$$\LCMo ^0_{(x_1,\dots ,x_h)R}(D)=\dots =\LCMo ^{h-1}_{(x_1,\dots
,x_h)R}(D)=0\ \ .$$ Because of this, an easy spectral sequence
argument (applied to the composed functor $\Gamma _I\circ \Gamma
_{(x_1,\dots ,x_h)R}$ and to the $R$-module $D$) shows that
$$\LCMo ^h_I(D)=\Gamma _I(\LCMo ^h_{(x_1,\dots ,x_h)R}(D))\subseteq
\LCMo ^h_{(x_1,\dots ,x_h)R}(D)=\InjH _R(k)\ \ .$$ The last equality
is theorem 7.4.1. But, from subsection 7.2 and from [Ly1, Example
2.1 (iv)], it is clear that $\LCMo ^h_I(D)$ has a $D$-module
structure and so, from [Ly1, Theorem 2.4 (b)], we deduce that $\LCMo
^h_I(D)$ is either $\InjH _R(k)$ or zero. Furthermore, the natural
injection
$$\LCMo ^h_I(R)\subseteq \LCMo ^h_{(x_1,\dots ,x_h)R}(R)$$
induces a surjection
$$D\to D(\LCMo ^h_I(R))$$
and hence, as $\LCMo ^h_I$ is right-exact, a surjection
$$\LCMo ^h_I(D)\to \LCMo ^h_I(D(\LCMo ^h_I(R)))\ \ .$$
But again, the last module has a $D$-module structure, and thus,
from [Ly1, Theorem 2.4 (b)] and from what we know already, we
conclude the statement. \vfil \eject {\bf 8 Attached prime ideals
and local homology}
\bigskip
\bigskip
{\bf 8.1 Attached prime ideals -- basics}
\bigskip
This subsection is a collection of definitions and facts about
primary and secondary representation, both in general situations (i.
e. we do not always assume that our modules have any finiteness
properties). We will make use of these facts in subsection 8.2. [BS]
is a reference for the notion of attached primes (of local
cohomology modules).
\bigskip
{\bf 8.1.1 Definition and remark}
\par Let $R$ be a ring, $M\neq 0$ an $R$-module and $N$ an
$R$-submodule of $M$. $M$ is coprimary iff the following condition
holds: For every $x\in R$ the endomorphism $M\buildrel x\over \to M$
given by multiplication by $x$ is injective or nilpotent (i. e.
$\exists N\in \Naturalsign : x^N\cdot M=0$, note that for finitely
generated $M$ this is equivalent to $\forall _{m\in M}\exists N\in
\Naturalsign : x^N\cdot m=0$). If $M$ is coprimary $\sqrt { \Ann
_R(M)}$ is a prime ideal of $R$. In general we say $N$ is a primary
submodule of $M$ iff $M/N$ is coprimary. Now let $U_1,\dots
,U_s\subseteq M$ be submodules of $M$. We say the $s$-tuple $(U_1,
\dots ,U_s)$ is a primary decomposition of (the zero ideal of) $M$
iff the following two conditions hold:
\par
(i) $U_1\cap \dots \cap U_s=0$.
\par
(ii) All $U_i$ are primary submodules of $M$.
\smallskip
In this case $(U_1, \dots ,U_s)$ is called a minimal primary
decomposition of $M$ iff, in addition, the following two statements
hold:
\par
(iii) Every $U_1\cap \dots \cap \hat {U_i}\cap \dots \cap U_s$ is
not zero.
\par
(iv) The ideals $\sqrt {\Ann _R(M/U_i)}$ (for $i=1,\dots ,s$) are
pairwise different.
\smallskip
It is clear that if there exists a primary decomposition of $M$
there is also a minimal one.
\bigskip
{\bf 8.1.2 Definition and remark} \par Let $R$ be a noetherian ring,
$M$ an $R$-module and assume there exists a minimal primary
decomposition $(U_1,\dots ,U_s)$ of $M$. Then the set
$$\{ \sqrt {\Ann _R(M/U_i)}\vert i=1,\dots ,s\} =: \Ass _R(M)$$
does not depend on the choice of a minimal primary decomposition of
$M$ (the proof of this goes just like the well-known proof in case
$M$ is finite). We say the prime ideals of $\Ass _R(M)$ are
associated to $M$
\bigskip
{\bf 8.1.3 Remark} \par Let $R$ be a noetherian ring and $M$ a
noetherian (i. e. finitely generated) $R$-module. Then it is
well-known that $M$ has a (minimal) primary decomposition. Note that
this holds without the hypothesis $R$ is noetherian, but anyway $M$
being noetherian implies that $R/\Ann _R(M)=:\overline R$ is
noetherian and $M$ is a $\overline R$-module.
\bigskip
{\bf 8.1.4 Definition} \par Let $R$ be a noetherian ring and $M$ an
$R$-module. One defines
$$\Ass _R(M):=\{ \goth p\subseteq R\hbox { prime ideal }\vert
\exists m\in M:\goth p=\Ann _R(m)\} .$$ It is easy to see that this
definition agrees with the above one whenever $M$ has a primary
decomposition.
\bigskip
{\bf 8.1.5 Definition and remark} \par Let $R$ be a ring and $M\neq
0$ an $R$-module. By definition, $M$ is secondary iff for every
$x\in R$ the endomorphism $M\buildrel x\over \to M$ given by
multiplication by $x$ is either surjective or nilpotent. Now let $M$
be arbitrary and $U_1,\dots, U_s\subseteq M$ $R$-submodules. We say
the $s$-tuple $(U_1,\dots ,U_s)$ is a secondary decomposition of $M$
iff the following two conditions hold: $U_1+\dots +U_s=M$ and all
$U_i$ are secondary. In this case the secondary decomposition
$(U_1,\dots ,U_s)$ is called minimal iff the following two
conditions hold: All $U_1+\dots +\hat {U_i}+\dots +U_s$ are proper
subsets of $M$ and all $\sqrt {\Ann _R(U_i)}$ are pairwise
different. Again, existence of a secondary decomposition implies
existence of a minimal one.
\bigskip
{\bf 8.1.6 Definition and remark} \par Let $R$ be a noetherian ring
and $M$ an $R$-module; assume there exists a minimal secondary
decomposition $(U_1,\dots ,U_s)$ of $M$. Then the set
$$\Att _R(M):=\{ \sqrt {\Ann _R(U_i)}\vert i=1,\dots s\} $$
does not depend on the choice of a minimal secondary decomposition
of $M$. We say the prime ideals in $\Att _R(M)$ are attached to $M$.
\bigskip
{\bf 8.1.7 Remark} \par Let $R$ be a noetherian ring and $M$ an
artinian $R$-module. Then there exists a (minimal) secondary
decomposition of $M$. The proof is simply a dual version of the
proof of 8.1.3 (which is, of course, well-known). Again this works
also if $R$ is not noetherian.
\bigskip
{\bf 8.1.8 Definition} \par Let $R$ be a noetherian ring and $M$ an
$R$-module. We define
$$\Att _R(M):= \{ \goth p\subseteq R \hbox { prime ideal }\vert
\exists \hbox { an }R\hbox {-submodule }U \subseteq M: \goth p=\Ann
_R(M/U)\} .$$ Is is not very difficult to see that this definition
agrees with the first one if $M$ has a secondary decomposition.
\bigskip
{\bf 8.1.9 Remark} \par Let $(R,\goth m)$ be a noetherian local
ring, $M$ an $R$-module and $(U_1,\dots ,U_s)$ a minimal primary
decomposition of $M$. The following implications are clear by
duality:
\smallskip
(i) $U_1\cap \dots \cap U_s=0 \Rightarrow D(M/U_1)+\dots +\dots
D(M/U_s)=D(M)$
\par
(ii) $M/U_i$ is coprimary $\Rightarrow D(M/U_i)$ is secondary (for
every $i$)
\par
(iii) The primary decomposition $(U_1,\dots ,U_s)$ of $M$ is minimal
$\Rightarrow $\par the secondary decomposition $(D(M/U_1),\dots
,D(M/U_s))$ of $D(M)$ is minimal.
\par
(iv) $\Ann _R(M/U_i)=\Ann _R(D(M/U_i))$ (for every $i$)
\smallskip
Thus we have
$$\Ass _R(M)=\Att _R(D(M))\ \ .$$
In a very similar way the following statement holds: Any (minimal)
secondary decomposition of $M$ induces a (minimal) primary
decomposition of $D(M)$. In particular, if $M$ has a secondary
decomposition:
$$ \Att _R(M)=\Ass _R(D(M))\ \ .$$
{\bf 8.1.10 Remark} \par It is true that if $U_1,\dots ,U_s$ are
arbitrary submodules of $R$ such that $(D(M/U_1),\dots ,D(M/Us))$ is
a (minimal) secondary decomposition of $D(M)$ then $(U_1,\dots
,U_s)$ is a (minimal) primary decomposition of $M$, but note that we
do not know that every submodule of $D(M)$ is of the form $D(M/U)$
for some submodule $U$ of $M$. Similarly, if $U_1,\dots ,U_s$ are
arbitrary submodules of $M$ such that $(D(M/U_1),\dots ,D(M/U_s))$
is a (minimal) primary decomposition of $D(M)$ then $(U_1,\dots
,U_s)$ is a (minimal) secondary decomposition of $M$.
\bigskip
{\bf 8.1.11 Remark} \par Let $(R,\goth m)$ be a noetherian local
ring, $\goth p$ a prime ideal of $R$ and $M$ an $R$-module. Then
$$\eqalign {\goth p\in \Ass _R(M)&\iff \exists \hbox { finitely generated
submodule }U \hbox { of } M: \goth p=\Ann _R(U),\cr \goth p\in \Att
_R(D(M))&\iff \exists \hbox { submodule }U^\prime \hbox { of }D(M):
\goth p=\Ann _R(D(M)/U^\prime ).\cr }$$ In particular the existence
of a submodule $U$ of $M$ satisfying $\goth p=\Ann _R(U)$ implies
$\goth p\in \Att _R(D(M))$. Therefore we have
$$\Ass _R(M)\subseteq \Att _R(D(M)).$$
This inclusion is strict in general: Take for example $M=\InjH
=\InjH _R(R/\goth m)$, an $R$-injective hull of $R/\goth m$: $\Ass
_R(\InjH _R(R/\goth m))=\{ \goth m\} $, but $D(\InjH )=\hat R$ and
so $\Att _R(D(\InjH ))=\Spec (R)$. But nevertheless a stronger
inclusion holds (plug in $D(M)$ for $M$ in theorem 8.1.12 to see
that it is actually stronger):
\bigskip
{\bf 8.1.12 Theorem}
\par
Let $(R,\goth m)$ be a noetherian local ring and $M$ an $R$-module.
Then
$$\Ass _R(D(M))\subseteq \Att _R(M)$$
and the sets of prime ideals maximal in each side respectively
coincide:
\par
$$\{ \goth p\vert \goth p\hbox { maximal in } \Ass _R(D(M))\}=\{ \goth p\vert \goth p\hbox { maximal in } \Att
_R(M)\}.$$ Proof: \par Let $\goth p\in \Ass _R(D(M))$ be arbitrary.
There exists a submodule $U^\prime $ of $D(M)$ such that $U^\prime
=R\cdot u^\prime \cong R/\goth p$ for some $u^\prime \in U^\prime
\subseteq D(M)$. $u^\prime $ induces a monomorphism $\overline
{u^\prime }:M/\ker (u^\prime )\to \InjH $ and so we have
$$\goth p=\Ann _R(U^\prime )=\Ann _R(u^\prime )=\Ann _R(\overline
{u^\prime })=\Ann _R(M/\ker (u^\prime ));$$ this implies $\goth p\in
\Att _R(M)$. Having proved this we only have to show that an
arbitrary prime ideal $\goth p$ of $R$ which is maximal in $\Att
_R(M)$ is associated to $D(M)$: $\goth p\in \Att _R(M)$ implies
$M/\goth pM\neq 0$ and so we must have $\Hom _R(R/\goth
p,D(M))=D(M/\goth pM)\neq 0$; but by the maximality hypothesis on
$\goth p$ implies $\goth p\in \Ass _R(D(M))$.
\bigskip
{\bf 8.1.13 Theorem}
\par
Let $(R,\goth m)$ be a noetherian local ring and $M$ an $R$-module.
Assume $(\goth p_i)_{i\in \Naturalsign }$ is a sequence of prime
ideals attached to $M$; assume furthermore that $\goth q:=\bigcap
_{i\in \Naturalsign }\goth p_i$ is a prime ideal of $R$. Then $\goth
q$ is also attached to $M$.
\par
Proof: \par For every $i$ we choose a quotient $M_i$ of $M$ such
that $\Ann _R(M_i)=\goth q_i$. Now the canonically induced map
$\iota :M\to \prod _{i\in \Naturalsign }M_i$ induces a surjection
$M\to \im (\iota )$; we obviously have $\bigcap _{i\in \Naturalsign
}\goth p_i\subseteq \Ann _R(\im (\iota ))$; on the other hand, for
every $i$ and every $s\in R\setminus \goth p_i$ there is a
$\overline {m_i}\in M_i$ coming from an element $m_i\in M$ that has
$s\cdot \overline {m_i}\neq 0$. But this implies that $s$ cannot
annihilate $\im (\iota )$; therefore
$$\Ann _R(\im (\iota))=\bigcap
_{i\in \Naturalsign }\goth p_i=\goth q$$
and the statement follows.
\bigskip
\bigskip
{\bf 8.2 Attached prime ideals -- results}
\bigskip
This subsection contains results on attached prime ideals (of local
cohomology modules). Our technique bases on subsection 8.1 where
some relations between attached primes of a module and associated
primes of the Matlis dual of the same module were established. This
method does not only lead to an easy proof of a known result
(theorem 8.2.1, see also remark 8.2.2), but also enables us to find
more attached prime ideals (of a local cohomology module, see
theorem 8.2.3 and corollary 8.2.4 for details). Furthermore, the
study of attached prime ideals leads to new evidence for conjecture
(*) (this evidence comes, essentially, from theorem 8.1.13 which
describes a property of the set of attached prime ideals that is
necessary for being closed under generalization).
\bigskip
There are some results on the set of attached primes of local
cohomology modules: In [MS, theorem 2.2] it was shown that if
$(R,\goth m)$ is a noetherian local ring and $M$ is a finitely
generated $R$-module then
$$\Att _R(\LCMo ^{\dim (M)}_\goth m(M))=\{ \goth p\in \Ass _R(M)\vert
\dim (R/\goth p)=\dim (M)\} $$ holds. In [DY, Theorem A] this was
generalized to
$$\Att _R(\LCMo ^{\dim (M)}_\goth a(M))=\{ \goth p\in \Ass
_R(M)\vert \cd (\goth a,R/\goth p)=\dim (M)\} ,$$ where $\goth
a\subseteq R$ is an ideal and $\cd (\goth a,R/\goth p):=\max \{ l\in
\Naturalsign \vert \LCMo ^l_\goth a(R/\goth p)\neq 0\} $. We are
going to show (theorem 8.2.1) that the results of section 8.1 lead
to a natural proof of this result and, furthermore, to new results
on the attached primes of local cohomology modules (8.2.3 -- 8.2.6).
\bigskip
Let $(R,\goth m)$ be a noetherian local $n$-dimensional ring and
$\goth a\subseteq R$ an ideal. Then $\LCMo ^n_\goth a(R)$ is an
artinian $R$-module and hence
$$\Ass _R(D(\LCMo ^n_\goth a(R)))=\Att _R(\LCMo ^n_\goth a(R)).$$
Now assume that we have ($\LCMo ^n_\goth a(R)\neq 0$ and) $\goth
p\in \Att _R(\LCMo ^n_\goth a(R))$; then we get
$$0\neq \LCMo ^n_\goth a(R)/\goth p\LCMo ^n_\goth a(R)=\LCMo
^n_\goth a(R/\goth p),$$ i. e. $\goth p\in \Assh (R) $($:=\{ \goth
q\in \Spec(R)\vert \dim (R/\goth q)=\dim (R)\}$) and $\cd (\goth
a,R/\goth p)=n$.
\par
Now suppose conversely that we have a prime ideal $\goth p$ of $R$
such that $\cd (\goth a,R/\goth p)=n$, equivalently $\LCMo ^n_\goth
a(\hat {R/\goth p})\neq 0$. By Hartshorne-Lichtenbaum vanishing we
get a prime ideal $\goth q\subseteq \hat R$ satisfying $\goth
p=\goth q\cap R$ and $\sqrt {\goth a\hat R+\goth q}=\goth m_{\hat
R}$($:=$maximal ideal of $\hat R$); this in turn implies
$$0\neq \LCMo ^n_{\goth a\hat R}(\hat R/\goth q)=\LCMo ^n_{\goth
m_{\hat R}}(\hat R/\goth q).$$ Matlis duality theory shows that
$\goth q\in \Ass _{\hat R}(D(\LCMo ^n_{\goth a\hat R}(\hat R)))$. It
is easy to see that
$$D(\LCMo ^n_{\goth a\hat R}(\hat R))=D(\LCMo ^n_\goth a(R)),$$
holds canonically, the $D$-functors taken over $\hat R$ resp. over
$R$. Thus we have shown
$$\Att _R(\LCMo ^n_\goth a(R))=\{ \goth p\vert \cd (\goth a,R/\goth
p)=n\} .$$ For every finitely generated $R$-module $M$ we can apply
this result to the ring $R/\Ann _R(M)$ and we get
\bigskip
{\bf 8.2.1 Theorem}
\par
Let $(R,\goth m)$ be a noetherian local ring, $\goth a\subseteq R$
ein Ideal and $M$ a finitely generated $n$-dimensional $R$-module.
Then
$$\Att _R(\LCMo ^n_\goth a(M))=\{ \goth p\in \Ass _R(M)\vert \cd
(\goth a,R/\goth p)=n\} $$ holds. \bigskip {\bf 8.2.2 Remark} \par
This is [DY, Theorem A], where it was proved by different means.
\bigskip
In subsection 8.1 we established several relations between attached
primes of a module and associated primes of the Matlis dual of the
same module; theorem 8.2.3 is a consequence of these relations; we
can retrieve more information from these relations to get new
theorems on the attached primes of top local cohomology modules
(remarks 8.2.5):
\bigskip
{\bf 8.2.3 Theorem} \par Let $(R,\goth m)$ be a $d$-dimensional
noetherian local ring.
\par
(i) If $J$ is an ideal of $R$ such that $\dim (R/J)=1$ and $\LCMo
^d_J(R)=0$ then
$$\Assh (R)\subseteq \Att _R(\LCMo ^{d-1}_J(R))$$
holds. If, in addition, $R$ is complete, one has
$$\Att _R(\LCMo ^{d-1}_J(R))=\{ \goth p\in \Spec (R)\vert \dim
(R/\goth p)=d-1, \sqrt {\goth p+J}=\goth m\} \cup \Assh (R).$$ \par
(ii) For any $x_1,\dots ,x_i\in R$ there is an inclusion
$$\{ \goth p\in \Spec (R)\vert x_1,\dots ,x_i \hbox { is a part of a
system of parameters of }R/\goth p\} \subseteq \Att _R(\LCMo
^i_{(x_1,\dots ,x_i)R}(R)).$$ Proof: \par (i) Note that theorems
3.2.6 and 3.2.7 show that one has $\Assh (R)=\Assh (D(\LCMo
^{d-1}_J(R)))$ in the given situation and, if $R$ is complete,
$$\Ass
_R(D(\LCMo ^{d-1}_J(R)))=\{ \goth p\in \Spec (R)\vert \dim (R/\goth
p)=1, \dim (R/(\goth p+J))=0\} \cup \Assh (R)\ \ .$$
Now we use
theorem 8.1.12 and remark: If $R$ is complete, given an arbitrary
$\goth p\in \Att _R(\LCMo ^{d-1}_J(R))$ it follows that $\LCMo
^{d-1}_J(R/\goth p)\neq 0$ and hence, by Hartshorne Lichtenbaum
vanishing, that $\dim (R/\goth p)\geq d-1$ and, if $\dim (R/\goth
p)=d-1$, that $\goth p+J$ is $\goth m$-primary.
\par
(ii) Follows from theorems 8.1.12 and 3.1.3 (ii).
\bigskip
{\bf 8.2.4 Corollary}
\par
Let $(R,\goth m)$ be a noetherian local ring. For every $x\in R$ one
has
$$\Att _R(\LCMo ^1_{xR}(R))=\Spec (R)\setminus {\frak V}(x).$$
Proof: \par "$\subseteq $" Let $\goth p\in \Att _R(\LCMo
^1_{xR}(R))$. Then
$$0\neq \LCMo ^1_{xR}(R)/\goth p\LCMo
^1_{xR}(R)=\LCMo ^1_{xR}(R/\goth p) \Rightarrow x\not\in \goth p\ \
.$$ "$\supseteq $" follows e. g. from 3.1.3 (ii).
\bigskip
{\bf 8.2.5 Remarks}
\par
(i) It was shown in remark 1.2.1 that, for any $x_1,\dots,x_i\in R$,
there is an inclusion
$$\Ass _R(D(\LCMo ^i_{(x_1,\dots ,x_i)R}(R)))\subseteq \{ \goth p\in
\Spec (R)\vert \LCMo ^i_{(x_1,\dots ,x_i)R}(R/\goth p)\neq 0\}.$$ By
what we have proved so far it is clear that there is a chain of
inclusions
$$\Ass _R(D(\LCMo ^i_{(x_1,\dots ,x_i)R}(R)))\subseteq \Att _R(\LCMo
^i_{(x_1,\dots ,x_i)R}(R))\subseteq \{ \goth p\in \Spec (R)\vert
\LCMo ^i_{(x_1,\dots ,x_i)R}(R/\goth p)\neq 0\} .$$ (ii) Conjecture
(*) says that, for any sequence $x_1,\dots ,x_i$ in $R$, the
inclusion $$\Ass _R(D(\LCMo ^i_{(x_1,\dots ,x_i)R}(R)))\subseteq \{
\goth p\in \Spec (R)\vert \LCMo ^i_{(x_1,\dots ,x_i)R}(R/\goth
p)\neq 0\}$$ is an equality; if this could be shown to be true, we
could conclude that
$$\Ass _R(D(\LCMo ^i_{(x_1,\dots ,x_i)R}))=\Att _R(\LCMo
^i_{(x_1,\dots ,x_i)R}(R))\ \ .$$ (iii) In the situation of theorem
8.2.3 (i) the attached primes of the top local cohomology module
coincide with the associated primes of the Matlis dual of the top
local cohomology module.
\bigskip
{\bf 8.2.6 Remarks} \par We now assume that $k$ is a field and
$R=k[[X_1,\dots ,X_n]]$ is a power series algebra over $k$ in $n$
variables $X_1,\dots ,X_n$; let $i\in \{ 1,\dots ,n\} $. Theorems
4.2.1, 4.3.4 and 8.1.12 imply the following statements:
\smallskip
(i) In the case $i=n$ we have
$$\Att _R(\LCMo ^n_{(X_1,\dots
,X_n)R}(R))=\{ 0\} \ \ .$$
\par
(ii) If $i=n-1$,
$$\Att _R(\LCMo ^{n-1}_{(X_1,\dots ,X_{n-1})}(R)=\{
0\} \cup \{ pR\vert p\in R\hbox{ prime element, }p\not\in (X_1,\dots
,X_{n-1})R\} $$ holds.
\par
(iii) Finally we concentrate on the case $i=n-2$, where we have the
following statements:
\par ($\alpha $) $\{ 0\} \in \Att _R(\LCMo ^{n-2}_{(X_1,\dots
,X_{n-2})R}(R))$;
\par ($\beta $) If $\goth p$ is a height-two prime ideal of $R$ such that
$\sqrt {(X_1,\dots ,X_{n-2})R+\goth p}=\goth m$ then $$\goth p\in
\Att _R(\LCMo ^{n-2}_{(X_1,\dots ,X_{n-2})R}(R))$$
holds.
\par
($\gamma $) Conversely, $\goth p\in \Att _R(\LCMo ^{n-2}_{(X_1,\dots
,X_{n-2})R}(R))$ implies that $\height (\goth p)\leq 2$; \par
($\delta $) If $p\in R$ is a prime element such that $p\not\in
(X_1,\dots ,X_{n-2})R$ then
$$pR\in \Att _R(\LCMo ^{n-2}_{(X_1,\dots
,X_{n-2})R}(R))$$ holds.\par ($\epsilon $) If $p\in R$ is a minimal
generator of $(X_1,\dots ,X_{n-2})R$ then
$$pR\not\in \Att _R(\LCMo ^{n-2}_{(X_1,\dots ,X_{n-2})R}(R))$$
holds.
\par
($\zeta $) Because of theorems 4.3.4 and 8.1.12, for every prime
element $p\in (X_1,\dots ,X_{n-2})R\cap (X_{n-1},X_n)R$ there exist
infinitely many (pairwise different) prime ideals $(\goth p_l)_{l\in
\Naturalsign }$ of height two attached to $\LCMo ^{n-2}_{(X_1,\dots
,X_{n-2})R}(R)$ and containing $p$. As any $q\in \bigcap _{l\in
\Naturalsign }\goth p_l$ must satisfy $\height (p,q)R<2$ it is clear
that we have $pR=\bigcap _{l\in \Naturalsign }\goth p_l$. Now
theorem 8.1.13 implies $pR\in \Att _R(\LCMo ^{n-2}_{(X_1,\dots
,X_{n-2})R}(R))$. But in view of theorems 1.2.3 and 8.1.12 it is
clear that $pR\in \Att _R(\LCMo ^{n-2}_{(X_1,\dots ,X_{n-2})R}(R)$
is a necessary condition for conjecture (*). This gives new evidence
for conjecture (*).
\bigskip
\bigskip
{\bf 8.3 Local homology and a necessary condition for
Cohen-Macaulayness}
\bigskip
Let $(R,\goth m)$ be a noetherian, local ring, $M$ an $R$-module and
$I$ an ideal of $R$. It is well-known that $\LCMo ^{\dim (M)}_I(M)$
is artinian for any proper ideal $I$ of $R$ provided $M$ is finitely
generated as $R$-module (cp. [Me]).
\par
There is a theory of local homology modules (cp. [T1] and [T2]): If
$X$ is an artinian $R$-module and $\underline x=x_1,\dots ,x_r$ is a
sequence of elements in $\goth m$, the $i$-th local homology module
$\LCMo ^{\underline x}_i(X)$ of $X$ with respect to $\underline x$
is defined by
$$\vtop{\baselineskip=1pt \lineskiplimit=0pt \lineskip=1pt\hbox{lim}
\hbox{$\longleftarrow $} \hbox{$^{^{n\in \bf N}}$}} H_i(K_\bullet
(x_1^n,\dots ,x_r^n;X))\ \ ,$$ where $K_\bullet (x_1^n,\dots
,x_r^n;X)$ is the Koszul complex of $X$ with respect to $x_1^n,\dots
,x_r^n$ and $H_i$ means taking the homology of this complex at the
$i$-th position; then $\LCMo ^{\underline x}_i(\ )$ is an
$R$-linear, covariant functor from the category of artinian
$R$-modules to the category of $R$-modules.
\par
We repeat the notions of Noetherian dimension $\Ndim (X)$ and width
of $X$, $\width (X)$: For $X=0$ one puts $\Ndim (X)=-1$, for $X\neq
0$ $\Ndim (X)$ denotes the least integer $r$ such that
$0:_X(x_1,\dots ,x_r)R$ has finite length for some $x_1,\dots
,x_r\in \goth m$. Now let $x_1,\dots ,x_n\in \goth m$. $x_1,\dots
,x_n$ is an $X$-coregular sequence if
$$0:_X(x_1,\dots ,x_{i-1})R\buildrel x_i\over \to 0:_X(x_1,\dots ,x_{i-1})R$$
is surjective for $i=1,\dots ,n$. $\width (X)$ is defined as the
length of a (in fact any) maximal $X$-coregular sequence in $\goth
m$. Details on $\Ndim (X)$ and $\width (X)$ can be found in [Oo] and
[Ro], here we cite one general fact: For any artinian $R$-module $X$
$$\width (X)\leq \Ndim (X)<\infty $$
holds and $X$ is co-Cohen-Macaulay if and only if $\width (X)= \Ndim
(X)$ holds (by definition). Tang has shown ([T1, Proposition 2.6])
that $\LCMo ^{\dim (M)}_\goth m(M)$ is co-Cohen-Macaulay (of
Noetherian dimension $\dim (M)$) if $M$ is a finitely generated
Cohen-Macaulay $R$-module and ([T1, Theorem 3.1]) that
$$\LCMo ^{x_1,\dots ,x_d}_{\dim (M)}(\LCMo ^{\dim
(M)}_\goth m(M))=\hat M$$ holds (here $x_1,\dots ,x_d$ is a s. o. p.
of $M$ and we still assume that $M$ is Cohen-Macaulay). Tang asks
([T1, Remark 3.5]) if one can show that $\LCMo ^{\underline x}_d(X)$
is finitely generated if $X$ is an artinian $R$-module of
N.dimension $d$ and $\underline x=x_1,\dots ,x_d$ is such that
$0:_X\underline x$ has finite length.
\par
In the example 8.3.1 below we give a negative answer to this
question. However, under the additional assumption that $R$ is
complete, we show that $\LCMo ^{\underline x}_d(X)$ is a finitely
generated $R$-module (theorem 8.3.3) and draw some consequences
establishing various duality results (theorem 8.3.5). As an
application we present a necessary condition for a given finite
$R$-module $M$ to be Cohen-Macaulay (corollary 8.3.6).
\par
$\underline x$ will always stand for a sequence $x_1,\dots ,x_d$ in
$\goth m$. The results of this and the next subsection can also be
found in [H5].
\bigskip
\parindent=0pt
{\bf 8.3.1 Example} \par Let $k$ be a field, $T$ a variable and $R$
the noetherian, local ring $k[T]_{(T)}$. Set $X:=T^{-1}\cdot
k[T^{-1}]:=\{ a_{-1}T^{-1}+\dots +a_{-n}T^{-n}\vert n\in
\Naturalsign ^+, a_{-1},\dots ,a_{-n}\in k\} $. $X$ has a $\hat
R=k[[T]]$-structure (such that $T^m\cdot T^{-n}=T^{m-n}$ if $m-n\leq
-1$ and $=0$ if $m-n\geq 0$, where $m\geq 0,n\geq -1$) and thus it
also has an $R$-module-structure. Every non-trivial $R$-submodule of
$X$ has the form $<T^{-n}>_X$ for some $n\geq 1$ and therefore $X$
is an artinian $R$-module. Furthermore $(0:_XT)=k\cdot T^{-1}$ is of
finite length (and so $\Ndim(X)=1$) and $\LCMo ^T_1(X)$ is the
indirect limit over all $(0:_XT^l)$, where the transition maps
$(0:_XT^{l+1})\to (0:_XT^l)$ are induced by multiplication by $T$.
An elementary calculation shows $\LCMo ^T_1(X)=k[[T]]=\hat R$ which
is not finite as an $R$-module.
\bigskip
But more can be said:
\bigskip
{\bf 8.3.2 Remark} \par Let $(R,\goth m)$ be a local noetherian
regular $d$-dimensional ring, $X$ an artinian co-Cohen-Macaulay
$R$-module, $\Ndim (X)=d$, $\underline x=x_1,\dots ,x_d\in \goth m$
such that $(0:_X\underline x)$ is of finite length. Then
$$\LCMo ^{\underline x}_d(X) \hbox { is a finite }R\hbox{-module}\iff R\hbox {
is complete}$$ holds. \par Proof: \par
$\Leftarrow$ follows from
theorem 8.3.3 below. $\Rightarrow $: From [T1, Remark 3.5] it
follows that $\depth (\LCMo ^{\underline x}_d(X))=d$ both as an $R-$
and as an $\hat R$-module; but now the Auslander-Buchsbaum formula
implies that $\LCMo ^{\underline x}_d(X)$ is a finite free $\hat
R$-module and so we must have $\hat R=R$.
\bigskip
From now on we assume that $R$ is complete and show at first that
the top local homology module is always finite; this is done,
essentially, by Matlis duality.
\bigskip
{\bf 8.3.3 Theorem}
\par
Let $(R,\goth m)$ be a noetherian, local, complete ring, $X$ an
artinian $R$-module of N.dimension $d$; let $x_1,\dots ,x_d\in \goth
m$ be such that $0:_X(x_1,\dots ,x_d)R$ has finite length. Then
$\LCMo ^{\underline x}_d(X)$ is a finitely generated $R$-module.
\par
Proof:
\par
$x_1,\dots ,x_d$ form a system of parameters for $D(X)$, because
$$D(X)/(x_1,\dots ,x_d)D(X)=D(0:_X(x_1,\dots ,x_d)R)$$
has finite
length and $\dim (D(X))=\Ndim (X)=d$. Using Matlis-duality we have
$$\eqalign {\LCMo ^{\underline x}_d(X)&=\LCMo ^{\underline x}_d(D(D(X)))\cr &=\vtop{\baselineskip=1pt \lineskiplimit=0pt \lineskip=1pt\hbox{lim}
\hbox{$\longleftarrow $} \hbox{$^{^{n\in \bf N}}$}} H_d(K_\bullet
(x_1^n,\dots ,x_d^n;D(D(X))))\cr &=\vtop{\baselineskip=1pt
\lineskiplimit=0pt \lineskip=1pt\hbox{lim} \hbox{$\longleftarrow $}
\hbox{$^{^{n\in \bf N}}$}}D(H^d(K^\bullet (x_1^n,\dots
,x_d^n;D(X))))\cr &=D(\vtop{\baselineskip=1pt \lineskiplimit=0pt
\lineskip=1pt\hbox{lim} \hbox{$\longrightarrow $} \hbox{$^{^{n\in
\bf N}}$}}H^d(K^\bullet (x_1^n,\dots ,x_d^n;D(X))))\cr &=D(\LCMo
^d_{(x_1,\dots ,x_d)R}(D(X)))\ \ ,\cr }$$ and the last module is
finitely generated because $\LCMo ^d_{(x_1,\dots ,x_d)R}(D(X))$ is
artinian.
\bigskip
{\bf 8.3.4 Corollary}
\par
Let $(R,\goth m)$ be a noetherian, local, complete ring and $X$ a
co-Cohen-Macaulay $R$-module of N.dimension $d$; let $x_1,\dots
,x_d\in \goth m$ be such that $0:_X(x_1,\dots ,x_d)R$ has finite
length. Then $\LCMo ^{x_1,\dots ,x_d}_d(X)$ is a Cohen-Macaulay
module. In particular if $d=\dim (R)$, $\LCMo ^{x_1,\dots
,x_d}_d(X)$ is a maximal Cohen-Macaulay module.
\par
Proof:
\par
The statements follow from theorem 8.3.3 and [T1, Remark 3.5].
\bigskip
Let $(R,\goth m)$ be a noetherian, local, complete ring. Let $\cal
N$ (resp. $\cal A$) denote the set of isomorphism classes of
noetherian (resp. of artinian) $R$-modules. We have maps $F_1$ and
$F_2$ from $\cal N$ to $\cal A$ induced by
$$M\buildrel F_1\over \mapsto \hbox{Matlis dual of }M$$
and
$$M\buildrel F_2\over \mapsto \LCMo ^{\dim (M)}_\goth m(M)$$
For $F_2$ it does not make any difference if we take $\LCMo ^{\dim
(M)}_{(x_1,\dots ,x_{\dim (M)})R}(M)$ instead of $\LCMo ^{\dim
(M)}_\goth m(M)$ (for any system of parameters $x_1,\dots ,x_{\dim
(M)}$ of $M$). Similarly we have maps $G_1$ and $G_2$ from $\cal A$
to $\cal N$ induced by
$$X\buildrel G_1\over \mapsto \hbox{Matlis-dual of }X$$
and
$$X\buildrel G_2\over \mapsto \LCMo ^{x_1,\dots ,x_{\Ndim
(X)}}_{\Ndim(X)}(X)$$ (here $x_1,\dots ,x_{\Ndim (X)}$ are such that
$0:_X(x_1,\dots ,x_{\Ndim (X)})R$ has finite length). By
Matlis-duality we have
$$F_1\circ G_1=\id _{\cal A}, G_1\circ F_1=\id _{\cal N}\ \ .$$
From the proof of theorem 8.3.3 one understands that
$$F_1\circ G_2=F_2\circ G_1=:T$$
and hence
$$G_1\circ F_2=G_2\circ F_1=:T^\prime \ \ ,$$
$$G_2=G_1\circ F_2\circ G_1=G_1\circ T, F_2=F_1\circ G_2\circ F_1=F_1\circ
T^\prime \ \ .$$
\par
{\bf 8.3.5 Theorem}
\par
Let $(R,\goth m$) be a noetherian, local, complete ring. Let $M$ be
a noetherian and $X$ an artinian $R$-module. Then
\par
(i) If $M$ is Cohen-Macaulay, then $F_2(M)$ is co-Cohen-Macaulay.
\par
(ii) If $M$ is Cohen-Macaulay, then $F_1(M)$ is co-Cohen-Macaulay.
\par
(iii) If $X$ is co-Cohen-Macaulay, then $G_2(M)$ is Cohen-Macaulay.
\par
(iv) If $X$ is co-Cohen-Macaulay, then $G_1(M)$ is Cohen-Macaulay.
\par
Proof:
\par
(ii) and (iv) are easily proved using Matlis-duality theory. (i) is
proved by [T1, Proposition 2.6])  and now (iii) follows from
$G_2=G_1\circ F_2\circ G_1$.
\bigskip
Let ${\cal N}_0$ (resp. ${\cal A}_0$) denote the set of isomorphism
classes of noetherian Cohen-Macaulay modules (resp. of artinian
co-Cohen-Macaulay modules). Then, by theorem 8.3.5,
$F_1,F_2,G_1,G_2$ induce maps between ${\cal N}_0$ and ${\cal A}_0$
in an obvious way. [T1, theorems 3.1 and 3.4] imply $F_2\circ
G_2=\id _{{\cal A}_0}$ and $G_2\circ F_2=\id _{{\cal N}_0}$. We
deduce $G_1=G_2\circ F_1\circ G_2$, $F_1=F_2\circ G_1\circ
F_2,T^2=\id ,{T^\prime }^2=\id $ on ${\cal N}_0$ and ${\cal A}_0$.
\bigskip
As an application we get a necessary condition for a finite module
to be Cohen-Macaulay:
\par
{\bf 8.3.6 Corollary}
\par
(i) Let $\omega _R$ be a dualizing module for $R$ (it exists
uniquely up to isomorphism since $R$ is complete). Assume that $M$
is Cohen-Macaulay. Then $\Ext ^{\dim (R)-\dim (M)}_R(M,\omega _R)$
is Cohen-Macaulay.
\par
(ii) In particular if there exists an ideal $I$ of $R$ such that
$I\subseteq \Ann _R(M)$, $\dim (R/I)=\dim (M)$ and $R/I$ is
Gorenstein, Cohen-Macaulayness of $M$ implies Cohen-Macaulayness of
$\Hom _{\overline R}(M,\overline R)$ (here $\overline R:=R/I$). Such
an ideal $I$ exists, for example, if $R$ itself is Gorenstein.
\par
Proof:
\par
The statements follow from local duality and theorem 8.3.5.
\bigskip
\bigskip
{\bf 8.4 Local homology and Cohen-Macaulayfications}
\bigskip
In the text following theorem 8.3.5 we have seen $G_2\circ
F_2=\id_{{\cal N}_0}$ and $F_2\circ G_2=\id _{{\cal A}_0}$. Now we
turn our interest to the question: What can be said about $G_2\circ
F_2$ in general, that is, on $\cal N$?
\bigskip
{\bf 8.4.1 Definition} \par Let $(R,\goth m)$ be a noetherian,
local, complete ring and $M$ a noetherian (i. e. finitely generated)
$R$-module. Let $\tilde M$ be a finitely generated  $R$-module
containing $M$ as a submodule. We say $\tilde M$ is a
Cohen-Macaulayfication of $M$ if the following three conditions
hold:
\par
(i) $\tilde M$ is Cohen-Macaulay.
\par
(ii) $\dim (\tilde M)=\dim (M)$.
\par
(iii) $\dim (\tilde M/M)\leq \dim M-2$ (this condition is equivalent
to $\LCMo ^{\dim (M)-1}_\goth m(\tilde M/M)=\LCMo ^{\dim (M)}_\goth
m(\tilde M/M)=0$).
\bigskip
In the sequel we won't always strictly distinguish between a module
$M$ and its isomorphism class, for reasons of simplicity.
\bigskip
{\bf 8.4.2 Theorem}
\par
Let $(R,\goth m)$ be a noetherian, local, complete ring and $M$ a
noetherian $R$-module. If $M$ has a Cohen-Macaulayfication, it has
(up to an $M$-isomorphism) only one Cohen-Macaulayfication, namely
$(G_2\circ F_2)(M)$.
\par
Proof:
\par
Let $\tilde M$ be a Cohen-Macaulayfication of $M$. We consider the
short exact sequence $0\to M\to \tilde M\to \tilde M/M\to 0$ and its
long exact cohomology sequence induced by applying $\Gamma _\goth m$
to it: Because of condition (iii) of definition 8.4.1 we get a
canonical isomorphism
$$\LCMo ^{\dim (M)}_\goth m(M)=\LCMo ^{\dim
(M)}_\goth m(\tilde M)\buildrel \hbox {8.4.1 (ii)}\over =\LCMo
^{\dim (\tilde M)}_\goth m(\tilde M)$$
and therefore $\tilde
M=(G_2\circ F_2)(\tilde M)=(G_2\circ F_2)(M)$.
\bigskip
{\bf 8.4.3 Remark}
\par
Goto (cp. [Go]) has shown: If $(A,\goth m)$ is a noetherian, local,
$d$-dimensional ring with total quotient ring $Q(A)$, the following
conditions are equivalent:
\par
(i) There is a Cohen-Macaulay ring $B$ between $A$ and $Q(A)$ such
that $B$ is finitely generated as an $A$-module, $\dim (B_\goth
n)=d$ for every maximal ideal $\goth n$ of $B$ and $\goth m\cdot
B\subseteq A$.
\par
(ii) $A$ is a Buchsbaum ring (see [SV] for details on Buchsbaum
rings) and $\LCMo ^i_\goth m(A)=0$ for $i\neq 1,d$.
\par
In this case, if $d\geq 2$, $B$ is uniquely determined and Goto
([Go]) calls it  the Cohen-Macaulayfication of $A$.
\bigskip
{\bf 8.4.4 Remark}
\par
Cohen-Macaulayfication in our sense is a generalization of Goto's
concept of Cohen-Macaulayfication:
\bigskip
{\bf 8.4.5 Theorem}
\par
Let $(R,\goth m)$ be a noetherian, local, complete ring, and assume
that the Cohen-Macaulayfication $B$ of $R$ (in the sense of Goto)
exists. Then $B$ is also a Cohen-Macaulayfication in our sense.
\par
Proof:
\par
Because of $\goth m\cdot B\subseteq R$ we have $\goth m\cdot
(B/R)=0$, which implies that $B/R$ is a finite-dimensional $R/\goth
m$-vector space. Because of $d=\dim (R)\geq 2$ we must have $\LCMo
^{d-1}_\goth m(B/R)=\LCMo ^d_\goth m(B/R)=0$.
\bigskip
{\bf 8.4.6 Remark}
\par
In particular if $(R,\goth m)$ is a noetherian, local, complete
Buchsbaum-ring of dimension $d\geq 2$ such that $\LCMo ^i_\goth
m(R)=0$ for $i\neq 1,d$, the $R$-module $R$ has a
Cohen-Macaulayfication.
\bigskip
{\bf 8.4.7 Example}
\par
An easy example is given by $R=k[[x_1,x_2,x_3,x_4]]/(x_1,x_2)\cap
(x_3,x_4)$. In the sense of Goto as well as in our sense $R$ has a
Cohen-Macaulayfication given by $(k[[x_1,\dots
,x_4]]/(x_1,x_2))\oplus (k[[x_1,\dots ,x_4]]/(x_3,x_4))$; this can
be seen either directly or by remarking that $R$ is a 2-dimensional
Buchsbaum ring with $\LCMo ^i_\goth m(R)=0$ for $i\neq 1,2$.\vfil
\eject {\bf References} \normal
\bigskip
\parindent=1.2cm
\def\litem{\par\noindent \hangindent=\parindent\ltextindent}
\def\ltextindent#1{\hbox to \hangindent{#1\hss}\ignorespaces}
\litem{[Ba]} Bass, H. On the ubiquity of Gorenstein rings, {\it
Math. Z.} {\bf 82}, (1963) 8--28.
\medskip
\litem{[BHe]} Brodmann, M. and Hellus, M. Cohomological patterns of
coherent sheaves over projective schemes, {\it Journal of Pure and
Applied Algebra} {\bf 172}, (2002) 165--182.
\medskip
\litem{[BH]} Bruns, W. and Herzog, J. Cohen-Macaulay Rings, {\it
Cambridge University Press}, (1993).
\medskip
\litem{[Bj]} Bjork, J.-E. Rings of Differential Operators, {\it
Amsterdam North-Holland}, (1979)
\medskip
\litem{[BS]} Brodmann, M. P. and Sharp, R. J. Local Cohomology, {\it
Cambridge studies in advanced mathematics} {\bf 60}, (1998).
\medskip
\litem{[CN]} Cowsik, R. and Nori, M. Affine curves in characteristic
$p$ are set-theoretic complete intersections, {\it Invent. Math.}
{\bf 45}, (1978) 111--114.
\medskip
\litem{[DY]} Dibaei, M. T. and Yassemi, S. Attached primes of the
top local cohomology modules with respect to an ideal, {\it Arch.
Math.} {\bf 84}, (2005) 292--297.
\medskip
\litem{[Ei]} Eisenbud, D. Commutative Algebra with A View Toward
Algebraic Geometry, {\it Springer Verlag}, (1995).
\medskip
\litem{[Go]} Goto, S. On the Cohen-Macaulayfication of certain
Buchsbaum rings, {\it Nagoya Math. J.} Vol. {\bf 80}, (1980)
107-116.
\medskip
\litem{[Gr]} Grothendieck, A. Local Cohomology, {\it Lecture Notes
in Mathematics, Springer Verlag}, (1967).
\medskip
\litem{[Ha1]} Hartshorne, R. Affine Duality and Cofiniteness, {\it
Invent. Math.} {\bf 9}, (1970) 145--164.
\medskip
\litem{[Ha2]} Hartshorne, R. Complete intersections in
characteristic $p>0$, {\it Amer. J. Math.} {\bf 101}, (1979)
380--383.
\medskip
\litem{[H1]} Hellus, M. On the set of associated primes of a local
cohomology module, {\it J. Algebra} {\bf 237}, (2001) 406--419.
\medskip
\litem{[H2]} Hellus, M. On the associated primes of Matlis duals of
top local cohomology modules, {\it Communications in Algebra} {\bf
33}, (2001), no. 11, 3997--4009.
\medskip
\litem{[H3]} Hellus, M. Matlis duals of top local cohomology modules
and the arithmetic rank of an ideal, to appear in {\it
Communications in Algebra}.
\medskip
\litem{[H4]} Hellus, M. Attached primes and Matlis duals of local
cohomology modules, submitted to {\it Archiv der Mathematik}
\medskip
\litem{[H5]} Hellus, M. Local Homology, Cohen-Macaulayness and
Cohen-Macaulayfications, to appear in {\it Algebra Colloquium}.
\medskip
\litem{[H6]} Hellus, M. Lokale Kohomologie, {\it Dissertation},
Regensburg (1999).
\medskip
\litem{[HS1]} Hellus, M. and St\"uckrad, J. Matlis duals of top
Local Cohomology Modules, submitted to {\it Proceedings of the
American Mathematical Society}.
\medskip
\litem{[HS2]} Hellus, M. and St\"uckrad, J. Generalization of an
example of Hartshorne concerning local cohomology, preprint.
\medskip
\litem{[Hu]} Huneke, C. Problems on Local Cohomology, {\it Res.
Notes Math.} {\bf 2}, (1992) 93--108.
\medskip
\litem{[HL]} Huneke, C. and Lyubeznik, G. On the vanishing of local
cohomology modules, {\it Invent. Math.} {\bf 102}, (1990) 73-93.
\medskip
\litem{[Ly1]} Lyubeznik, G. Finiteness properties of local
cohomology modules (an application of $D$-modules to Commutative
Algebra), {\it Invent. Math.} {\bf 113}, (1993), 41--55.
\medskip
\litem{[Ly2]} Lyubeznik, G. A survey of problems and results on the
number of defining equations, Commutative Algebra, {\it Math. Sci.
Res. Inst. Publ.} {\bf 15}, Springer, (1989), 375--390.
\medskip
\litem{[Ma]} Matsumura, H. Commutative ring theory, {\it Cambridge
University Press}, (1986).
\medskip
\litem{[MD]} MacDonald I. G. Secondary representation of modules
over a commutative ring, {\it Symp. Math.} {\bf XI}, (1973) 23--43.
\medskip
\litem{[Me]} Melkersson, L. Some applications of a criterion for
artinianness of a module, {\it J. Pure and Appl. Alg.} {\bf 101},
(1995) 293-303.
\medskip
\litem{[Ms]} Matlis, E. Injective modules over Noetherian rings,
{\it Pacific J. Math.} {\bf 8}, (1958) 511--528.
\medskip
\litem{[MS]} MacDonald I. G. and Sharp, R. Y. An elementary proof of
the non-vanishing of certain local cohomology modules, {\it Quart.
J. Math. Oxford} {\bf 23}, (1972) 197--204.
\medskip
\litem{[MV]} Marley, M. and Vassilev, J.C. Local cohomology modules
with infinite dimensional socles, {\it Proc. Amer. Math. Soc.} Vol.
{\bf 132}, No {\bf 12}, (2004) 3485--3490
\medskip
\litem{[Og]} Ogus, A., Local Cohomological Dimension of Algebraic
Varieties, {\it Ann. Math.} Vol. {\bf 98}, (1973) 327--365
\medskip
\litem{[Oo]} Ooishi, A. Matlis duality and width of a module, {\it
Hiroshima Math. J.} {\bf 6}, (1976) 573--587.
\medskip
\litem{[Ro]} Roberts R. N. Krull dimension for Artinian modules over
quasi local commutative rings, {\it Quart. J. Math. (Oxford)(3)}
{\bf 26}, (1975) 269--273.
\medskip
\litem{[SS]} Scheja, G. and Storch, U. Regular Sequences and
Resultants, {\it AK Peters}, (2001).
\medskip
\litem{[SV]} St\"uckrad, J. and Vogel, W. Buchsbaum Rings and
Applications. {\it VEB Deutscher Verlag der Wissenschaften, Berlin,
Germany}.
\medskip
\litem{[T1]} Tang Z. M. Local Homology and Local Cohomology, {\it
Algebra Colloquium} {\bf 11}:4, (2004) 467--476.
\medskip
\litem{[T2]} Tang Z. M. Local homology theory for Artinian modules,
{\it Comm. Alg.} {\bf 22}, (1994) 2173--2204. \vfil \eject
\parindent=0pt
{\bf Deutsche Zusammenfassung}
\bigskip
\bigskip
{\bf Lokale Kohomologie und Matlis-Dualit\"at}
\bigskip
\bigskip
Eine algebraische Menge $X$ hei{\ss}t (mengentheoretisch)
vollst\"andiger Durchschnitt, wenn sie von $\codim (X)$ vielen
algebraischen Gleichungen "ausgeschnitten" werden kann (etwa in
einem affinen oder projektiven Raum). Es ist bekannt, dass, im Falle
positiver Charakteristik, jede Kurve im $n-$dimensionalen affinen
Raum mengentheoretisch vollst\"andiger Durchschnitt ist ([CN]). Auf
der anderen Seite sind im Zusammenhang mit mengentheoretisch
vollst\"andigen Durchschnitten bemerkenswert viele Fragen
unbeantwortet. Als Beispiele seien angef\"uhrt: Ist jeder
abgeschlossene Punkt in ${\bf P}^2_{\bf Q}$ (zweidimensionaler
projektiver Raum \"uber den rationalen Zahlen) mengentheoretisch
vollst\"andiger Durchschnitt? Ist jede Kurve in ${\bf A}^3_{\bf C}$
(dreidimensionaler affiner Raum \"uber den komplexen Zahlen)
mengentheoretisch vollst\"andiger Durchschnitt? Zu diesen und vielen
weiteren verwandten Fragen enth\"alt [Ly2] eine \"Ubersicht.
\medskip Ein weiteres Beispiel ist die Kurve $C_4\subseteq {\bf
P}^3_k$, die durch
$$(u^4:u^3v:uv^3:v^4)$$
parametrisiert ist. Es ist, zumindest im Falle der Charakteristik
Null, unbekannt, ob $C_4$ mengentheoretisch vollst\"andiger
Durchschnitt ist; eine offensichtliche Obstruktion w\"are $\LCMo
^3_{I_{C_4}}(k[X_0,X_1,X_2,X_3])\neq 0$ (wobei $I_{C_4}$ das
Verschwindungsideal von $C_4$ bezeichnet). Es ist aber bekannt, dass
$$\LCMo ^3_{I_{C_4}}(k[X_0,X_1,X_2,X_3])=0$$
ist. Es ist sogar so, dass das (Nicht-)Verschwinden von lokalen
Kohomologien im Allgemeinen nicht die Minimalzahl algebraischer
Gleichungen, die die gegebene algebraische Menge "ausschneiden",
bestimmt. Algebraisch ausgedr\"uckt, bedeutet dies, dass (f\"ur ein
Ideal $I$) die Ungleichung
$$\cd (I)<\ara (I)$$
gelten kann (hier bezeichnen $\cd (I)$ die (lokale) kohomologische
Dimension von $I$ und
$$\ara(I):=\min \{ l\in \Naturalsign \vert \exists r_1,\dots ,r_l\in
R:\sqrt {I}=\sqrt {(r_1,\dots ,r_l)R}\} $$ die Minimalzahl
algebraischer Gleichungen, die die zu $I$ geh\"orende algebraische
Menge "ausschneiden"). \"Ubrigens enth\"alt 5.1 ein konkretes
Beispiel f\"ur das Vorliegen dieser Ungleichung. Auf der anderen
Seite enthalten die Matlis-Duale gewisser lokaler Kohomologiemoduln
Informationen dar\"uber, ob ein mengentheoretisch vollst\"andiger
Durchschnitt vorliegt oder nicht -- dies ist der Inhalt von
\bigskip
{\bf 1.1.4 Korollar}
\par
Seien $(R,\goth m)$ ein noetherscher lokaler Ring, $I\subsetneq R$
ein echtes Ideal, $h\in \Naturalsign $ und $\underline f=f_1,\dots
,f_h\in I$ eine $R$-regul\"are Folge. Dann sind \"aquivalent:
\par
(i) $\sqrt {\underline fR}=\sqrt I$ (d. h. $I$ ist mengentheoretisch
vollst\"andiger Durchschnitt).
\par
(ii) $\LCMo ^l_I(R)=0$ f\"ur jedes $l>h$ und $\underline f$ ist eine
$D(\LCMo ^h_I(R))$-quasi-regul\"are Folge.
\par
(ii) $\LCMo ^l_I(R)=0$ f\"ur jedes $l>h$ und $\underline f$ ist eine
$D(\LCMo ^h_I(R))$-regul\"are Folge.
\bigskip
Dieses Ergebnis legt es nahe, Matlis-Duale von lokalen
Kohomologiemoduln zu studieren, insbesondere ihre assoziierten
Primideale; dies sind auch die Hauptziele dieser Arbeit. Die
erhaltenen Ergebnisse und verwendeten Methoden f\"uhren auch zu
verschiedenen Anwendungen, die in Kapitel 6 versammelt sind.
Dar\"uber hinaus ergeben sich Zusammenh\"ange zur (lokalen)
Kohomologie formaler Schemata (7.1), zu sogenannten "attached"
Primidealen von lokalen Kohomologiemoduln (8.1, 8.2) und zum
Begriffe der lokalen Homologie (8.3, 8.4).
\medskip
Folgende Bezeichnungen seien vereinbart: Sind $R$ ein Ring,
$I\subseteq R$ ein Ideal und $M$ ein $R$-Modul, so bezeichnet $\LCMo
^l_I(M)$ die $l$-te lokale Kohomologie von $M$ mit Tr\"ager in $I$;
ist $(R,\goth m)$ ein lokaler Ring, so ist $\InjH _R(R/\goth m)$
eine (fixierte) $R$-injektive H\"ulle des $R$-Moduls $R/\goth m$.
Schlie{\ss}lich bezeichnet (\"uber dem lokalen Ring $(R,\goth m)$)
$D$ den Matlis-Dualisierungsfunktor, d. h. $D(M):=\Hom _R(M,\InjH
_R(R/\goth m))$ f\"ur jeden $R$-Modul $M$. Zur Vermeidung von
Missverst\"andnissen werden wir gegebenenfalls $D_R$ statt $D$
schreiben.\medskip Es folgt eine chronologische \"Ubersicht des
Inhalts der einzelnen Kapitel:
\medskip
Ziel von Abschnitt 1.1 ist der Beweis des eingangs zitierten
Korollars 1.1.4; dies geschieht, indem zun\"achst die folgenden
S\"atze 1.1.2 und 1.1.3 bewiesen werden, aus denen dann, im
Wesentlichen durch Spezialisierung, Korollar 1.1.4 folgt:
\bigskip
{\bf 1.1.2 Satz}
\par
Seien $(R,\goth m)$ ein noetherscher lokaler Ring, $I\subseteq R$
ein Ideal, $h\geq 1$ und $\underline f=f_1,\dots ,f_h\in I$ eine
Folge mit $\sqrt {\underline fR}=\sqrt {I}$ und so, dass
$$\LCMo ^{h-1-l}_I(R/(f_1,\dots ,f_l)R)=0\ \ (l=0,\dots ,h-3)$$
gilt (f\"ur $h\leq 2$ ist diese Bedingung leer). Dann ist
$\underline f$ eine $D(\LCMo ^h_I(R))$-quasi-regul\"are Folge.
\bigskip
{\bf 1.1.3 Satz}
\par
Seien $(R,\goth m)$ ein noetherscher lokaler Ring, $I\subseteq R$
ein Ideal, $h\geq 1$ und $\underline f=f_1,\dots ,f_h\in I$ so, dass
$$\LCMo ^l_I(R)=0\ \ (l>h)$$
und
$$\LCMo ^{h-1-l}_I(R/(f_1,\dots ,f_h)R)=0\ \ (l=0,\dots ,h-2)$$
gelten (f\"ur $h\leq 1$ ist diese Bedingung leer) und so, dass
$\underline f$ eine $D(\LCMo ^h_I(R))$-quasi-regul\"are Folge ist.
Dann gilt $\sqrt I=\sqrt{(f_1,\dots ,f_h)R}$.
\bigskip
Korollar 1.1.4 (siehe oben) legt es nahe, zu untersuchen, f\"ur
welche $f\in R$ die Multiplikation mit $x$ auf einem Matlis-Dual
eines lokalen Kohomologiemoduls injektiv ist, mit anderen Worten,
die Menge der Nullteiler auf einem solchen Modul zu bestimmen. Eine
genauere Frage ist die nach der Menge der zu diesem Modul
assoziierten Primideale. In diesem Zusammenhang verweisen wir auf
\bigskip
{\bf 1.2.2 Vermutung}
\par
Sind $(R,\goth m)$ ein noetherscher lokaler Ring, $h>0$ und
$x_1,\dots ,x_h$ Elemente von $R$, so gilt
$$\Ass _R(D(\LCMo ^h_{(x_1,\dots ,x_h)R}(R)))=\{ \goth p\in \Spec
(R)\vert \LCMo ^h_{(x_1,\dots ,x_h)R}(R/\goth p)\neq 0\} \ \ .$$
Diese Vermutung bezeichnen wir mit (*). Die Inklusion $\subseteq $
ist stets richtig, dies ist (unter anderem) der Inhalt von
\bigskip
{\bf 1.2.1 Bemerkung}
\par
Sind $(R,\goth m)$ ein noetherscher lokaler Ring, $h>0$ und
$x_1,\dots ,x_h$ Elemente von $R$, so gilt
$$\Ass _R(D(\LCMo ^h_{(x_1,\dots ,x_h)R}(R)))\subseteq \{ \goth p\in \Spec
(R)\vert \LCMo ^h_{(x_1,\dots ,x_h)R}(R/\goth p)\neq 0\} \ \ .$$ Es
gibt zu (*) \"aquivalente Aussagen:
\bigskip
{\bf 1.2.3 Satz} \par Die folgenden Aussagen sind \"aquivalent:
\par
(i) Vermutung (*) ist richtig, d. h. f\"{u}r jeden noetherschen lokalen
Ring $(R,\goth m)$, jedes $h>0$ und jede Folge $x_1,\dots ,x_h\in R$
gilt
$$\Ass _R(D(\LCMo ^h_{(x_1,\dots ,x_h)R}(R)))=\{ \goth p\in \Spec
(R)\vert \LCMo ^h_{(x_1,\dots ,x_h)R}(R/\goth p)\neq 0\} \ \ .$$
(ii) F\"ur jeden noetherschen lokalen Ring $(R,\goth m)$, jedes
$h>0$ und jede Folge $x_1,\dots ,x_h\in R$ ist die Menge
$$Y:=\Ass _R(D(\LCMo ^h_{(x_1,\dots ,x_h)}(R)))$$
abgeschlossen unter Generalisierung, d. h. aus $\goth p_0,\goth
p_1\in \Spec(R),\goth p_0\subseteq \goth p_1,\goth p_1\in Y$ folgt
$\goth p_0\in Y$.
\par
(iii) F\"ur jeden noetherschen lokalen Integrit\"atsring $(R,\goth
m)$, jedes $h>0$ und jede Folge $x_1,\dots ,x_h\in R$ gilt die
Implikation
$$\LCMo ^h_{(x_1,\dots ,x_h)}(R)\neq 0\Longrightarrow \{ 0\} \in
\Ass _R(D(\LCMo ^h_{(x_1,\dots ,x_h)R}(R)))\ \ .$$ (iv) F\"ur jeden
noetherschen lokalen Ring $(R,\goth m)$, jeden endlich erzeugten
$R$-Modul $M$, jedes $h>0$ und jede Folge $x_1,\dots ,x_h\in R$ gilt
die Gleichheit
$$\Ass _R(D(\LCMo ^h_{(x_1,\dots ,x_h)R}(M)))=\{ \goth p\in
\Supp_R(M)\vert \LCMo ^h_{(x_1,\dots ,x_h)R}(M/\goth pM)\neq 0\} \ \
.$$
\bigskip
Aussage (iv) ist also formal allgemeiner als Aussage (i), aber
inhaltlich \"aquivalent dazu.
\medskip
[HS1, Kapitel 0] enth\"alt eine weitere Vermutung zur Struktur der
Menge der assoziierten Primideale von $D(\LCMo ^h_{(x_1,\dots
,x_h)R}(R))$: Alle Primideale $\goth p$, die maximal in $\Ass
_R(D(\LCMo ^h_{(x_1,\dots ,x_h)R}(R)))$ sind, haben die Dimension
$h$: $\dim (R/\goth p)=h$; diese Vermutung ist falsch, Bemerkung
1.2.4 enth\"alt ein Gegenbeispiel (mit $\dim (R)-h=2$).
\medskip
Indem wir uns mit (quasi-)regul\"aren Folgen auf Moduln der Form
$D(\LCMo ^h_I(R))$ besch\"aftigen, stellt sich folgende Frage: Im
allgemeinen ist $D(\LCMo ^h_I(R))$ nicht endlich erzeugt (viele
Ergebnisse dieser Arbeit zeigen, dass dieser Modul im Allgemeinen
unendlich viele assoziierte Primideale hat), der Begriff der
regul\"aren Folge auf nicht-endlichen Moduln l\"asst manche
Eigenschaften vermissen: Beispielsweise gilt (\"uber einem lokalen
noetherschen Ring ($R,\goth m)$) f\"ur einen endlichen $R$-Modul $M$
und eine $M$-regul\"are Folge $r_1,\dots ,r_h\in R$, dass auch die
Folge $r_1^\prime ,\dots ,r_h^\prime \in R$ $M$-regul\"ar ist, falls
nur $(r_1,\dots ,r_h)R=(r_1^\prime ,\dots ,r_h^\prime )R$
vorausgesetzt ist; f\"ur nicht-endliche Moduln stimmt diese Aussage
im Allgemeinen nicht. Die eingangs erw\"ahnte Frage lautet: Stimmt
die Aussage f\"ur Moduln der Form $D(\LCMo ^h_{(x_1,\dots
,x_h)R}(R))$? Die Antwort ist (unter gewissen Voraussetzungen)
positiv: \vfil \eject {\bf 1.3.1 Satz}
\medskip
Seien $(R,\goth m)$ ein noetherscher lokaler Ring, $h\geq 1$ und
$I\subseteq R$ ein Ideal mit $\LCMo ^h_I(R)\neq 0\iff l=h$. Weiter
seien $1\leq h^\prime \leq h$ und $r_1,\dots ,r_{h^\prime }\in I$
eine $R$-regul\"are Folge, die auch $D(\LCMo ^h_I(R))$-regul\"ar
ist. Es seien $r_1^\prime, \dots ,r_{h^\prime }^\prime \in I$ mit
$(r_1,\dots ,r_{h^\prime })R=(r_1^\prime ,\dots ,r_{h^\prime
}^\prime )R$. Dann ist auch $r_1^\prime ,\dots ,r_{h^\prime }^\prime
$ eine $D(\LCMo ^h_I(R))$-regul\"are Folge.
\medskip
In Abschnitt 1.4 ist $R_0$ ein lokaler Unterring von $R$ und wir
untersuchen Beziehungen zwischen
$$D_R(\LCMo ^i_{(y_1,\dots
,y_i)R}(R))$$ und
$$D_{R_0}(\LCMo ^i_{(y_1,\dots ,y_i)R}(R))\ \ :$$
Ein Ergebnis ist
\bigskip
{\bf 1.4.3 Bemerkung (ii), zweite Aussage}
\par
Seien $(R,\goth m)$ ein noetherscher lokaler
\"aquicharakteristischer kompletter Ring mit Koeffizientenk\"orper
$k$ und $\underline y=y_1,\dots ,y_i\in R$ eine Folge in $R$ so,
dass $R_0:=k[[y_1,\dots ,y_i]]\ (\subseteq R)$ regul\"ar und
$i$-dimensional ist (dies ist z. B. der Fall, wenn $\LCMo
^i_{(y_1,\dots ,y_i)R}(R)\neq 0$ ist). Wenn Vermutung (*) richtig
ist, gilt
$$\Ass _R(D_{R_0}(\LCMo ^i_{(y_1,\dots ,y_i)R}(R)))=\Ass _R(D_R(\LCMo
^i_{(y_1,\dots ,y_i)R}(R)))\ \ .$$ In Kapitel 2 werden Eigenschaften
der Menge
$$\Ass _R(D(\LCMo ^i_{(x_1,\dots ,x_i)R}(R)))$$
untersucht ($(R,\goth m)$ ein noetherscher lokaler Ring, $\underline
x=x_1,\dots ,x_i$ eine Folge in $R$). Die verwendeten Methoden sind
konstruktiv in dem Sinne, dass zun\"achst in dem $R$-Modul
$$E:=k[X_1^{-1},\dots ,X_n^{-1}]$$
($k$ ein K\"orper) gewisse Elemente konstruiert werden (Lemmata 2.1
-- 2.3); bekanntlich ist $E$ eine $R$-injektive H\"ulle von $k$,
falls $R=k[[X_1,\dots ,X_n]]$ eine formale Potenzreihenalgebra
\"uber $k$ ist. Ein zentrales Ergebnis in diesem Kapitel (und eine
Folgerung aus Lemma 2.1) ist
\bigskip
{\bf 2.4 Satz} \par Seien $(R,\goth m)$ ein noetherscher lokaler
\"aquicharakteristischer Ring, $i\geq 1$ und $x_1,\dots ,x_i$ eine
Folge in $R$. Dann ist
$$\{ \goth p\in \Spec (R)\vert x_1,\dots ,x_i\hbox { ist Teil eines Paramtersystem von }R/\goth p\} \subseteq \Ass _R(D(\LCMo ^i_{(x_1,\dots
,x_i)R}(R)))\ \ .$$ Satz 2.5 enth\"alt ein \"ahnliches Ergebnis im
gemischt-charakteristischen Fall.
\medskip
Satz 2.4 erm\"oglicht es, im Falle $i=1$ die Menge der assoziierten
Primideale vollst\"andig zu berechnen:
\bigskip
{\bf 2.6 Korollar}
\par
Seien $(R,\goth m)$ ein noetherscher lokaler
\"aquicharakteristischer Ring und $x\in R$. Dann ist
$$\Ass _R(D(\LCMo ^1_{xR}(R)))=\Spec (R)\setminus ({\cal V}x)\ \ .$$
Insbesondere ist die Menge der assoziierten Primideale des
Matlis-Duals einen lokalen Kohomologiemoduls im Allgemeinen nicht
endlich. \bigskip
Andererseits zeigen wir in Bemerkung 2.7 (ii), dass
die in Satz 2.4 bewiesene Inklusion im Allgemeinen echt ist, dass
also nicht alle assoziierten Primideale von $D(\LCMo ^i_{(x_1,\dots
,x_i)R}(R))$ von der in Satz 2.4 angegebenen Form sind.
Schlie{\ss}lich untersuchen wir (in Bemerkung 2.7 (iii)) die
Teilmengen
$$Z_1:=\{ \goth p\in \Spec (R)\vert \LCMo ^i_{(x_1,\dots ,x_i)R}(R/\goth
p)\neq 0\} $$ und
$$Z_2:=\{ \goth p\in \Spec (R)\vert x_1,\dots ,x_i\hbox { ist Teil eines Parametersystems von }R/\goth p\}$$
von $\Spec (R)$ im Hinblick auf ihre Abgeschlossenheit unter
Generalisierung (man beachte, dass gem\"a{\ss} Satz 2.4 (bzw. 2.5)
und Bemerkung 1.1.2)
$$Z_2\subseteq \Ass _R(D(\LCMo ^i_{(x_1,\dots ,x_i)R}(R)))\subseteq
Z_1$$ gilt). Dabei zeigt sich, dass $Z_1$ abgeschlossen unter
Generalisierung ist, $Z_2$ hingegen im Allgemeinen nicht, selbst
dann nicht, wenn $R$ regul\"ar ist. Immerhin gilt die schw\"achere
Aussage
$$Z_2\neq \emptyset \Longrightarrow \{ 0\} \in Z_2\ \ .$$
In Kapitel 3 wird die Untersuchung von $\Ass _R(D(\LCMo
^i_{(x_1,\dots ,x_i)R}(R)))$ fortgesetzt, wobei nun keine
Voraussetzungen \"uber die (Gleich-)Charakteristik gemacht werden.
Dabei ist nachfolgendes Lemma ein entscheidender Ausgangspunkt:
\bigskip
{\bf 3.1.1 Lemma}
\par
Seien $R$ ein Ring, $x,y\in R$ und $U$ ein $R$-Untermodul von $R_x$
mit $\im (\iota _x)\subseteq U$ ($\iota: R\to R_x$ bezeichnet die
kanonische Abbildung). Weiter bezeichne $S:=\im (\iota _y)\subseteq
R_y$. Dann existiert ein Epimorphismus
$$R_x/U\to R_{xy}/(S_x+U_y)$$
von $R$-Moduln. \bigskip Daraus folgt unter Verwendung von
\v{C}ech-Kohomologie leicht
\bigskip
{\bf 3.1.2 Satz} \par
Seien $R$ ein noetherscher Ring, $x_1,\dots
,x_m,y_1,\dots ,y_n\in R$ ($m\in \Naturalsign ^+, n\in \Naturalsign
$) und $M$ ein $R$-Modul. Dann existiert ein Epimorphism
$$\LCMo ^m_{(x_1,\dots ,x_m)R}(R)\to \LCMo ^{m+n}_{(x_1,\dots
,x_m,y_1,\dots ,y_n)R}(R)$$ von $R$-Moduln.
\bigskip
Die Idee ist nun, diesen Epimorphismus zu dualisieren; man erh\"alt
einen Monomorphismus und folglich eine Inklusionsbeziehung zwischen
den jeweiligen Mengen von assoziierten Primidealen. Wir erhalten:
\vfil \eject {\bf 3.1.3 Satz}
\par
Seien $(R,\goth m)$ ein noetherscher lokaler Ring, $m\in
\Naturalsign ^+$, $x_1,\dots ,x_m\in R$ und $M$ ein endlich
erzeugter $R$-Modul. Dann gelten:
\par
(i) F\"ur jedes $\goth p\in \Ass _R(D(\LCMo ^m_{(x_1,\dots
,x_m)R}(M)))$ ist $\dim (M/\goth pM)\geq m$.
\par
(ii) $\{ \goth p\in \Supp _R(M)\vert x_1,\dots ,x_m\hbox{ ist Teil
eines Parametersystem von }R/\goth p\} \subseteq \Ass _R(D(\LCMo
^m_{(x_1,\dots ,x_m)R}(M)))$.
\par
(iii) F\"ur jedes $x\in R$ gilt $\Ass _R(D(\LCMo ^1_{xR}(R)))=\Spec
(R)\setminus {\cal V}(x)$.
\par
(iv) Ist $x_1,\dots ,x_m$ Teil eines Parametersystems von $M$, so
gilt $\Assh (M)\subseteq \Ass _R(D(\LCMo ^m_{(x_1,\dots
,x_m)R}(M)))$; im Falle $m=\dim (M)$ gilt sogar Gleichheit: $\Assh
(M)=\Ass _R(D(\LCMo ^m_{(x_1,\dots ,x_m)R}(M)))$ (dabei ist
$\Assh(M)$ definiert als die Menge der h\"ochstdimensionalen zu $M$
assoziierten Primideale).
\par
(v) Falls $R$ komplett ist, gilt f\"ur jeded $\goth p\in \Supp
_R(M)$ mit $\dim (R/\goth p)=m$ die \"Aquivalenz
$$\goth p\in \Ass _R(D(\LCMo ^m_{(x_1,\dots ,x_m)R}(M)))\iff x_1,\dots
,x_m\hbox { ist ein Parametersystem von }R/\goth p\ \ .$$ Im
Abschnitt 3.2 wird die Menge
$$\Ass _R(D(\LCMo ^{\dim (R)-1}_I(R)))$$
untersucht; dabei ist $I$ (zun\"achst) ein beliebiges Ideal von $R$,
wir setzen also nicht voraus, dass $I$ (bis auf Radikal) von $\dim
(R)-1$ Elementen erzeugt wird. Die wichtigsten Ergebnisse sind die
beiden folgenden S\"atze:
\bigskip
{\bf 3.2.6 Satz}
\par
Seien $(R,\goth m)$ ein noetherscher lokaler $d$-dimensionaler Ring
und $J\subseteq R$ ein Ideal mit $\dim (R/J)=1$ und $\LCMo
^d_J(R)=0$. Dann gilt
$$\Assh (D(\LCMo ^{d-1}_J(R)))=\Assh (R)\ \ .$$
\bigskip
{\bf 3.2.7 Satz}
\par
Seien $(R,\goth m)$ ein noetherscher lokaler kompletter
$d$-dimensionaler Ring und $J\subseteq R$ ein Ideal mit $\dim
(R/J)=1$ und $\LCMo ^d_J(R)=0$. Dann gilt
$$\Ass _R(D(\LCMo ^{d-1}_J(R))=\{ P\in \Spec(R)\vert \dim (R/P)=d-1,\dim
(R/(P+J))=0\} \cup \Assh (R)\ \ .$$ Die Beweise sind etwas technisch
und beruhen, unter anderem, auf \bigskip {\bf 3.2.1 Lemma}
\par
Seien $(S,\goth m)$ ein noetherscher lokaler kompletter
Gorenstein-Ring der Dimension $n+1$ und $\goth P\subseteq S$ ein
Primideal der H\"ohe $n$. Dann gilt kanonisch
$$D(\LCMo ^n_\goth P(S))=\widehat {S_\goth P}/S\ \ .$$
In Kapitel 4 untersuchen wir einen Spezialfall, den wie als \lq\lq
regul\"aren Fall\rq\rq bezeichnen: $k$ ein K\"orper, $R=k[[X_1,\dots
,X_n]]$ eine Potenzreihenalgebra \"uber $k$ in $n$ Variablen und $I$
das Ideal $(X_1,\dots ,X_h)R$ von $R$ ($1\leq h\leq n$). Zum Beweis
von Vermutung (*) kann man sich auf den regul\"aren Fall
zur\"uckziehen:
\bigskip
{\bf 4.1.2 Satz}
\par
Seien $(R,\goth m)$ ein noetherscher lokaler kompletter Ring mit
einem Koeffizientenk\"orper $k$, $l\in \Naturalsign ^+$ und
$x_1,\dots ,x_l$ Teil eines Parametersystems von $R$. $I:=(x_1,\dots
,x_l)R$. Seien $x_{l+1},\dots ,x_d\in R$ so, dass $x_1,\dots ,x_d$
ein Parametersystem von $R$ ist. $R_0$ bezeichne den
($d$-dimensionalen, regul\"aren) Unterring $k[[x_1,\dots ,x_d]]$ von
$R$. Ist $\Ass _{R_0}(D(\LCMo ^l_{(x_1,\dots ,x_l)R_0}(R_0)))$
abgeschlossen unter Generalisierung, so auch $\Ass _R(D(\LCMo
^l_{(x_1,\dots ,x_l)R}(R)))$.
\bigskip
In Abschnitt 4.2 behandeln wir den regul\"aren Fall. Satz 4.2.1
fasst (im Wesentlichen) zusammen, was die bisher gezeigten S\"{a}tze im
regul\"aren Fall bedeuten. Ein weiteres Ergebnis ist
\bigskip
{\bf 4.2.3 Satz}
\par
Seien $(R_0,\goth m_0)$ ein noetherscher lokaler kompletter
\"aquicharakteristischer Ring, $\dim (R_0)=n-1$, $k\subseteq R_0$
ein Koeffizientenk\"orper und $h\in \{ 1,\dots ,n\} $. Weiter seien
$x_1,\dots ,x_n\in R$ Elemente mit $\sqrt {(x_1,\dots ,x_n)R}=\sqrt
{\goth m_0}$. $I_0:=(x_1,\dots ,x_h)R_0$. $R:=k[[X_1,\dots ,X_n]]$
sei eine Potenzreihenalgebra \"uber $k$ in $n$ Unbestimmten,
$I:=(X_1,\dots ,X_h)R$. Der durch $X_i\mapsto x_i$ ($i=1,\dots ,n$)
festegelegte $k$-Algebrahomomorphismus $R\to R_0$ induziert einen
modul-endlichen Homomorphismus $\iota :R/fR\to R_0$ mit einem
geeigneten Primideal $f$ von $R$. Wir setzen
$$D:=D(\LCMo ^h_I(R))\ \ .$$
Dann gelten:
\par
(i) $D$ hat ein $f$ enthaltendes assoziiertes Primideal genau dann,
$\LCMo ^h_{I_0}(R_0)\neq 0$ gilt.\smallskip Wenn wir (zus\"atzlich)
annehmen, dass $R_0$ regul\"ar ist und dass $\height (I_0)<h$ ist,
so gelten:
\smallskip
(ii) Es gibt keine zu $D$ assoziiertes Primideal, dass $f$ enth\"alt
und die H\"ohe $n-h$ hat.
\par
(iii) Wenn $\LCMo ^h_{I_0}(R_0)\neq 0$ ist (dann ist $f$ in einem zu
$D$ assoziierten Primideal enthalten), gilt $\dim (R/\goth q)>h$
f\"ur jedes $f$ enthaltende maximale Element in $\Ass _R(D)$.
\bigskip
Aussage (iii) h\"angt eng mit dem Gegenbeispiel zur Vermutung (+)
aus [HS1, Kapitel 0] zusammen -- vgl. dazu den Abschnitt nach Satz
1.2.3 in dieser Zusammenfassung. \medskip Abschnitt 4.3 behandelt
den Fall $h=n-2$ (in den F\"allen $h=n-1$ und $h=n$ wurde $\Ass
_R(D(\LCMo ^h_I(R)))$ vollst\"andig bestimmt -- vgl. Satz 4.2.1):
Unter anderem zeigen wir:
\bigskip
{\bf 4.3.1 Korollar}
\par
In der Situation von Satz 4.2.3 seien $R_0$ regul\"ar und $\height
(I_0)<n-2=:h$. Dann gilt
$$fR\in \Ass _R(D)\iff \LCMo ^{n-2}_{I_0}(R_0)\neq 0\ \ .$$
Falls dieses Bedingungen zutreffen, ist $fR$ maximal in $\Ass
_R(D)$.
\bigskip
Bekanntlich (vgl. [HL, Theorem 2.9) ist $\LCMo ^{n-2}_{I_0}(R_0)$
genau dann trivial, wenn $\Spec (R_0/I_0)\setminus \{ \goth
m_0/I_0\} $ formal-geometrisch zusammenh\"angend ist.
\medskip
Der folgende Satz ist f\"ur sich genommen interessant und wird sich
sp\"ater (im Abschnitt 6.2: Verallgemeinerung eines Beispiels von
Hartshorne) als n\"utzlich erweisen:
\bigskip
{\bf 4.3.4 Satz}
\par
Es seien $k$ ein K\"orper, $R=k[[X_1,\dots ,X_n]]$ ($n\geq 3$) eine
Potenzreihenalgebra \"uber $k$ in $n$ Variablen und $I$ das Ideal
$(X_1,\dots ,X_{n-2})R$. Ausserdem sei $p\in R$ ein Primelement mit
$p\in I\cap (X_{n-1},X_n)R$. Dann ist die Menge
$$\{ \goth p\in \Spec (R)\vert \goth p\in \Ass _R(D(\LCMo ^{n-2}_I(R))),p\in
\goth p, \height(\goth p)=2\} $$ unendlich.
\medskip
Kapitel 5 behandelt die Frage: Was bedeutet es, dass ein gegebenes
Ideal arithmetischen Rang eins oder zwei hat? Unter anderem werden
Kriterien f\"ur diese Bedingungen bewiesen. Zu Beginn jedoch
pr\"asentieren wir ein Beispiel, bei dem arithmetischer Rang und
kohomologische Dimension nicht \"ubereinstimmen:
\bigskip
{\bf 5.1 Beispiel}
\par
Seien $k$ ein K\"orper und $R=k[[x,y,z,w]]$ eine Potenzreihenalgebra
\"uber $k$ in $4$ Variablen. Es sei
$$I:=\sqrt {(xw-yz,y^3-x^2z,z^3-w^2y)R}\ \ .$$
Dann gilt
$$\cd (I/(xw-yz)R)=1\neq 2=\ara (I/(xw-yz)R)\ \ .$$
Die Hauptergebnisse in Abschnitt 5.2 sind Kriterien f\"ur $\ara
(I)\leq 1$ bzw. $\ara (I)\leq 2$:
\bigskip
{\bf Definition}
\par
Seien $(R,\goth m)$ ein noetherscher lokaler Ring und $X$ eine
Teilmenge von $\Spec (R)$. Wir sagen, dass $X$ Primvermeidung
erf\"ullt, wenn f\"ur jedes Ideal $J$ von $R$ die Implikation
$$J\subseteq \bigcup _{\goth p\in X}\goth p\Longrightarrow \exists
\goth p_0\in X:J\subseteq \goth p_0$$ gilt.
\bigskip
{\bf 5.2.5 Satz (i)}
\par
Es sei $I$ ein Ideal in einem noetherschen lokalen Ring $(R,\goth
m)$ mit $0=\LCMo ^2_I(R)=\LCMo ^3_I(R)=\dots $. Dann gilt:
$$\ara (I)\leq 1\iff \Ass _R(D(\LCMo ^1_I(R)))\hbox { erf\"ullt
Primvermeidung}\ \ .$$
\bigskip
{\bf 5.2.6 Korollar (i)}
\par
Sei $I$ ein Ideal in einem noetherschen lokalen Ring $(R,\goth m)$.
Genau dann gilt $\ara (I)\leq 2$, wenn ein $g\in I$ existiert mit
$0=\LCMo ^2_I(R/gR)=\LCMo ^3_I(R/gR)=\dots $ und so, dass $\Ass
_R(D(\LCMo ^1_I(R/gR)))$ Primvermeidung erf\"ullt.
\bigskip
Zu beiden Kriterien gibt es analoge Aussagen im graduierten Fall
(Satz 5.2.5 (ii) und Korollar 5.2.6 (ii)), auch die Beweise sind
analog.
\medskip
Abschnitt 5.3 behandelt subtile Unterschiede zwischen der
graduierten und der lokalen Situation.
\medskip
Kapitel 6 enth\"alt verschiedene Anwendungen der in den
vorangehenden Kapitel entwickelten Theorie: Als erste Anwendung
verweisen wir auf zwei neue Beweise (Satz 6.1.2 und Satz 6.1.4) des
Satzes von Hartshorne-Lichtenbaum; der besagt bekanntlich, dass
f\"ur einen noetherschen lokalen kompletten Integrit\"atsring
$(R,\goth m)$ und ein Ideal $I\subseteq R$ genau dann $\LCMo ^{\dim
(R)}_I(R)\neq 0$ gilt, wenn $\sqrt I=\goth m$ ist. Der Beweis von
6.1.2 verwendet die Normalisierung von $R$ und Matlis-Duale von
lokalen Kohomologiemoduln; beim zweiten Beweis (6.1.4) verwenden wir
die Tatsache, dass \"uber einem noetherschen lokalen kompletten
Gorenstein-Ring $(S,\goth m)$ mit $\dim (S)=n+1$ f\"ur jedes
Primideal $\goth P$ von $S$ der H\"ohe $n$
$$D(\LCMo ^n_\goth P(S))=\widehat {S_\goth P}/S$$
gilt (dies ist Lemma 3.2.1); besonders bemerkenswert ist dabei wohl,
dass der Beweis von 6.1.4 die Ring-Struktur von $\widehat {S_\goth
P}$ verwendet (n\"amlich im Beweis von Lemma 6.1.3).
\medskip
Hartshorne ([Ha1, section 3]) untersuchte (im Wesentlichen)
folgendes Beispiel: Seien $k$ ein K\"orper, $R=k[[X_1,X_2,X_3,X_4]]$
eine Potenzreihenalgebra \"uber $k$ in vier Unbestimmten,
$I=(X_1,X_2)R$ und $p=X_1X_4+X_2X_3\in R$. Dann ist $\Supp _R(\LCMo
^2_I(R/pR))=\{ \goth m\} $, aber $\LCMo ^2_I(R/pR)$ ist nicht
artinsch als $R$-Modul. In Abschnitt 6.2 der vorliegenden Arbeit
zeigen wir zun\"achst, dass $D(\LCMo ^2_I(R/pR))$ unendliche viele
assoziierte Primideale hat; somit ist $\LCMo ^2_I(R/pR)$ nicht
artinsch. Mit anderen Worten: Die Untersuchung der assoziierten
Primideale von $D(\LCMo ^2_I(R/pR))$ f\"uhrt zu einem einfachen
Beweis der Tatsache, dass $\LCMo ^2_I(R/pR)$ nicht artinsch ist.
Diese Idee wird nun verallgemeinert zu:
\bigskip
{\bf 6.2.3 Satz}
\par
Seien $k$ ein K\"orper, $n\geq $, $R=k[[X_1,\dots ,X_n]]$,
$I=(X_1,\dots ,X_{n-2})R$ und $p$ ein Primelement in $R$ mit $p\in
(X_{n-1},X_n)R$. Dann ist
$$\LCMo ^{n-2}_I(R/pR)$$
nicht artinsch.
\bigskip
Auch Marley und Vassilev ([MV, theorem 2.3]) haben Hartshornes
Beispiel verallgemeinert; man kann [MV, theorem 2.3] und Satz 6.2.3
nur in einem Spezialfall vergleichen: Dies machen wir in Bemerkung
6.2.5 und erhalten als Ergebnis, dass (in diesem Spezialfall) Satz
6.2.3 mit schw\"acheren Voraussetzungen auskommt als [MV, theorem
2.3].
\medskip
Im Abschnitt 6.3 ist $(R,\goth m)$ ein noetherscher lokaler Ring.
Ist nun $I=(x_1,\dots ,x_i)R\subseteq R$ ein Ideal, das
mengentheoretisch vollst\"andiger Durchschnitt ist (im Sinne von
$\height (I)=i$), so ergibt sich sofort $\LCMo ^i_I(R)\neq 0$, z. B.
indem man lokalisiert. Hauptergebnis in 6.3 ist nun eine gewisse
notwendige Bedingung f\"ur $\LCMo ^i_I(R)\neq 0$:
\bigskip
{\bf 6.3.1 Satz (partiell)}
\par
Seien $(R,\goth m)$ ein noetherscher lokaler kompletter
Integrit\"atsring, der einen K\"orper $k$ enthalte und $x_1,\dots
,x_i\in R$ ($i\geq 1$) eine Folge in $R$. Bezeichne $R_0$ den
Unterring $k[[x_1,\dots ,x_i]]$ von $R$. Dann gilt die Implikation
$$\LCMo ^i_I(R)\neq 0\Longrightarrow R\cap Q(R_0)=R_0$$
(dabei ist $Q(R_0)$ der Quotientenk\"orper von $R_0$ und die
Durchschnittsbildung ist in $Q(R)$ gemeint).
\medskip
\"Uber gewissen Ringen (z. B. kompletten Cohen-Macaulay Ringen) gibt
es eine Korrespondenz zwischen $\Ext $-Moduln auf der einen und
lokalen Kohomologiemoduln auf der anderen Seite; diese wird als
lokale Dualit\"at bezeichnet, vgl. dazu etwa [BS, section 11]. Alle
in dieser Korrespondenz vorkommenden lokalen Kohomologiemoduln haben
$\goth m$ als Tr\"agerideal. In Abschnitt 6.4 verallgemeinern wir
dieses Prinzip auf beliebige Ideale: \bigskip {\bf 6.4.1 Satz}
\par
Seien $(R,\goth m)$ ein noetherscher lokaler Ring, $I\subseteq R$
ein Ideal, $h\in \Naturalsign $ mit
$$\LCMo ^l_I(R)\neq 0\iff l=h$$
und sei $M$ ein $R$-Modul. Dann gilt f\"ur jedes $i\in \{ 0,\dots
,h\} $ kanonisch
$$\Ext ^i_R(M,D(\LCMo ^h_I(R)))=D(\LCMo ^{h-i}_I(M))\ \ .$$
Die nachfolgende Bemerkung 6.4.2 zeigt, dass Satz 6.4.1 wirklich
eine verallgemeinerte lokale Dualit\"at ist.
\medskip
In Abschnitt 7.2 zeigen wir, dass $D(\LCMo ^i_I(R))$ eine kanonische
$D$-Modul-Struktur hat; damit ist folgendes gemeint: Seien $k$ ein
K\"orper und $R=k[[X_1,\dots ,X_n]]$ eine Potenzreihenalgebra \"uber
$k$ in $n$ Unbestimmten. Sei
$$D(R,k)\subseteq \End _k(R)$$
der (nicht-kommutative) Unterring, der von allen Multiplikationen
mit allen Elementen aus $R$ und allen $k$-linearen Derivationen
erzeugt wird. $D:=D(R,k)$ wird als Ring der $k$-linearen
Derivationen auf $R$ bezeichnet. (Links-)$D$-Moduln im Zusammenhamg
mit lokaler Kohomologie wurden in [Ly1] studiert; darin wurde auch
gezeigt, dass (in sehr allgemeinen Situationen) lokale
Kohomologiemoduln eine kanonische (Links-)$D$-Modul-Struktur tragen.
Wir zeigen nun (in 7.2), dass f\"ur jedes Ideal $I\subseteq
R=k[[X_1,\dots ,X_n]]$ und f\"ur jedes $i\in \Naturalsign $ auch
$$D(\LCMo ^i_I(R))$$
eine kanonische (Links-)$D$-Modul-Struktur hat; weiter zeigen wir,
dass $D(\LCMo ^i_I(R))$ als $D$-Modul im Allgemeinen nicht endlich
erzeugt ist, insbesondere nicht holonom (siehe [Bj, sections 1,3]
f\"ur den Begriff der holonomen $D$-Moduln).
\medskip
Seien $(R,\goth m)$ ein noetherscher lokaler Integrit\"atsring und
$x_1,\dots ,x_i\in R$ ($i\geq 1$). In zahlreichen Situationen (vgl.
etwa Satz 3.1.3 (ii)) gilt dann
$$\{ 0\} \in \Ass _R(D(\LCMo ^i_{(x_1,\dots ,x_i)R}(R)))$$
(sogar immer falls $\LCMo ^i_{(x_1,\dots ,x_i)R}(R)\neq 0$, wenn
Vermutung (*) zutrifft). Es ist nat\"urlich, nach der
$Q(R)$-Vektorraum-Dimension von
$$D(\LCMo ^i_{(x_1,\dots ,x_i)R}(R))\otimes _RQ(R)$$
zu fragen (dies ist eine sogenannte Bass-Zahl von $D(\LCMo
^i_{(x_1,\dots ,x_i)R}(R))$). Es zeigt sich, dass diese Dimension im
Allgemeinen nicht endlich ist; genauer gilt:
\bigskip
{\bf 7.3.2 Satz}
\par
Seien $k$ ein K\"orper und $R=k[[X_1,\dots ,X_n]]$ eine
Potenzreihenalgebra \"uber $k$ in $n\geq 2$ Unbestimmten, $1\leq
i<n$ und $I$ das Ideal $(X_1,\dots ,X_i)R$ von $R$. Dann ist
$$\dim _{Q(R)}(D(\LCMo ^i_I(R))\otimes _RQ(R))=\infty \ \ .$$
In Abschnitt 7.4 untersuchen wir Moduln der Form $\LCMo ^h_I(D(\LCMo
^h_I(R)))$. Das Hauptergebnis ist \bigskip {\bf 7.4.1 Satz und 7.4.2
Satz (Spezialfall)}
\par
Seien $(R,\goth m)$ ein noetherscher lokaler kompletter regul\"arer
Ring der \"Aquicharakteristik Null, $I\subseteq R$ ein Ideal der
H\"ohe $h\geq 1$ mit $\LCMo ^l_I(R)=0$ ($l>h$); weiter sei
$\underline x=x_1,\dots ,x_h$ eine $R$-regul\"are Folge in $I$. Dann
ist
$$\LCMo ^h_I(D(\LCMo ^h_I(R)))$$
entweder gleich Null oder isomorph zu $\InjH _R(R/\goth m)$. Im
Falle $I=(x_1,\dots ,x_h)R$ trifft letzteres zu, i. e.
$$\LCMo ^h_I(D(\LCMo ^h_I(R)))=\InjH _R(R/\goth m)\ \ .$$
In den Abschnitten 8.1 und 8.2 werden sogenannte "attached"
Primideale studiert, und zwar im Hinblick auf lokale
Kohomologiemoduln; 8.1 versammelt zahlreiche grundlegende (und
teilweise nat\"urlich bekannte) Eigenschaften von "attached"
Primidealen, 8.2 enth\"alt unsere Ergebnisse, d. h. Informationen
\"uber "attached" Primideale von lokalen Kohomoliemoduln. Schon in
8.1 zeigt sich ein enger Zusammenhang zwischen assoziierten den
assoziierten Primidealen vom Matlis-Dual eines $R$-Moduls $M$
einerseits und den "attached" Primidealen von $M$ andererseits. Hier
eine Auswahl unserer Ergebnisse:
\bigskip
{\bf 8.2.1 Satz}
\par
Seien $(R,\goth m)$ ein noetherscher lokaler Ring und $M$ ein
endlich erzeugter $n$-dimensionaler $R$-Modul. Dann gilt
$$\Att _R(\LCMo ^n_\goth a(M))=\{ \goth p\in \Ass _R(M)\vert \cd
(\goth a,R/\goth p)=n\} \ \ .$$ Dies war urspr\"unglich ein Ergebnis
von Dibaei und Yassemi ([DY, Theorem A], vgl. auch [MS, theorem
2.2]); hier wird es mit anderen Methoden bewiesen. Neue Ergebnisse
sind
\bigskip
{\bf 8.2.3 Satz}
\par
Es sei $(R,\goth m)$ ein noetherscher lokaler $d$-dimensionaler
Ring.
\par
(i) Ist $J$ ein Ideal von $R$ mit $\dim (R/J)=1$ und $\LCMo
^d_J(R)=0$, so gilt
$$\Assh (R)\subseteq \Att _R(\LCMo ^{d-1}_J(R))\ \ .$$
Ist $R$ (zus\"atzlich) komplett, so gilt sogar
$$\Att _R(\LCMo ^{d-1}_J(R))=\{ \goth p\in \Spec (R)\vert \dim
(R/\goth p)=d-1,\sqrt {\goth p+J}=\goth m\} \cup \Assh (R)\ \ .$$
(ii) F\"ur jede Folge $x_1,\dots ,x_i$ in $R$ gilt
$$\{ \goth p\in \Spec (R)\vert x_1,\dots ,x_i\hbox{ ist Teil eines
Parametersystems von }R/\goth p\} \subseteq \Att _R(\LCMo
^i_{(x_1,\dots ,x_i)R}(R))\ \ .$$ {\bf 8.2.4 Korollar}
\par
Sei $(R,\goth m)$ ein noetherscher lokaler Ring. Dann gilt f\"ur
jedes $x\in R$
$$\Att _R(\LCMo ^1_{xR}(R))=\Spec (R)\setminus {\cal V}(x)\ \ .$$
Die weitere Untersuchung zeigt, dass ein Ergebnis aus Abschnitt 8.1
(n\"amlich Satz 8.1.13) als zus\"atzliche Evidenz f\"ur unsere
Vermutung (*) aufgefasst werden kann; die Details dazu sind etwas
technisch -- vgl. Bemerkung 8.2.6 (iii) ($\zeta $). \medskip Es gibt
eine Theorie der lokalen Homologie ([T1], [T2]); diese ist in
gewisser Weise dual zur lokalen Kohomologietheorie. Seien $(R,\goth
m)$ ein noetherscher lokaler Ring, $\underline x=x_1,\dots ,x_r$
eine Folge in $\goth m$ und $X$ ein artinscher $R$-Modul. Dann ist
der $i$-te lokale Homologiemodul $\LCMo ^{\underline x}_i(X)$ von
$X$ bez\"uglich $\underline x$ definiert als
$$\vtop{\baselineskip=1pt \lineskiplimit=0pt \lineskip=1pt\hbox{lim}
\hbox{$\longleftarrow $} \hbox{$^{^{n\in \bf N}}$}} H_i(K_\bullet
(x_1^n,\dots ,x_r^n;X))\ \ ,$$ wobei $K_\bullet (x_1^n,\dots
,x_r^n;X)$ der Koszul-Komplex von $X$ bez\"uglich $x_1^n,\dots
,x_r^n$ ist und wobei $H_i$ f\"ur die $i$-te Homologie dieses
Komplexes steht. Man beachte, dass diese Homologien bez\"uglich $n$
in naheliegender Weise ein projektives System bilden. Es ist leicht
zu sehen, dass $\LCMo ^{\underline x}_I$ ein $R$-linearer
kovarianter Funktor von der Kategorie der artinschen $R$-Moduln in
die Kategorie der $R$-Moduln ist. Den Begriffen der (Krull-)
Dimension und der Tiefe (von noetherschen, also endlich erzeugten
Moduln) entsprechen hier die Begriffe der noetherschen Dimension
$\Ndim (X)$ und der Weite $\width (X)$ eines artinschen $R$-Moduls
$X$: F\"ur $X=0$ setzt man $\Ndim (X)=-1$, andernfalls bezeichnet
$\Ndim (X)$ die kleinste Zahl $r\in \Naturalsign $, zu der
$x_1,\dots ,x_r\in \goth m$ mit
$$\length (0:X(x_1,\dots ,x_r)R)<\infty $$
existieren. Eine Folge $\underline x=x_1,\dots ,x_n\in \goth m$
heisst $X$-koregul\"ar, wenn f\"ur jedes $i=1,\dots ,n$
$$(0:X(x_1,\dots ,x_{i-1})R)\buildrel x_i\over \to (0_X(x_1,\dots
,x_{i-1})R)$$ surjektiv ist. $\width (X)$ ist definiert als die
L\"ange (irgend)einer maximalen $X$-koregul\"aren Folge. [Oo] und
[Ro] sind Referenzen f\"ur diese Begriffe. Allgemein gilt
$$\width (X)\leq \Ndim (X)<\infty $$
f\"ur jeden artinschen $R$-Modul $X$; man nennt $X$
ko-Cohen-Macaulay, wenn $\width (X)=\Ndim (X)$ gilt.
\medskip
Sei $M$ ein endlich erzeugter Cohen-Macaulay $R$-Modul. Dann ist
$\LCMo ^{\dim (M)}_\goth m(M)$ ko-Cohen-Macaulay mit $\Ndim (\LCMo
^{\dim (M)}_\goth m(M))=\dim (M)$ ([T1, Proposition 2.6]). Ausserdem
gilt
$$\LCMo ^{x_1,\dots ,x_d}_{\dim (M)}(\LCMo ^{\dim
(M)}_\goth m(M))=\hat M$$ (wobei $x_1,\dots ,x_d$ ein
Parametersystem f\"ur $M$ sei). Seien nun $X$ ein artinscher
$R$-Modul mit $\Ndim (X)=d$ und $\underline x=x_1,\dots ,x_d\in
\goth m$ so, dass $(0:_X\underline x)$ endliche L\"ange hat. Tang
stellt die Frage nach der endlichen Erzeugbarkeit von $\LCMo
^{\underline x}_d(X)$ ([T1, Remark 3.5]). Wir zeigen zun\"achst mit
einem Gegenbeispiel (8.3.1), dass diese Antwort negativ zu
beantworten ist; die anschlie{\ss}ende Bemerkung 8.3.2 beantwortet
-- unter zus\"atzlichen Vorausstzungen -- die Frage dann
vollst\"andig:
\bigskip
{\bf 8.3.2 Bemerkung}
\par
Seien $(R,\goth m)$ ein noetherscher lokaler regul\"arer
$d$-dimensionaler Ring, $X$ ein artinscher ko-Cohen-Macaulay
$R$-Modul mit $\Ndim (X)=d$ und $\underline x=x_1,\dots ,x_d\in
\goth m$ so, dass $(0_X(x_1,\dots ,x_d)R)$ endliche L\"ange hat.
Dann gilt
$$\LCMo ^{\underline x}_d(X)\hbox{ ist endlich erzeut als
}R\hbox{-Modul}\iff R\hbox{ ist komplett.}$$ In einer allgemeineren
Situation gilt
\bigskip
{\bf 8.3.3 Satz}
\par
Seien $(R,\goth m)$ ein noetherscher lokaler kompletter Ring und $X$
ein artinscher $R$-Modul mit $\Ndim (X)=d$; seien $x_1,\dots ,x_d\in
\goth m$ so, dass $(0:_X(x_1,\dots ,x_d)R)$ endliche L\"ange hat.
Dann ist $\LCMo ^{\underline x}_d(X)$ als $R$-Modul endlich erzeugt.
\bigskip
Aus Satz 8.3.3 zusammen mit [T1, Remark 3.5] folgt leicht
\bigskip
{\bf 8.3.4 Korollar}
\par
Seien $(R,\goth m)$ ein noetherscher lokaler kompletter Ring und $X$
ein ko-Cohen-Macaulay $R$-Modul, $\Ndim (X)=d$; seien $x_1,\dots
,x_d\in \goth m$ so, dass $(0:X(x_1,\dots ,x_d)R)$ endliche L\"ange
hat. Dann ist $\LCMo ^{x_1,\dots ,x_d}_d(X)$ Cohen-Macaulay
(insbesondere endlich erzeugt). Im Falle $d=\dim (R)$ ist also
$\LCMo ^{x_1,\dots ,x_d}_d(X)$ ein maximaler Cohen-Macaulay
$R$-Modul.
\bigskip
Nun ordnen wir jedem endlich erzeugten $R$-Modul $M$ ordnen wir den
(artinschen) $R$-Modul
$$F_2(M)=\LCMo ^{\dim (M)}_\goth m(M)$$
und jedem artinschen $R$-Modul $X$den (endlich erzeugten, 8.3.3)
$R$-Modul
$$G_2(X)=\LCMo ^{x_1,\dots ,x_{\Ndim (X)}}(X))$$
zu. $F_2$ bzw. $G_2$ induzieren Abbildungen von der Menge der
Isomorphieklassen aller noetherschen in die Menge der
Isomorphieklassen aller artinschen Moduln bzw. umgekehrt. Auf der
anderen Seite induziert auch der Matlis-Dualit\"atsfunktor $D$
Abbildungen zwischen diesen beiden Mengen. Eine Untersuchungen der
Beziehungen zwischen diesen vier Abbildungen (Anmerkungen nach 8.3.4
und Satz 8.3.5) liefert das Ergebnis \bigskip {\bf 8.3.6 Korollar,
Aussage (ii)}
\par
Seien $(R,\goth m)$ ein noetherscher lokaler kompletter Ring und $I$
ein Ideal von $R$ mit $I\subseteq \Ann _R(M)$, $\dim (R/I)=\dim (M)$
und so, dass $R/I$ Gorenstein ist. Dann gilt
$$M\hbox { ist Cohen-Macaulay}\Longrightarrow \Hom _R(M,R/I)\hbox {
ist Cohen-Macaulay.}$$ Unterabschnitt 8.4 ist eine weitere Anwendung
des Zusammenspiels der weiter oben erw\"ahnten vier Abbildungen. Wir
definieren zun\"achst den Begriff einer Cohen-Macaulayfizierung:
\bigskip
{\bf 8.4.1 Definition}
\par
Seien $(R,\goth m)$ ein noetherscher lokaler kompletter Ring und $M$
ein endlich erzeugter $R$-Modul. Ein Obermodul $\tilde M$ von $M$
heisst Cohen-Macaulayfizierung von $M$, wenn folgende drei
Bedingungen gelten:
\par
(i) $\tilde M$ ist Cohen-Macaulay.
\par
(ii) $\dim (\tilde M)=\dim (M)$.
\par
(iii) $\dim (\tilde M/M)\leq \dim M-2$ (diese Bedingung ist
\"aquivalent zu $\LCMo ^{\dim (M)-1}_\goth m(\tilde M/M)=\LCMo
^{\dim (M)}_\goth m(\tilde M/M)=0$).
\bigskip
{\bf 8.4.2 Satz}
\par
Jede Cohen-Macaulayfizierung von $M$ (falls existent) ist zu
$(G_2\circ F_2)(M)$ isomorph.
\bigskip
In [Go] wird ein anderes Konzept des Begriffes
"Cohen-Macaulayfizierung" verwendet. Unser Begriff ist eine
Verallgemeinerung dieses Konzeptes (siehe Bemerkung 8.4.3 und Satz
8.4.5 in der vorliegenden Arbeit). Abschlie{\ss}end behandeln 8.4.6
und 8.4.7 (einfache) Beispiele von Cohen-Macaulayfizierungen. \vfil
\eject
\parindent=0pt {\bf Selbst\"andigkeitserkl\"arung}
\bigskip
\bigskip
Hiermit erkl\"are ich, die vorliegende Habilitationsschrift
selbst\"andig und ohne unzul\"assige fremde Hilfe angefertigt zu
haben. Ich habe keine anderen als die angef\"uhrten Quellen und
Hilfsmittel benutzt und s\"amtliche Textestellen, die w\"ortlich
oder sinngem\"a{\ss} aus ver\"offentlichten oder
unver\"offentlichten Schriften entnommen wurden, und alle Angaben,
die auf m\"undlichen Ausk\"unften beruhen, als solche kenntlich
gemacht. Ebenfalls sind alle von anderen Personen bereitgestellten
Materialien oder erbrachten Dienstleistungen als solche
gekennzeichnet.
\bigskip
Leipzig, den 20.05.2006
\bigskip
\bigskip
\bigskip
(Dr. Michael Hellus)

\end